\documentclass[11pt]{amsart}
\usepackage[utf8]{inputenc}
\usepackage[T1]{fontenc}
\usepackage{lmodern}
\usepackage[a4paper,margin=29mm]{geometry}
\usepackage{microtype}
\usepackage{amsmath,amsfonts,amssymb,amsthm,mathrsfs,braket,mathtools}
\usepackage{graphicx,epsfig}
\usepackage{array,makecell}
\usepackage{xcolor}
\usepackage{tikz}
\usepackage{accents}
\usepackage[authoryear,round]{natbib}
\usepackage[colorlinks=true,linkcolor=blue!55!black,citecolor=blue!55!black,urlcolor=blue!65!black]{hyperref}
\usepackage[capitalise,sort]{cleveref}

\crefname{enumi}{item}{items}
\crefname{figure}{Figure}{Figures}
\crefname{equation}{}{}
\crefname{subsection}{Subsection}{Subsections}
\makeatletter
\@ifundefined{theorem}{
  \theoremstyle{plain}
  \newtheorem{theorem}{Theorem}[section]
}{}
\@ifundefined{lemma}{
  \theoremstyle{plain}
  \newtheorem{lemma}[theorem]{Lemma}
}{}
\@ifundefined{proposition}{
  \theoremstyle{plain}
  \newtheorem{proposition}[theorem]{Proposition}
}{}
\@ifundefined{corollary}{
  \theoremstyle{plain}
  \newtheorem{corollary}[theorem]{Corollary}
}{}
\@ifundefined{assumption}{
  \theoremstyle{plain}
  \newtheorem{assumption}[theorem]{Assumption}
}{}
\@ifundefined{definition}{
  \theoremstyle{definition}
  \newtheorem{definition}[theorem]{Definition}
}{}
\@ifundefined{example}{
  \theoremstyle{definition}
  \newtheorem{example}[theorem]{Example}
}{}
\@ifundefined{xca}{
  \theoremstyle{definition}
  
}{}
\@ifundefined{remark}{
  \theoremstyle{remark}
  \newtheorem{remark}[theorem]{Remark}
}{}
\makeatother
\numberwithin{equation}{section}
\newcommand{\abs}[1]{\lvert#1\rvert}

\newcommand{\mc}[1]{\mathcal{#1}}
\newcommand{\mb}[1]{\mathbb{#1}}
\newcommand{\mf}[1]{\ifcat\noexpand#1\relax\boldsymbol{#1}\else\mathbf{#1}\fi}
\newcommand{\mr}[1]{\mathrm{#1}}

\newcommand{\ubar}[1]{\underaccent{\bar}{#1}}
\title[Spectral convergence of random feature method]{Spectral convergence of random feature method in one dimension}
\author[P. Ming]{Pingbing Ming}
\author[H. Yu]{Hao Yu}
\address{SKLMS, Institute of Computational Mathematics and Scientific/Engineering Computing, Academy of Mathematics and Systems Science, Chinese Academy of Sciences, Beijing 100190, China}
\address{School of Mathematical Sciences, University of Chinese Academy of Sciences, Beijing 100049, China}
\email{mpb@lsec.cc.ac.cn}
\email{yuhao@amss.ac.cn}
\thanks{Corresponding author: Pingbing Ming.}
\hypersetup{
  pdftitle={Spectral convergence of random feature method in one dimension},
  pdfauthor={Pingbing Ming and Hao Yu},
  pdfsubject={Random feature methods for numerical partial differential equations},
  pdfkeywords={random feature method, spectral accuracy, partition of unity method, exponential ill-conditioning}
}
\begin{document}

\begin{abstract}
We first prove the spectral convergence of the random feature method (RFM) when used to solve second-order elliptic equations and eigenvalue problems in one dimension, provided that the solutions belong to Gevrey classes or Sobolev spaces. Second, we derive the convergence rate of RFM when integrated with the Partition of Unity Method (PUM) in terms of the patch size. Finally, we show that the singular values of the resulting random feature matrix decay exponentially, leading to exponential growth of the condition number. We also prove that PUM can mitigate this excessive singular-value decay.
\end{abstract}

\maketitle

\section{Introduction}
Recently, deep learning has emerged as an active paradigm for solving partial differential equations (PDEs), leading to numerous prominent frameworks (see, e.g., \cite{EweinanYB2017DRM, raissi2019physics, sirignano2018dgm, Xu2020FiniteNeuronM, zang2020weak}). Inherently meshfree, these methods enjoy the flexibility to handle complex geometries and high-dimensional problems, facilitating the integration of experimental data. However, current deep learning-based PDE solvers still face significant challenges, such as high computational costs and a lack of strategies to consistently improve accuracy.

A parallel line of recent work has shown that learning-based and data-driven discretizations can be analyzed using classical numerical-analysis concepts such as consistency, stability, approximation, and conditioning. Rigorous error estimates have been developed for physics-informed neural networks (PINNs) and related neural residual methods~\citep{DeRyckJagtapMishra2023PINNNS, ZeinhoferMasriMardal2024PINNFramework, BonitoDeVorePetrovaSiegel2025ConsistentPINN}; least-squares formulations have been used to connect residual minimization with quasi-best approximation in finite-dimensional and machine-learning trial spaces~\citep{MonsuurSmeetsStevenson2025QLS}; sampling-based Fourier collocation and kernel methods have clarified how smoothness, sampling and approximation interact~\citep{WangBrugiapaglia2024CompressiveFourier, WenzelSantinHaasdonk2024KernelInterpolation}; and recent work on Vandermonde least-squares approximation has highlighted the central role of conditioning in high-order approximation from sampled data~\citep{ZhuNakatsukasa2025VandermondeArnoldi}. These developments point to a common numerical-analysis question: once the trial space and the sampling mechanism are specified, what are the approximation rate, the stability of the residual map, and the conditioning of the resulting algebraic system? For RFM-type PDE solvers, however, the observed spectral accuracy and the severe ill-conditioning of the random feature matrix have not yet been explained in a unified manner.

Against this background, recent efforts have sought to bridge traditional numerical solvers with machine learning, including extreme learning machines (ELMs), random feature methods, randomized neural networks, and transferable neural networks~\citep{yang2018novel,dwivedi2020physics,fabiani2021numerical,dong2021LocalELM,Chen2022BridgingTA, shang2023randomized,zhang2024transferable}. The core idea is to approximate the solution by neural networks with prescribed inner-layer weights, thereby reducing the training process to a linear least-squares problem for the outer-layer weights, which can be solved efficiently by standard linear solvers. The Partition of Unity Method~\citep{babuvska1997partition} is employed to enhance local accuracy and resolve intricate solution structures. A hallmark of these approaches is their exponential or near-exponential convergence with respect to the number of features and the refinement of the patch size~\citep{dong2021LocalELM,Chen2022BridgingTA,wang2024extreme}. This property renders these approaches particularly attractive for fields such as solid mechanics, fluid mechanics, and interface problems. Despite their diverse terminology, these methods share a common underlying framework. Although empirical spectral convergence has been observed, a rigorous theoretical proof remains lacking. In this work, we undertake a systematic theoretical analysis, using RFM as a representative case.

A critical issue in RFM is the ill-conditioning of the random feature matrix (RFMtx) within its least-squares formulation. Numerical evidence in \cite{chen2024optimization} reveals that the singular values of the RFMtx decay exponentially fast, while a rigorous theoretical account of the decay rate remains open.

\subsection{Our contributions}
We summarize the main contributions as follows.

1. \textbf{Spectral convergence with respect to the number of features.} Under the assumption on the weight distribution in Assumption~\ref{Assumption: extended sampling distribution}, we establish approximation bounds for RFM in the $W^{\sigma,p}$ norm with $\sigma\ge 0$ and $1 \le p \le \infty$. The dependence of the approximation bounds on the weight distribution, the regularity of the solution and the domain size has been clarified. Approximation bounds in Sobolev norms are crucial to derive the convergence rate when RFM is employed to solve partial differential equations. By exploiting these approximation bounds, we prove the convergence rates ranging from super-exponential to algebraic, depending on the regularity of the solution; see Table~\ref{Table: Convergence rate of random feature method proved in this paper}.

2. \textbf{Convergence rates with respect to the patch size $r$.} We demonstrate how the PUM accelerates RFM convergence by quantifying the error decay as $\mathcal{O}(r^\alpha)$; see Table~\ref{Table: Convergence rate of random feature method proved in this paper}. Specifically,  $\alpha$ scales as $N_{\mathrm{loc}} - \sigma$ for target functions in the Gevrey classes, and scales as $\min(N_{\mathrm{loc}} - \sigma, s - \sigma)$ for those in the Sobolev spaces.

3. \textbf{Singular-value decay of RFMtx.} For $N$ hidden weights in $(0,S)$ and $n$ collocation points in $[-R, R]$, we prove that the $m$-th singular value of the RFMtx for second-order PDEs decays as $\sqrt{n N} (S R)^{m-3}/(m-3)!$ and its condition number grows exponentially with $\min(n, 2N)$. This ill-conditioning is fundamentally driven by the spectral approximation among features. By establishing two-sided bounds for the singular values of the global RFMtx in terms of its sub-blocks, we show that the $m$-th singular value is closely approximated by the $\lceil m/(\text{number of patches}) \rceil$-th singular value of the sub-block RFMtx, which explains why the PUM mitigates excessive singular-value decay.
\begin{table}[htbp]
\centering
\caption{Convergence rates of RFM in $W^{\sigma,p}$ norm ($\sigma\ge 0, 1 \le p \le \infty$). $N$ ($N_{\mathrm{loc}}$) denotes the number of global (local) features, and $r$ is the maximum patch radius. Here, $\epsilon, \kappa > 0$, $0 \le \sigma \le s$.
$^*$The exponent of $r$ improves to $N_{\mathrm{loc}} + 1/p - \sigma$ for bounded hidden weights.} \label{Table: Convergence rate of random feature method proved in this paper}
\begin{tabular}{ccc}%
\hline\noalign{\smallskip}
Solution space & Rate without PUM (w.r.t. $N$) & Rate with PUM (w.r.t. $r$) \\
\noalign{\smallskip}\hline\noalign{\smallskip}
Gevrey class, $s<1$ & $\Gamma(N)^{-(1-s) + \epsilon}$ & $r^{N_{\mathrm{loc}} + 1/p - \sigma - 1} |\ln r|^{N_{\mathrm{loc}}-1}$ $^{*}$\\
Analytic class & $\exp (-\kappa N)$ & $r^{N_{\mathrm{loc}} + 1/p - \sigma - 1} |\ln r|^{N_{\mathrm{loc}}-1}$ $^{*}$\\
Gevrey class, $s>1$ & $\exp (-\kappa N^{1/s})$ & $r^{N_{\mathrm{loc}} + 1/p - \sigma - 1} |\ln r|^{N_{\mathrm{loc}}-1}$ $^{*}$ \\
Sobolev space, $s\geq \sigma$ & $N^{-s+\sigma}$ & $r^{\frac{(N_{\mathrm{loc}} - \sigma - 1 + 1/p)(s-\sigma)}{N_{\mathrm{loc}} + s - \sigma + 1/p}} |\ln r|^{N_{\mathrm{loc}}-1}$ \\
\noalign{\smallskip}\hline
\end{tabular}
\end{table}

\subsection{Related work}
\textbf{Learning-based PDE Solvers.}
Recent numerical analysis work has developed rigorous foundations for several learning-based and sampling-based PDE solvers. For the Deep Ritz method, extensive error estimates for elliptic problems may be found in~\cite{lu2022machine,JMLR:v26:24-1258}. For PINNs, error estimates have been obtained for the Navier--Stokes equations~\citep{DeRyckJagtapMishra2023PINNNS}, for broad classes of linear PDEs through abstract stability frameworks~\citep{10.1093/imanum/drab032,ZeinhoferMasriMardal2024PINNFramework}, and for elliptic problems through consistent residual losses and optimal-recovery arguments~\citep{BonitoDeVorePetrovaSiegel2025ConsistentPINN}. Least-squares and residual-based formulations have also been analyzed in settings where the trial space may be generated by machine-learning models~\citep{MonsuurSmeetsStevenson2025QLS,BachmayrDahmenOster2025NeuralResidual}. On the spectral side, compressive Fourier collocation combines random sampling with sparse recovery for high-dimensional diffusion equations~\citep{WangBrugiapaglia2024CompressiveFourier}, while kernel and RBF methods provide another meshfree approximation framework in which smoothness, scaling and sampling strongly influence convergence~\citep{LarssonSchaback2023RBFScaling, WenzelSantinHaasdonk2024KernelInterpolation}. Our analysis is complementary to these works: the hidden parameters are fixed randomly, the remaining problem is a linear least-squares system, and the main issue is to quantify both the spectral approximation power of the random feature space and the singular-value decay of the corresponding RFMtx.

\textbf{Theoretical foundations and applications of RFMs.}
Rooted in the broader paradigm of randomized architectures like ELMs~\citep{huang2006extreme,gallicchio2020deep}, RFMs possess well-established approximation capabilities. Fundamental results encompass universal approximation theorems for specific activations~\citep{huang2006universal,huang2006extreme,sun2019approximation}, uniform convergence of random Fourier features~\citep{rahimi2007random}, reproducing kernel Hilbert space (RKHS) approximation bounds~\citep{Rahimi2008UniformAO}, and inverse approximation theorems~\citep{Weinan2020TowardsAM}. Furthermore, in the infinite width limit, these networks converge to Gaussian processes~\citep{neal1996priors,williams1996computing}. RFMs have also been applied to high-dimensional PDEs~\citep{gonon2023random,gonon2023approximation}, complemented by theoretical developments~\citep{deryck2025approximation,neufeld2025full}. More broadly, they have been extended to operator learning to approximate infinite-dimensional maps~\citep{nelsen2021random,nelsen2024operator}.%

Prior work focuses mainly on high-dimensional settings and their convergence rates rarely exceed $1/N$, where $N$ denotes the number of features. In contrast, our primary interest lies in establishing the spectral accuracy of RFM approximation. Spectral accuracy, rooted in spectral methods, has also been proved for $hp$-finite element methods when the solution belongs to the Gevrey class~\citep{gui1986h1Dimension,guo1986hpversion,melenk2004hp,feischl2020exponential}. \cite{mhaskar1996neural} proved that shallow neural networks achieve a geometric approximation rate for analytic functions. See also~\cite{weinan2018exponential,montanelli2019deep,DeRyck2021app_tanh} for extensions to deep neural networks and other classes of target functions. All these proofs proceed by constructive selection of the inner and outer parameters, whereas the analysis of RFMs requires deriving convergence rates for randomly sampled features. We establish an explicit, log-term-free spectral convergence rate in the $W^{\sigma,p}$ norm for RFMs applied to functions in the Gevrey classes and the Sobolev spaces. Notably, most previous works operate within the RKHS framework and, in particular, do not impose explicit smoothness assumptions on the target function. A further key distinction of our approach is that, rather than tuning the sampling distribution to obtain an unbiased kernel estimator, our theory accommodates general prescribed distributions for the hidden weights. In particular, recent work~\citep{Fabiani2025Randomprojection} has demonstrated spectral convergence for infinitely differentiable target functions using random projection neural networks, which are essentially a version of RFM. However, that work provides neither an explicit convergence rate nor a bound on the output-weight norm.

\textbf{Condition number of RFMtx.}
The spectral behavior of RFMtxs is governed by the sampling strategy and the smoothness of the activation functions. \cite{barnett2022exponentially} has established that sampling hidden weights and collocation points on uniform grids, which corresponds to submatrices of the discrete Fourier transform, yields an exponentially ill-conditioned system. Under Gaussian sampling with a lower bound on variance products, RFMtx can be well-conditioned in strictly over- and underparameterized regimes, though it deteriorates significantly at the interpolation threshold~\citep{chen2024conditioning}. In contrast to \cite{chen2024conditioning}, our work samples weights and collocation points from bounded domains, which leads to exponential or even super-exponential singular-value decay. Recently, \cite{Zhang2025shallow} proved that Gram matrix eigenvalues decay algebraically for the ReLU activation, and exponentially for the $\tanh$ activation.
\subsection{Notation and preliminaries}
We denote by $\mb{N}$, $\mb{N}_+$, $\mb{R}$, and $\mb{R}_+$ the sets of nonnegative integers, positive integers, real numbers, and positive real numbers, respectively. The notation $a \lesssim b$ means $a \le C b$ for some constant $C$ independent of $b$. An unsubscripted $\lesssim$ implies that $C$ is universal, while subscripts indicate its dependence on specific parameters. Denote by $\Gamma$ the Gamma function. The Fourier transform of $\varphi$ in the sense of distributions is denoted by $\hat{\varphi}$. In particular, for any $\varphi \in L^1(\mb{R}^d)$, the Fourier transform $\hat{\varphi}$ is defined as
\begin{equation*}
\hat{\varphi}(\xi) = \frac{1}{(2\pi)^{d/2}} \int_{\mb{R}^d} \varphi(x) e^{-i x \cdot \xi} \, \mr{d} x.
\end{equation*}
For $s \ge 0,1 \le p \le \infty$ and a Lipschitz domain $\Omega$, the standard Sobolev space $W^{s,p}(\Omega)$ is defined as in~ \cite{Adams2003Sobolev}. We adopt the standard multi-index notation. For $\alpha = (\alpha_1, \ldots, \alpha_d) \in \mathbb{N}^d$, we denote $|\alpha| = \sum_{j=1}^d \alpha_j$, $\alpha! = \prod_{j=1}^d \alpha_j!$ and $D^\alpha = \partial_{x_1}^{\alpha_1} \cdots \partial_{x_d}^{\alpha_d}$.

As in ~\cite{rodino1993linear}, the Gevrey class $G^s(\Omega)$ is defined as
\begin{definition}
Given $s\ge0$ and $\Omega \subset \mb{R}^d$, the Gevrey class $G^{s}(\Omega)$ is the set of all $f \in C^{\infty}(\Omega)$ such that for every compact subset $K \subset \Omega$, there exist $M_{f, K}, C_{f, K}>0$ satisfying
\begin{equation*}
\max_{x \in K}\left|D^\alpha f(x)\right| \leq M_{f, K} C_{f, K}^{|\alpha|}(\alpha!)^{s}\quad \text{for all } \alpha \in \mb{N}^d.
\end{equation*}
$G^{s}_{0}(\Omega)$ is the space of all $f \in G^{s}(\Omega)$ with compact support in $\Omega$. When $\Omega$ is compact, we abbreviate $M_{f, \Omega}, C_{f, \Omega}$ to $M_{f}, C_{f}$.
\end{definition}

$G^{s}(\Omega)$ characterizes functions with $n$-th derivative growth bounded by $(n!)^{s}$. Clearly, $G^{s}(\Omega) \subseteq G^{t}(\Omega)$ for $s \leq t$, and $G^{1}(\Omega)$ is the real-analytic function space on $\Omega$. We summarize basic properties of $G^{s}(\Omega)$ in the following lemma, and defer their quantitative versions to Lemmas \ref{Lemma: G^s is a linear space and is closed to the product, dimension d} and \ref{Lemma: trigonometric functions, exponential functions, polynomials in G^0, dimension d}. The following lemma indicates that the exact solutions for some numerical examples in \cite{yang2018novel,dwivedi2020physics,dong2021LocalELM,Chen2022BridgingTA,shang2023randomized,zhang2024transferable,wang2024extreme} lie in $G^{1}(\Omega)$ or even $G^{0}(\Omega)$.
\begin{lemma} \label{Lemma: preliminary facts about G^s(Omega)}
(1) \cite[Prop.\ 1.4.5]{rodino1993linear} The space $G^{s}(\Omega)$ is a linear space, closed under pointwise multiplication and differentiation. Furthermore, $G^{1}(\Omega)$ is closed under composition. \\
(2) Trigonometric, exponential, polynomial, and bandlimited functions all belong to $G^{0}(\mb{R}^{d})$.
\end{lemma}

The remainder of this paper is organized as follows. We establish the approximation rates of the RFM in \S \ref{Section: Main results}, while \S \ref{Section: RFM solvers error estimates and conditioning} provides error estimates for RFM solvers for PDEs and analyzes the conditioning of the RFMtx. The proofs for \S \ref{Section: Main results} and \S \ref{Section: RFM solvers error estimates and conditioning} are detailed in \S \ref{Section: Fast convergence of approximation errors} and \S \ref{section: Exponential ill conditionality of random feature matrices}, respectively. We conclude this work in \S \ref{Section: Conclusion and outlook}. Proofs of auxiliary results are postponed to Appendices~\ref{Appendix section: Bounds for the Gamma function and sampling constants}--\ref{Appendix section: Auxiliary result in the proof of Theorem: decay rate of random feature matrices and their exponential ill conditionality}.
\section{Approximation rates for random feature methods} \label{Section: Main results}
\subsection{Spectral convergence with respect to the number of features} \label{subsection of main results: Spectral convergence with respect to the number of features}
We first estimate the approximation rate of RFM on the domain $\Omega=(-R,R)$ without PUM. Specifically, let $\{k_{i}\}_{i = 1}^{N}$ be i.i.d.\  samples drawn from a Borel probability distribution $\mathcal{P} \in \mathscr{P}(\mb{R})$. We employ the random Fourier feature expansion as in~\cite{rahimi2007random}
\begin{equation} \label{random feature model used in the fast convergence part}
\begin{aligned}
u_{N}(x) = \sum_{i=1}^{N} \alpha_{i} \cos \left(k_{i} x\right) + \sum_{i=1}^{N} \alpha_{N+i} \sin \left(k_{i} x\right),
\end{aligned}
\end{equation}
where $\boldsymbol{\alpha} = (\alpha_{i})_{i=1}^{2N}$ is the vector of trainable parameters. Denote the feature space
\begin{equation*}
U_N:=\operatorname{span}\{\cos(k_i\cdot),\,\sin(k_i\cdot):1\le i\le N\},
\end{equation*}
and the random feature (RF) operator $S_N: \mathbb{R}^{2N} \to U_N$ such that  $S_N\boldsymbol{\alpha} = u_{N}$ as in (\ref{random feature model used in the fast convergence part}). We make the following assumption on $\mathcal{P}$. %
\begin{assumption} \label{Assumption: extended sampling distribution}
$\mathcal{P}$ is absolutely continuous with respect to the Lebesgue measure and $\mu := \sqrt{\mb{E}|k_1|^2}<\infty$. Let $\pi$ be the probability density function of $|k_1|/\mu$.
Assume $$\varpi := \max\left\{1,\sup_{x>0} \max(1,x^3)\pi(x)\right\}<\infty.$$
\end{assumption}
Observe that $\varpi$ is scale invariant and $\varpi\ge 1$. Some commonly used distributions satisfy Assumption~\ref{Assumption: extended sampling distribution} with moderately large $\varpi$.
\begin{example} \label{Example: constants mu and varpi for normal, exponential and Student's t distributions}
(1) Uniform distribution: If $k_1\sim\operatorname{Unif}(-S,S)$,
then $\mu = S/\sqrt{3}$ and $\varpi = 3.$

(2) Normal distribution: If $k_1\sim \mc{N}(0,\sigma^2)$,  then $\mu = \sigma$ and  $\varpi = 1.$

(3) Exponential distribution: %
If $k_1 \sim \operatorname{Exp}(\theta)$, then $\mu = \sqrt{2}/\theta$, and $\varpi = \sqrt{2}.$

(4) Student's t distribution with $\nu > 2$ degrees of freedom: If $k_1\sim t_{\nu}$,
then
$$\mu=\sqrt{\frac{\nu}{\nu-2}}, \ \ \text{and}\ \
\varpi = \max\left\{1,\frac{2\Gamma((\nu+1)/2)}{\sqrt{(\nu-2)\pi}\Gamma(\nu/2)}
\max\left(1,  3\sqrt{3} \left(\frac{\nu+1}{\nu-2} \right)^{-(\nu+1)/2} \right)\right\}.$$
For Student's $t$ distribution, Lemma \ref{Lemma: bound for varpi when k_1 follows a Student's t distribution} gives $\varpi<1.28$ when $\nu\ge 3$ and $\varpi=1$ when $\nu\ge 4.69$.
\end{example}

We shall construct in \S\ref{Section: Fast convergence of approximation errors} a random linear operator $\mathcal{R}_N[\{k_{i}\}_{i = 1}^{N},\mu,C,c]: C^\infty(\Omega)\to U_{N}$. The following theorem gives an approximation bound, upon which we derive the convergence rate of RFM for different scenarios.
\begin{theorem} \label{Theorem: main theorem of RFM, Gevrey classes, W^(l,p) error norm, extended distribution}
Assume  $u \in G^{s}(\bar{\Omega})$ with $s \leq 1$ and Assumption \ref{Assumption: extended sampling distribution} holds. Let $c>2$ and $u_{N} = \mathcal{R}_N[\{k_{i}\}_{i = 1}^{N},\mu,C_u,c] u$ be the RF approximation for $u$. Let $l\in \mb{N}$ with $l\le 2N$. For $l<2N$, define
$$\tau %
:= \sqrt{5}\varpi \max(C_{u}, \mu) R \Gamma(2N-l+1)^{(s-1)/(2N-l)}, $$
and set $\tau:=\sqrt{5}\varpi \max(C_{u}, \mu) R$ for $l=2N$. Then, for any $1\le p\le\infty$,
\begin{equation*}
\begin{aligned}
\mb{E} \abs{u-u_{N}}_{W^{l,p}(\bar{\Omega})}
& \lesssim_{l}  M_{u} \varpi^{l}\max(C_{u}, \mu)^{l} R^{1/p} N\\
& \quad \times\left(e^{-[c-1-\ln(1+c)](N-1)}\mf{1}_{\{N\ge2\}}   + N^{(l-2)s-1/p}e^{c(N-1)}(1+\tau) \tau^{2N-l}\right). \\
\end{aligned}
\end{equation*}
\end{theorem}

By absorbing constants independent of $N$, the above error bound changes to
\begin{equation}
\begin{aligned}  \label{error bound in Theorem:Theorem: main theorem of RFM, Gevrey classes, W^(l,p) error norm}
\mb{E}\abs{u-u_{N}}_{W^{l,p}(\overline{\Omega})}  %
& \lesssim N e^{-[c-1-\ln(1+c)](N-1)} + N^{(l-2)s-1/p+1} e^{c(N-1)} \tau^{2N-l}  .
\end{aligned}
\end{equation}

\begin{remark}[Sufficient number of features]
By Stirling's formula, $\tau \ll 1$ when
$$
2N -l\gg_{s,\varpi} \bigl(   R\max(C_{u}, \mu)\bigr)^{1/(1-s)}.
$$ %
Thus, $(R\max(C_{u}, \mu))^{1/(1-s)}$ reflects the difficulty of approximation and gives a sufficient scale for the number of features. For example, when $s = 0$, it indicates that the number of features is proportional to $\max(C_{u}, \mu)$, i.e., sampling in the frequency space reaches a certain threshold of ``density''.
\end{remark}
\begin{remark}[Optimal choice of \ensuremath{\mu}] %
It is natural to choose \ensuremath{\mu} comparable to, or slightly larger than \ensuremath{C_{u}}. Taking \ensuremath{\mu} too large unnecessarily enlarges \ensuremath{\max(C_{u},\mu)} and introduces excessively high-frequency features. Conversely, if \ensuremath{\mu} is too small, the high-frequency components of $u$ may not be sufficiently captured, and RFMtx is severely ill-conditioned; see~\S\ref{subsection of main result: Exponential ill conditionality of random feature matrices}. Numerical results in Tables A.4 and A.5 of \cite{Chen2022BridgingTA} support this prediction.
\end{remark}

In what follows, we shall prove the spectral accuracy of RFM. The smoother target functions admit faster convergence rates, ranging from algebraic to super-exponential. In Corollary~\ref{Corollary: spectral convergence when solution u in Analytic class, Gevrey class s>1 or Sobolev space, W^(l,p) error norm}, $\mc{T}_{\mu}$ denotes the case-dependent linear operator constructed in the proof, which maps $u$ into $G^{0}(\bar{\Omega})$ with Gevrey constant $C_{\mc{T}_{\mu}u} \le \mu$.
\begin{corollary}[Gevrey class $s<1$] \label{Corollary: Gevrey classes s<1, super exponential convergence for N, W^(l,p) error norm}
Let $\sigma\ge0$. Under the assumptions of Theorem \ref{Theorem: main theorem of RFM, Gevrey classes, W^(l,p) error norm, extended distribution} with $s < 1$, for any $\epsilon > 0$, the following estimate holds with an implied constant independent of $N$:
\begin{equation*}
\begin{aligned}
\mb{E} \|u-u_{N}\|_{W^{\sigma,p}(\bar{\Omega})}
& \lesssim \Gamma(N)^{-(1-s) + \epsilon} .
\end{aligned}
\end{equation*}
\end{corollary}

\begin{corollary}  \label{Corollary: spectral convergence when solution u in Analytic class, Gevrey class s>1 or Sobolev space, W^(l,p) error norm}
Under Assumption \ref{Assumption: extended sampling distribution}, let $\mu = N/({34 \varpi R})$ and $u_{N} = \mc{R}_{N}\mc{T}_{\mu}u$.  %
Let $\sigma\ge0$ and $1\le p\le \infty$. Denote $\mc{M}(t) = 1+t^{\sigma}$, $\mc{S}_{R}=R^{1/p}(R^{-\sigma}+1) $ and $\zeta_s = \sum_{k=1}^\infty k^{-s}$. \\
(a) \textbf{Analytic class.} If there exists  $\rho>0$  such that  $u$  can be extended to an analytic function in the complex strip  $|\operatorname{Im} z| < \rho$  with  $\|u(\cdot+i y)\|_{L^{2}(\mb{R})}$  uniformly bounded for all  $y \in(-\rho, \rho)$, then with $\kappa = \min(\rho/(68\varpi R),1/2)$, %
\begin{equation*}
\begin{aligned}
& \mb{E}\|u-u_{N}\|_{W^{\sigma,p}(\Omega)}
\lesssim_{\sigma, p} \mc{M}(R/\rho) \mc{S}_{R} \rho^{-1/2} \|e^{\rho |\cdot|} \hat{u}\|_{L^{2}(\mb{R})} e^{-\kappa N}  .
\end{aligned}
\end{equation*}
(b) \textbf{Gevrey class $s>1$.}  If there exists  $0<\rho\le R$  such that  $u$ has an extension in $G^{s}(-R-\rho,R+\rho)$ with $s>1$, then with $\kappa = (3s/8) \left[34 \varpi R(C_u + 2 \zeta_s /\rho)\right]^{-1/s}$, %
\begin{equation*}
\begin{aligned}
\mb{E}\|u-u_{N}\|_{W^{\sigma,p}(\Omega)} &  \lesssim_{s, \sigma, p} M_{u} \mc{M}(C_u R + R/\rho) (C_u R + R/\rho) \mc{S}_{R} e^{-\kappa N^{1/s}}  .
\end{aligned}
\end{equation*}
(c) \textbf{Sobolev space.} If $u \in W^{s,p}(\Omega)$ and $0\le \sigma\le s$, then
\begin{equation*}
\begin{aligned}
\mb{E}\|u-u_{N}\|_{W^{\sigma,p}(\Omega)}
& \lesssim_{s,\sigma,p} (R^s+R^{-\sigma}) \varpi^{s-\sigma} \|u\|_{W^{s,p}(\Omega)} N^{-s+\sigma} .
\end{aligned}
\end{equation*}
\end{corollary}

\subsection{Convergence rates with respect to the patch size} \label{subsection in main result: The combination of partition of unity}

The solution of the PDE typically exhibits local variations, possibly at small scales. To accommodate this, following~\cite{Chen2022BridgingTA}, we construct local random feature approximants and assemble them by PUM. We consider a general class of partition of unity (PoU) functions, and make the following standard assumptions on the PoU functions and their corresponding patches.
\begin{assumption} \label{Assumption: PoU with Finite Overlapping Condition and Derivative Bounds}
Let $\{\phi_j\}_{j=1}^J \subset W^{m, \infty}(\Omega)$ form a partition of unity. Let the patches $\Omega_j := \left(\operatorname{supp}\phi_j\right)^{\circ}$ be open intervals with size $r_j := \operatorname{diam}(\Omega_j)/2$. We assume

1. {(Partition of Unity)} $\sum_{j=1}^J \phi_j(x) \equiv 1$ on $\Omega$.

2. {(Finite Overlapping Condition)} There exists a constant $M \in \mathbb{N}$ such that
    \begin{equation*}
    \sup_{x\in\Omega}\operatorname{card}\{j \mid x \in \Omega_j\} \le M.
    \end{equation*}

3. {(Derivative Bounds)} There exists a constant $C_{\phi} > 0$, independent of $j$, such that
    \begin{equation} \label{eq:pou_deriv_bound}
        \|D^k \phi_j\|_{L^\infty(\Omega_j)} \le C_{\phi}r_j^{-k}, \quad \text{for all integers } 0 \le k \le m.
    \end{equation}
\end{assumption}
On each patch $\Omega_j$, %
the local trial function $v_j$ is taken as the local random feature expansion
\begin{equation*}
v_{j}(x) = \sum_{q=1}^{N_{j}} \left[ \alpha_{j,q} \cos\left(k_{j,q} (x - x_{j})\right) + \alpha_{j,q+N_{j}}\sin\left(k_{j,q} (x - x_{j})\right) \right],
\end{equation*}
where $x_j$ is the center of $\Omega_j$, and $\{k_{j,q}\}$ are i.i.d.\  hidden weights satisfying Assumption \ref{Assumption: extended sampling distribution} with parameters $\mu_j$ and $\varpi$. Assembling the local approximants $v_j$ via the PoU functions, we get the global approximation $u_{RF}$ as
\begin{equation} \label{random feature model with partition of unity}
u_{RF}(x) = \sum_{j=1}^J \phi_j(x) v_j(x).
\end{equation}
In this part, we prove global approximation estimates for RFM with PUM. On each patch, we construct $v_j$ via the randomized linear approximation in Theorem \ref{Theorem: main theorem of RFM, Gevrey classes, W^(l,p) error norm, extended distribution}. Since the dependence of the local errors on the number of features and other parameters has been detailed in Theorem \ref{Theorem: main theorem of RFM, Gevrey classes, W^(l,p) error norm, extended distribution}, we mainly focus on the convergence rate with respect to the patch size $r_j$. To state the result, we define the local scaling factor $A_j = \varpi \max(C_{u,\bar{\Omega}_j}, \mu_j) r_j$.
\begin{theorem} \label{Theorem: error bound for RFM combined with PUM, Gevrey class s leq 1}
Assume $u \in G^{s}(\bar{\Omega})$ and $A_j \le e^{-2}$ for all $j$. Then, for any $m \in \mathbb{N}$ and $1 \le p \le \infty$ such that $m \le \min_j N_j-1$, there holds%
\begin{equation*}
\mb{E}\|u-u_{RF}\|_{W^{m,p}(\Omega)} %
\lesssim  M_u \max_{1\le j\le J} \left( r_j^{{1}/{p}-m-1} A_j^{N_j} |\ln A_j|^{N_j-1} \right),
\end{equation*}
where the implied constant depends only on $m$, $p$, $s$, $C_{\phi}$, $M$, $R$ and $\mathbf{N} = (N_1, \dots, N_J)$.
\end{theorem}
The factor $A_j$ motivates adaptive refinement, allowing smaller $r_j$ where $C_{u,\bar{\Omega}_j}$ is large to equilibrate and reduce the error.  Dropping the constant $\max_j(C_{u,\bar{\Omega}}, \mu_j)$, Theorem \ref{Theorem: error bound for RFM combined with PUM, Gevrey class s leq 1} reveals an $\mathcal{O}(r^{N_{\mathrm{loc}}+{1}/{p}-m-1}|\ln r|^{N_{\mathrm{loc}}-1})$ convergence rate.  %
Furthermore, if $k_{j,q}$ are bounded random variables, a common scenario in practice, this rate improves by one order:
\begin{corollary}[Gevrey class] \label{Corollary: error bound for RFM combined with PUM, Gevrey class, bounded hidden weights}
Under the assumptions of Theorem \ref{Theorem: error bound for RFM combined with PUM, Gevrey class s leq 1}, if $|k_{j,q}|\le S$ a.s. for some $S\ge 0$ and all $j,q$, then for any $m \in \mathbb{N}$ and $1 \le p \le \infty$ such that $m \le \min_j N_j$, %
\begin{equation*}
\mb{E}\|u-u_{RF}\|_{W^{m,p}(\Omega)}
\lesssim \max_{1\le j\le J} \left( r_j^{N_j+{1}/{p}-m} |\ln r_j|^{N_j-1} \right),
\end{equation*}
where the implied constant does not depend on $r_j$.
\end{corollary}

As a direct consequence, we obtain the convergence rate in Sobolev spaces. For simplicity, we assume $N_j = N_{\mathrm{loc}}$ for all patches.
\begin{corollary}[Sobolev space]  \label{Corollary: convergence for PUM when solution u in Analytic class, Gevrey class s>1 or Sobolev space}
Let $u \in W^{s,p}(\Omega)$ with $s \ge 0$ and $1\le p\le\infty$. Assume %
the PoU is quasi-uniform, i.e. $r := \max_{1 \le j \le J} r_j \lesssim \min_{1 \le j \le J} r_j $. Suppose $\mu_j = \nu r^{t}$ with $\nu > 0$ and $t \ge (s+1)/(N_{\mathrm{loc}}+s-\sigma+1/p)-1$. Then, for sufficiently small $r$ and any $0 \le \sigma \le \min(s,N_{\mathrm{loc}}-1)$,
\begin{equation} \label{eq in Corollary: convergence for PUM when solution u in Sobolev space}
\begin{aligned}
\mb{E}\|u-u_{RF}\|_{W^{\sigma,p}(\Omega)} & \lesssim r^{\frac{(N_{\mathrm{loc}} - \sigma - 1 + 1/p)(s-\sigma)}{N_{\mathrm{loc}} + s - \sigma + 1/p}} |\ln r|^{N_{\mathrm{loc}}-1} \|u\|_{W^{s,p}(\Omega)} ,
\end{aligned}
\end{equation}
where the implied constant does not depend on $r$.
\end{corollary}

For sufficiently large $N_{\mathrm{loc}}$, the rate in (\ref{eq in Corollary: convergence for PUM when solution u in Sobolev space}) approaches $s-\sigma$; for sufficiently large $s$, it tends to $N_{\mathrm{loc}} - \sigma - 1 + 1/p$, which coincides with the rate in Theorem~\ref{Theorem: error bound for RFM combined with PUM, Gevrey class s leq 1}. Since $N_{\mathrm{loc}} - \sigma \ge 1$, the corollary always holds for $t=0$. This practically significant setting corresponds to employing a fixed sampling law for the local hidden weights during the patch size refinement~\citep{Chen2022BridgingTA,sun2024local}.

\section{RFM solvers: error estimates and conditioning} \label{Section: RFM solvers error estimates and conditioning}
\subsection{Abstract error estimates for RFM solvers}\label{subsection: Abstract Error Estimates for RFM Solvers}
This part establishes abstract error estimates for RF discretization. These bounds are independent of RF realization, relying solely on the stability of the underlying PDEs and the approximation capacity of the trial space. Applying these abstract estimates, taking the expectation, and invoking approximation results from the previous section, we obtain the final error estimates.
\subsubsection{Strong-form RFM}
Consider the two-point boundary value problem for a general second-order linear differential operator
\begin{equation}
\left\{
\begin{aligned}\label{1-d elliptic equation with variable coefficients, Dirichlet BC}
Lu &:= a u^{\prime \prime} + b u^{\prime}+ c u = f, \quad &&x \in \Omega :=(-R,R), \\
Bu &:= g_{1} u^{\prime} + g_{2} u = g,   \quad &&x \in \partial\Omega,
\end{aligned}\right.
\end{equation}
where $a\in L^{\infty}(\Omega)$, $b, c, f \in L^{2}(\Omega)$ and $g_{1}, g_{2}, g \in L^{2}(\partial \Omega).$
$\Lambda_{1}$ denotes the $L^{\infty}(\Omega)$-norm of $a$, $\Lambda_{2}, \Lambda_{3}$ denote the $L^{2}(\Omega)$-norms of $b, c$, and $\Lambda_{4}, \Lambda_{5}$ denote the $L^{2}(\partial\Omega)$-norms of $g_1, g_2$.
The operator $B$ accommodates general boundary conditions (BCs). If Dirichlet BC is imposed on $\partial \Omega_{D} \subset \partial \Omega$ and Neumann or Robin BC is imposed on $\partial \Omega_{N} = \partial \Omega \setminus \partial \Omega_{D}$, then $g_{1} = \tilde{g}_{1}\mathbf{1}_{\partial \Omega_{N}}$ and $g_{2} = \mathbf{1}_{\partial \Omega_{D}} + \tilde{g}_{2}\mathbf{1}_{\partial \Omega_{N}}$.

To solve Problem~\eqref{1-d elliptic equation with variable coefficients, Dirichlet BC}, we employ a random linear trial space $V_{RF}$ spanned by the ansatz functions in (\ref{random feature model used in the fast convergence part}) or (\ref{random feature model with partition of unity}). The loss function is defined by
\begin{equation}
\begin{aligned} \label{loss function for solving 1D problems}
\mathcal{L}(u)= & \left\|L u-f\right\|_{L^{2}(\Omega)}^{2} + \gamma \left\|B u-g\right\|_{L^{2}(\partial \Omega)}^{2} ,
\end{aligned}
\end{equation}
where $\gamma > 0$ is the boundary penalty parameter. The RF solution is then obtained by %
\begin{equation*}
u_{RF}^* = \mathop{\arg\min}_{u \in V_{RF}} \mathcal{L}(u).
\end{equation*}
\begin{theorem}[Strong-form RFM]%
\label{Theorem: abstract error estimate for Strong-form RFM}
Assume that $1/a\in L^\infty(\Omega)$. Let $u^*$ be a solution to Problem (\ref{1-d elliptic equation with variable coefficients, Dirichlet BC}). %
Then, there exists a constant $C$ independent of $u^*$ and $u_{RF}^*$ such that
\begin{equation*}
\inf_{v \in V} \|u^{\ast}+v-u_{RF}^{\ast}\|_{H^{2}(\Omega)} \leq C \inf_{w \in V_{RF}} \|u^{\ast}-w\|_{H^{2}(\Omega)} ,
\end{equation*}
where the null space $V:=\{v\in H^2(\Omega): Lv=0 \text{ in }\Omega,\ Bv=0
\text{ on }\partial\Omega\}$. In particular, if the solution is unique, then $\|u_{RF}^* - u^*\|_{H^{2}(\Omega)} \leq C \inf_{w \in V_{RF}} \|w - u^*\|_{H^{2}(\Omega)}.$
\end{theorem}
This result accommodates various BCs and problems with a non-trivial null space, such as the one-dimensional Laplace equation with Neumann BC.
\subsubsection{Weak-form RFM}
Let $V$ be a Hilbert space with norm $\|\cdot\|_V$. Consider the variational problem: %
find $u^* \in V$ such that
\begin{equation}\label{eq: weak variational problem}
a(u^*,v)=F(v),\qquad \forall v\in V,
\end{equation}
where $a:V\times V\to\mathbb R$ is a bilinear form and $F:V\to\mathbb R$ is a bounded linear functional. Assume that $a$ is continuous and coercive, i.e., there exist constants $M,\alpha>0$ such that
\begin{equation}\label{continuity and coercivity of a}
|a(w,v)|\le M\|w\|_V\|v\|_V\qquad a(v,v)\ge \alpha\|v\|_V^2,\qquad \forall w,v\in V.
\end{equation}
By the Lax--Milgram theorem, Problem \eqref{eq: weak variational problem} admits a unique solution. Let $V_{RF}\subset V$ be the RF trial space. The weak-form RF solution $u_{RF}\in V_{RF}$ is defined by
\begin{equation*}\label{weak form RFM}
a(u_{RF},v)=F(v)\qquad \forall v\in V_{RF}.
\end{equation*}
A direct application of C\'ea's lemma~\citep{Cea1964Approximation} gives %
\begin{theorem}[Weak-form RFM]\label{Theorem: abstract error estimate for Weak-form RFM}
Under Assumption \eqref{continuity and coercivity of a}, it holds that
\begin{equation*}\label{eq:weak form RFM error estimate}
\|u^*-u_{RF}\|_V\le \frac{M}{\alpha} \inf_{w\in V_{RF}}\|u^*-w\|_V.
\end{equation*}
\end{theorem}

\subsubsection{RFM for Eigenvalue Problems}
We now consider the RFM approximation for elliptic eigenvalue problems. Let $a(\cdot,\cdot)$ be a symmetric, continuous and coercive bilinear form as in the weak formulation, and let $m(\cdot,\cdot)$ be a symmetric, continuous and nonnegative bilinear form on $V\times V$. Assume that the solution operator $T:V\to V$ defined by $a(Tf,v)=m(f,v)$ for all $v\in V$ is compact. The continuous eigenvalue problem is to find $(\lambda,u)\in\mathbb R\times V$, $u\ne0$, such that
\begin{equation*}\label{continuous eigenvalue problem}
a(u,v)=\lambda m(u,v),\qquad \forall v\in V.
\end{equation*}
Let $\lambda$ be an eigenvalue of the above problem with multiplicity $q$, and let $E$ be the corresponding eigenspace. Define the best approximation error of $E$ by
\begin{equation*}\label{RFM eigenspace best approximation error}
\eta_{RF}(E):=\sup_{u\in E, \ \|u\|_V=1}\inf_{w\in V_{RF}}\|u-w\|_V.
\end{equation*}
The RF eigenvalue approximation is to find $(\lambda_{RF},u_{RF})\in\mathbb R\times V_{RF}$, $u_{RF}\ne0$, such that
\begin{equation}\label{RFM eigenvalue problem}
a(u_{RF},v)=\lambda_{RF}m(u_{RF},v),\qquad \forall v\in V_{RF}.
\end{equation}
The following theorem is a direct consequence of \citet[Theorems 9.12 and 9.13]{Boffi2010FEAeigenvalue}.

\begin{theorem}[RFM for eigenvalue problems]\label{Theorem: abstract error estimate for RFM eigenvalue problems}
Under the above assumptions, for $V_{RF}$ sufficiently rich, the discrete problem \eqref{RFM eigenvalue problem} has exactly $q$ eigenvalues $\{\lambda_{RF,j}\}_{j=1}^{q}$ converging to $\lambda$, counted with multiplicity. Let $E_{RF}$ be the space spanned by the corresponding discrete eigenfunctions. Then, there exists a constant $C>0$, independent of $V_{RF}$, such that
\begin{equation*}
\begin{aligned}
& \sup_{u\in E,\  \|u\|_V=1}\inf_{v\in E_{RF}}\|u-v\|_V\le C\eta_{RF}(E), \\
& \max_{1\le j\le q}|\lambda-\lambda_{RF,j}|\le C\eta_{RF}^2(E).
\end{aligned}
\end{equation*}
\end{theorem}
Theorem \ref{Theorem: abstract error estimate for RFM eigenvalue problems} indicates that the convergence rate of the eigenspace matches the best approximation error, whereas the eigenvalues converge at twice this rate.

\subsection{Fast singular-value decay and exponential ill-conditioning of RFMtx} \label{subsection of main result: Exponential ill conditionality of random feature matrices}
\begin{figure}%
    \centering
    \includegraphics[width=0.38\textwidth]{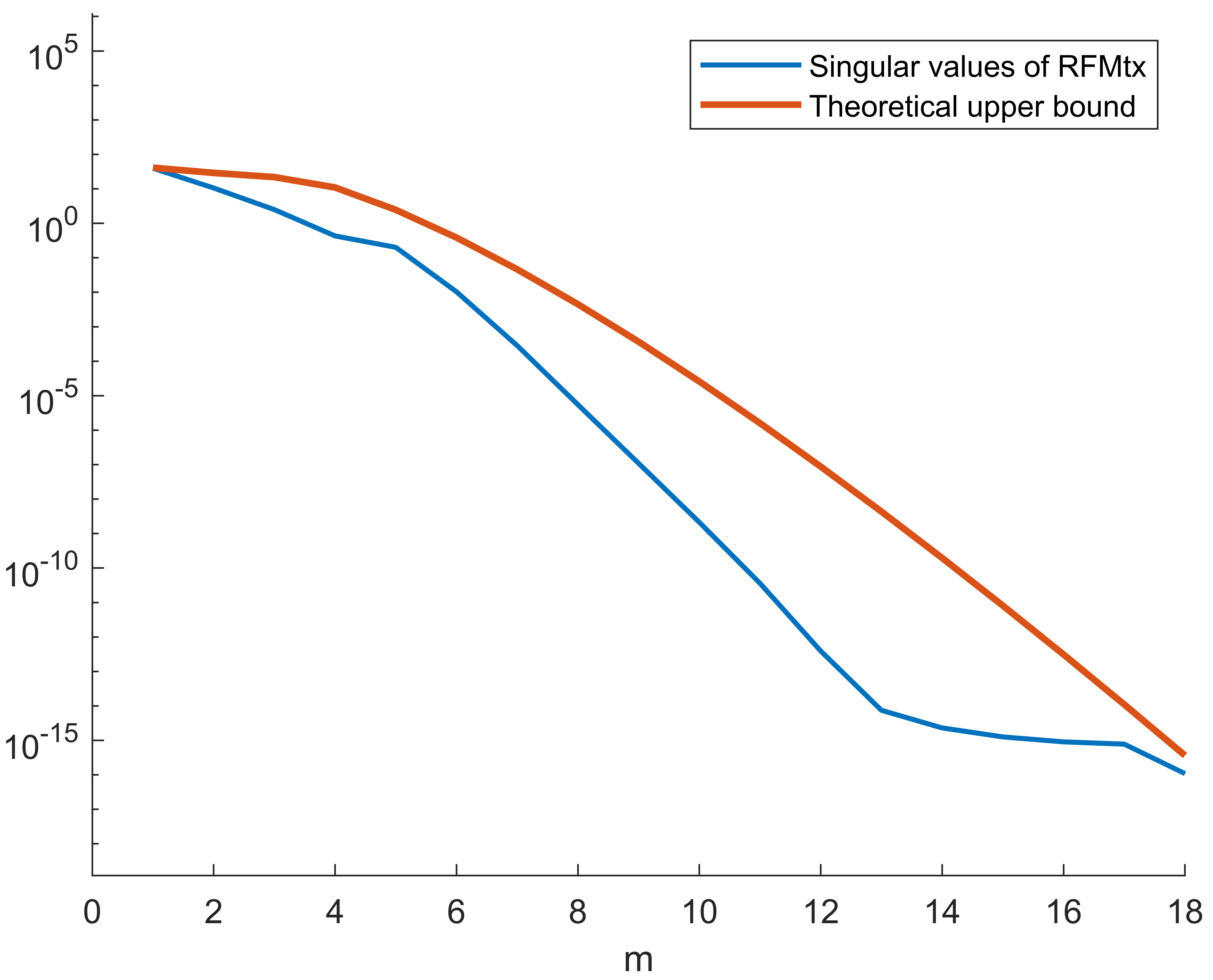}
    \includegraphics[width=0.108\textwidth]{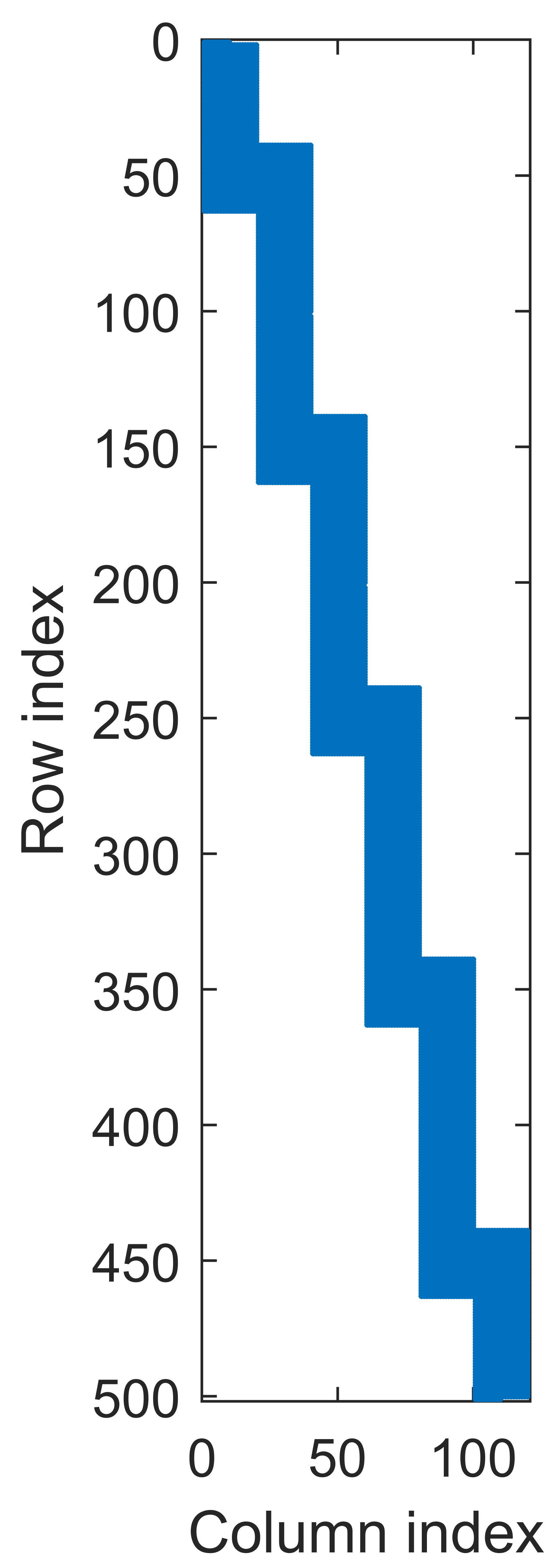}
    \includegraphics[width=0.38\textwidth]{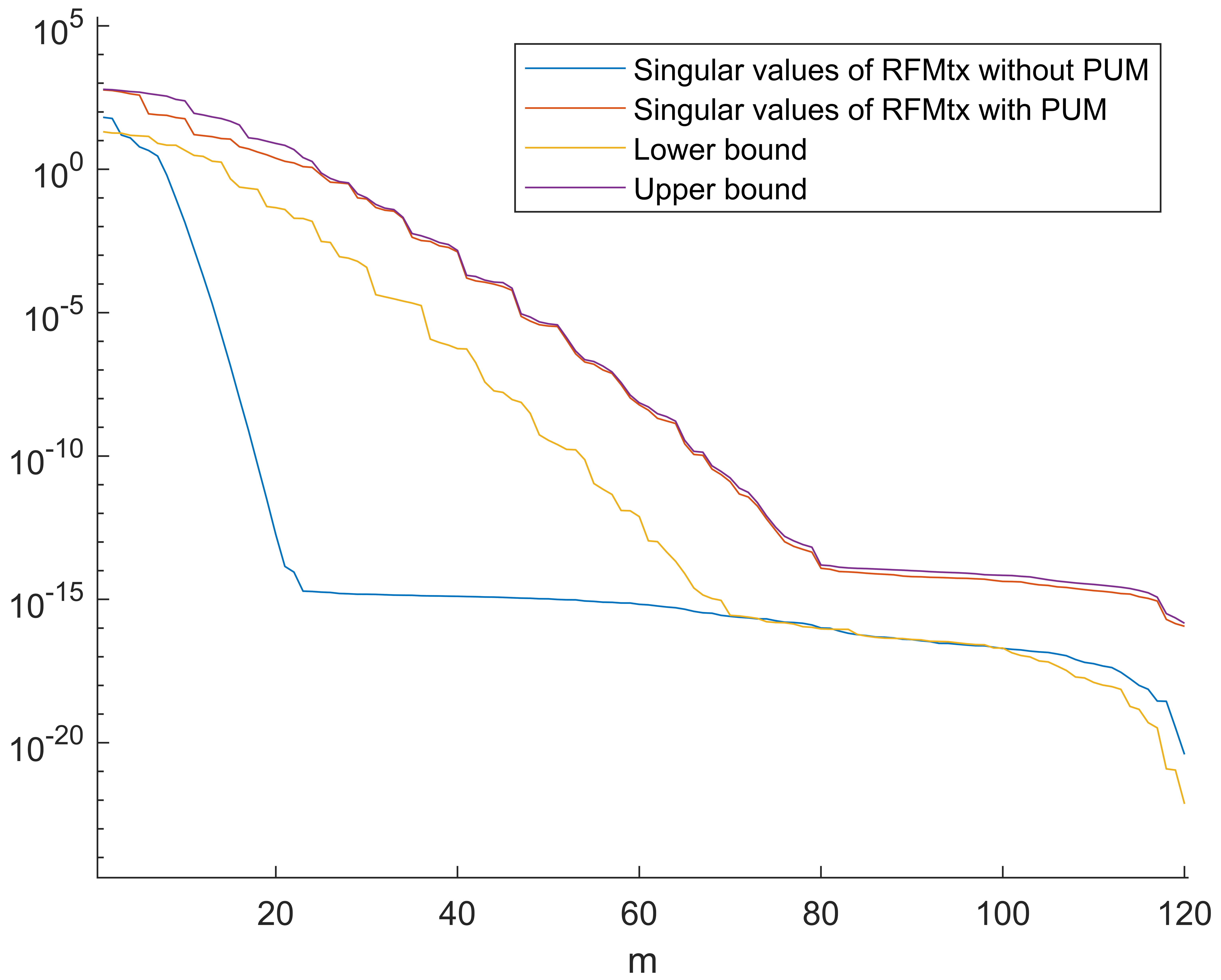}
    \caption{Singular values and sparsity pattern of the RFMtx. \textbf{Left:} Singular values of an RFMtx without PUM ($18$ features, $S = 1$, $R = 1/2$) compared against the upper bound in (\ref{ineq: fast decay of sigma_m(Z)}). \textbf{Middle:} Sparsity pattern of an RFMtx using $501$ equidistant collocation points, $6$ patches, and $20$ local features per patch ($N_{\mathrm{loc}} = 10$). \textbf{Right:} Singular values of RFMtxs with and without PUM, alongside theoretical bounds given by (\ref{ineq: basic lower and upper bounds for RFMtx with PUM}). Here, $S = 1$, $R = 4$, and all other settings match those in the middle panel. The operator $\widetilde{L}$ is the second derivative with Dirichlet BC.} \label{Figure: Sparsity pattern of RFMtx}
\end{figure}

In practical computations, the collocation points $\{x_{i}\}_{i = 1}^{n-2} \subset (-R, R)$, $x_{n-1}=-R$, and $x_{n}=R$ are selected to discretize $\mc{L}(u)$. %
We first study RFM without PUM. Define the operator $\widetilde{L}$ by $\widetilde{L} = L$ in $\Omega$ and $\widetilde{L} = \sqrt{\gamma} B$ on $\partial\Omega$. The trainable parameter vector $\boldsymbol{\alpha}$ is obtained by solving the least-squares problem $\mf{\Phi}\boldsymbol{\alpha} = \mf{F}$, where $\mf{\Phi} = \left(\Phi_{i j}\right) \in \mb{R}^{n \times 2 N}$ is the RFMtx with
\begin{equation*}
\Phi_{i j}=\left\{\begin{aligned}
& \widetilde{L}\cos \left(k_{j} x_{i}\right), \quad &&1 \leq j \leq N, \\
& \widetilde{L}\sin \left(k_{j-N} x_{i}\right), \quad &&N+1 \leq j \leq 2 N,
\end{aligned}\right.
\end{equation*}
and $\mf{F} = (f(x_{1}), \ldots, f(x_{n-2}), \sqrt{\gamma}g(x_{n-1}), \sqrt{\gamma}g(x_{n}))^{\top}$. Each column of $\mf{\Phi}$ corresponds to a feature and each row corresponds to a collocation point. Let $\{\sigma_m(\mf{\Phi})\}_{m = 1}^{\min(n,2N)}$ be the singular values of $\mf{\Phi}$ in a descending order. Denote constants $A,C,D>0$ and $B\in\mb{R}$ such that $A^{2} := \sum_{i=1}^{n-2} a^{2}(x_{i})$, $B := \sum_{i=1}^{n-2} a(x_{i}) c(x_{i})$, $C^2 :=  \sum_{i=1}^{n-2} c^{2}(x_{i})$ and $D^{2} =\sum_{i=1}^{n-2}b^{2}(x_{i})$. Let $\tilde{\Lambda}_{4}=\sqrt{\gamma}\Lambda_{4}$ and $\tilde{\Lambda}_{5}=\sqrt{\gamma}\Lambda_{5}$. The following theorem gives the decay rate of $\sigma_{m}$ and reveals the exponentially ill-conditioning of $\mf{\Phi}$. %
\begin{theorem} \label{Theorem: decay rate of random feature matrices and their exponential ill conditionality}
Assume that the hidden weights $\{k_{j}\}_{j = 1}^{N} \subset (0, S)$ are pairwise different\footnote{Otherwise, an orthogonal transformation can be used to leave the same columns in the matrix $\mf{\Phi}$ with a constant multiple of only one of them, and the others eliminated to zero.}. Let $\hat{\tau}(3) := S R$ and $\hat{\tau}(m) := S R \Gamma(m-2)^{-{1}/({m-3})}$ for $m > 3$. \\
(1) \textbf{Singular values decay.} For $m \le2$, $\sigma_{m}(\mf{\Phi}) \leq \|\mf{\Phi}\|_{F} /\sqrt{m}$.
For  $m \geq 3$, let $\hat{\tau} = \hat{\tau}(m).$ Then,
\begin{equation}
\begin{aligned} \label{ineq: fast decay of sigma_m(Z)}
\sigma_{m}(\mf{\Phi}) & \leq  \sqrt{N} \left(\sqrt{3}A S^{2} + \sqrt{3D^{2} + 2\tilde{\Lambda}_{4}^{2}} S \hat{\tau} + \sqrt{3C^{2} + 2\tilde{\Lambda}_{5}^{2}} \hat{\tau}^{2}\right) (1+ \hat{\tau})\hat{\tau}^{m-3}.
\end{aligned}
\end{equation}
(2) \textbf{Lower bounds for condition number.}
Set $M:=\min(n,2N)$ and assume $M\ge3$. Let $\hat{\tau}:=\hat{\tau}(M)$. Then,\footnote{When $\sigma_{M}(\mf{\Phi}) = 0$ (i.e., $\mf{\Phi}$ is rank-deficient), we define $\kappa(\mf{\Phi}) := +\infty$.}
\begin{equation*}
\begin{aligned}
\kappa(\mf{\Phi}) := \frac{\sigma_{1}(\mf{\Phi})}{\sigma_{M}(\mf{\Phi})}
& %
\geq\sqrt{\dfrac{\rho}{6M}}
\dfrac{\hat{\tau}^{3-M}}{\max(1,\hat{\tau}^3)},
\end{aligned}
\end{equation*}
where $\rho \geq 1 - 2\max(B, 0)\ubar{S}^{2}\left[{A^{2} \ubar{S}^{4} + (D^{2} + \tilde{\Lambda}_{4}^{2})\ubar{S}^{2} + C^{2} + \tilde{\Lambda}_{5}^{2}}\right]^{-1}$ with $\ubar{S} = \min \Big(S, A^{-1/2}  \sqrt[4]{C^{2} + \tilde{\Lambda}_{5}^{2}}\Big).$

Moreover, if $\{k_{j}\}_{j = 1}^{N} \overset{i.i.d.}{\sim}  \operatorname{Unif}(0, S)$, then $\mb{P}(\rho \geq 1/200)\ge 1-0.95^N$. %
\end{theorem}
The proof of Theorem \ref{Theorem: decay rate of random feature matrices and their exponential ill conditionality} relies on the spectral approximation between random features. Numerical evidence in~\cite{chen2024optimization} suggests that the exponential singular-value decay carries over to two- and three-dimensional problems, and to other analytic activation functions. Under a typical experimental setup, Fig.~\ref{Figure: Sparsity pattern of RFMtx} demonstrates that the upper bound (\ref{ineq: fast decay of sigma_m(Z)}) is sharp.

\begin{remark}
Theorem \ref{Theorem: decay rate of random feature matrices and their exponential ill conditionality} accommodates general sampling distributions, provided that $\{k_{j}\}_{j = 1}^{N} \subset (0, S)$, and allows the collocation points $\{x_{i}\}_{i = 1}^{n}$ to be selected randomly or deterministically. Via numerical integration, its final assertion can be refined to $\mb{P}(\rho \geq 1/11) \geq 1-0.72^N$ (see Appendix \ref{Appendix section: Auxiliary result in the proof of Theorem: decay rate of random feature matrices and their exponential ill conditionality}). To guarantee a unique solution for Problem (\ref{1-d elliptic equation with variable coefficients, Dirichlet BC}), $a$ and $c$ are typically assumed to have opposite signs, yielding $B \leq 0$ and $\rho \geq 1$.
\end{remark}

\begin{remark}
The proof of Theorem \ref{Theorem: decay rate of random feature matrices and their exponential ill conditionality}, particularly the construction of $\boldsymbol{\alpha}$, holds for any smooth activation function $\eta$. It reveals that the decay rate of the $m$-th singular value of the RFMtx depends heavily on $(\eta - p_{m-2})^{(L)}$, where $p_{m-2}$ is the Taylor polynomial of degree $m-2$ for $\eta$, and $L$ is the highest derivative order applied to random features (e.g., $L = 2$ in Theorem~\ref{Theorem: decay rate of random feature matrices and their exponential ill conditionality}). For trigonometric activation functions, singular values decay exponentially when $m$ and $SR$ are of the same order, and super-exponentially (scaling as $\Gamma(m)^{-1}$) for larger $m$. For the $\tanh$ activation function, exponential decay occurs when $SR < \pi/2$ (see Fig. 1 in \cite{chen2024optimization}, where $S = R = 1$). The same mechanism applies to RFMtxs involving derivatives of RFs, with the same rapid singular-value decay; see also Fig.~3 in~\cite{chen2024optimization} for numerical evidence.
\end{remark}

Next, we analyze the effect of applying PUM. For simplicity, we fix a PoU. Let $\{\xi_{j}\}_{j=1}^{J}$ be $J$ equispaced points given by $\xi_{j} = -R + 2r(j-1)$ with $r = R/(J-1)$. The PoU function centered at $\xi_{j}$ is defined as $\phi_{j}(x) = \phi\left((x - \xi_{j})/{r}\right)$, where
\begin{equation*}
\phi(t) = \textbf{1}_{[-\frac{5}{4}, -\frac{3}{4})}(t) \frac{1+\sin (2 \pi t)}{2} + \textbf{1}_{[-\frac{3}{4}, \frac{3}{4})}(t) + \textbf{1}_{[\frac{3}{4}, \frac{5}{4})}(t) \frac{1 - \sin (2 \pi t)}{2}.
\end{equation*}
Let $u_{\text{RF}}$ be defined by~\eqref{random feature model with partition of unity}. The resulting $\mf{\Phi}$ exhibits a block-sparse structure (see Fig. \ref{Figure: Sparsity pattern of RFMtx}). Based on their associated patches, the columns of $\mf{\Phi}$ can be partitioned into $J$ blocks. For the $j$-th block column, its rows are partitioned into submatrices $\mf{A}_{j}$, $\mf{B}_{j}$ and $\mf{C}_{j}$, corresponding to collocation points $x_{i}$ located in the intervals $(\xi_{j}-5r/4, \xi_{j}-3r/4)$, $[\xi_{j}-3r/4, \xi_{j}+3r/4]$ and $(\xi_{j}+3r/4, \xi_{j}+5r/4)$, respectively. Consequently, $\mf{\Phi} = [\mf{\Phi}_{1}, \mf{\Phi}_{2}, \ldots, \mf{\Phi}_{J}]$ and its nonzero block columns $\mf{D}_{j}$ can be expressed as
\begin{equation*}
\mf{\Phi} =
\begin{bmatrix}
 \mf{B}_{1} &            &        & \\
 \mf{C}_{1} & \mf{A}_{2} &        & \\
            & \mf{B}_{2} & \ddots & \\
            & \mf{C}_{2} & \ddots & \mf{A}_{J} \\
            &            &        & \mf{B}_{J}
\end{bmatrix},
\quad \text{with} \quad
\begin{aligned}
&\mf{D}_{1} = [\mf{B}_{1}^\top, \mf{C}_{1}^\top]^\top, \\
&\mf{D}_{j} = [\mf{A}_{j}^\top, \mf{B}_{j}^\top, \mf{C}_{j}^\top]^\top \\
& \text{for } 2 \leq j \leq J-1, \\
& \mf{D}_{J} = [\mf{A}_{J}^\top, \mf{B}_{J}^\top]^\top.
\end{aligned}
\end{equation*}
The entries of $\mf{D}_{j}$ correspond to the differential operator acting on a local feature $v$, evaluated at the collocation points $x_{i}$:
\begin{equation*}
\begin{aligned}
L\left(\phi_{j} v\right) &  = a \phi_{j} v^{\prime \prime} + \left(2 a \phi_{j}^{\prime} +b \phi_{j}\right) v^{\prime}+\left(a \phi_{j}^{\prime \prime} +b \phi_{j}^{\prime}+c\phi_{j}\right)  v.
\end{aligned}
\end{equation*}
The singular-value decay of $\mf{B}_{j}$ and $\mf{D}_{j}$ has been analyzed in Theorem~\ref{Theorem: decay rate of random feature matrices and their exponential ill conditionality}. Below, we analyze the relationship between the singular values of $\mf{\Phi}$ and those of its submatrices $\mf{B}_{j}$ and $\mf{D}_{j}$. This connection illustrates how PUM mitigates the rapid singular-value decay. Without loss of generality, we assume that $\mf{B}_{j}$ has at least as many rows as columns; i.e., there are sufficient collocation points per patch.
\begin{theorem} \label{Theorem: The effect of PUM on the singular value of RFMtx}
Let block diagonal matrices $\mf{B} = \operatorname{diag}(\mf{B}_{j})_{j =1}^{J}$, $\mf{D}_{o} = \operatorname{diag}(\mf{D}_{j})_{j \text{ is odd}}$,
$\mf{D}_{e} = \operatorname{diag}(\mf{D}_{j})_{j \text{ is even}},$
where $\mf{D}_{o}$ (resp. $\mf{D}_{e}$) has $N_{o}$ (resp. $N_{e}$) columns. Then,
\begin{equation}  \label{ineq: basic lower and upper bounds for RFMtx with PUM}
\sigma_{m}(\mf{B}) \leq \sigma_{m}(\mf{\Phi}) \leq  \min_{\max(1, m-N_{e}) \leq k\leq \min(m, N_{o}+1) }  \sqrt{\sigma_{k}^{2}(\mf{D}_{o}) + \sigma_{m+1-k}^{2}(\mf{D}_{e})},
\end{equation}
for any integer $1\le m\le N_e+N_o$.\footnote{We use the convention that $\sigma_j(\mf{D}_o)=0$ for $j>N_o$ and $\sigma_j(\mf{D}_e)=0$ for $j>N_e$.} In particular,  if all $N_{j}$ are equal, then
\begin{equation}  \label{ineq: simplified lower and upper bounds for RFMtx with PUM}
\min_{j} \sigma_{\lceil m/J\rceil}(\mf{B}_{j}) \leq \sigma_{m}(\mf{\Phi}) \leq \sqrt{2} \max_{j} \sigma_{\lceil m/J\rceil}(\mf{D}_{j}).
\end{equation}
\end{theorem}

In practice, applying PUM or domain decomposition causes row overlaps between adjacent submatrices (e.g., $\mf{C}_{j-1}$ and $\mf{A}_{j}$), preventing the RFMtx from being strictly block-diagonal. We analyze this practical but nontrivial case. Theorem \ref{Theorem: The effect of PUM on the singular value of RFMtx} provides simple, effective bounds relating the singular values of the global RFMtx to those of its submatrices. As shown in Fig. \ref{Figure: Sparsity pattern of RFMtx}, the singular values of the RFMtx are close to the upper bound in \eqref{ineq: basic lower and upper bounds for RFMtx with PUM}. The lower bound captures the correct decay trend but remains loose. This is because the PoU functions exhibit large first and second derivatives in the overlapping regions, causing the entries of $\mf{A}_j$ and $\mf{C}_j$ to dominate those of $\mf{B}_j$. If $\widetilde{L}$ is the identity operator, both bounds become tight (see Fig. \ref{Figure: integralbound}).

\cite{chen2024quantifytraindifficulty} showed that under identical patch sizes, local feature counts, and sampling distributions, the singular values of $\mf{D}_j$ (and $\mf{B}_j$) are highly consistent across different $j$ with high probability. Thus,~\eqref{ineq: simplified lower and upper bounds for RFMtx with PUM} indicates that the $m$-th singular value of the RFMtx approximates the $\lceil m/J \rceil$-th singular value of the submatrices. This explains why PUM reduces the singular-value decay rate by a factor slightly below the number of patches in log-scale plots (see Figs.~2 and~3 in~\citet{chen2024optimization} and Figs.~\ref{Figure: Sparsity pattern of RFMtx}, \ref{Figure: integralbound}).
Furthermore, it is straightforward to extend the proof of Theorem \ref{Theorem: The effect of PUM on the singular value of RFMtx} to high dimension and to more general PoU.
\section{Proofs of approximation rates} \label{Section: Fast convergence of approximation errors}
\subsection{Proofs of spectral convergence with respect to feature count}
\label{Subsection: Proofs of spectral convergence with respect to feature count}
We prove the assertions of \S \ref{subsection of main results: Spectral convergence with respect to the number of features} in this section. %
The approximant $u_N$ is constructed by linearly combining features as in (\ref{random feature model used in the fast convergence part}) to match the Taylor expansion of $u$ up to degree $2N-1$. This $u_N$ is adopted if its coefficients remain well-controlled; otherwise, we simply set $u_N=0$. Specifically, we approximate $u$ by the truncated Maclaurin series,  whose even and odd parts are
\begin{equation*}
\begin{aligned}
p_{u}(x) = \sum_{n=0}^{N-1} \frac{u^{(2 n)}(0)}{(2 n)!} x^{2 n} \quad\text{and}\quad  q_{u}(x) = \sum_{n=0}^{N-1} \frac{u^{(2 n+1)}(0)}{(2 n+1)!} x^{2 n+1}.
\end{aligned}
\end{equation*}
In parallel, we approximate $\cos(k_{i} \cdot)$ and $\sin(k_{i} \cdot)$ by the truncated Maclaurin series
\[
p_{i}(x) = \sum_{n=0}^{N-1}(-1)^{n} \frac{\left(k_{i} x\right)^{2 n}}{(2 n)!} \quad\text{and}\quad  q_{i}(x)=\sum_{n=0}^{N-1}(-1)^{n} \frac{\left(k_{i} x\right)^{2 n+1}}{(2 n+1)!}.
\]
Hence, we decompose the approximation error $u-u_N$ as
\begin{equation*}
\begin{aligned}
u-u_N & = (u-p_{u}-q_{u}) + \left(p_{u}-\sum_{i=1}^{N} \alpha_{i} p_{i}\right) + \left(q_{u}-\sum_{i=1}^{N} \alpha_{N+i} q_{i}\right) \\
& \quad +  \sum_{i=1}^{N} \alpha_{i} \left[p_{i}-\cos(k_{i} \cdot)\right] + \sum_{i=1}^{N} \alpha_{N+i} \left[q_{i}-\sin(k_{i} \cdot)\right].
\end{aligned}
\end{equation*}

Next, we set $p_{u}=\sum_{i=1}^{N} \alpha_{i} p_{i}$ and $q_{u}=\sum_{i=1}^{N} \alpha_{N+i} q_{i}$ to eliminate the middle two terms in the above decomposition. Define linear operators $\mathscr{T}_{e}, \mathscr{T}_{o}: C^{2N-1}(\Omega) \to \mb{R}^N$ componentwise by
$(\mathscr{T}_{e}u)_i = u^{(2i-2)}(0)$ and $(\mathscr{T}_{o}u)_i = u^{(2i-1)}(0).$
Denote the Vandermonde matrix $\mf{V}_N = \left((-k_j^2)^{i-1}\right)_{1\le i,j\le N}$  and the diagonal matrix $\mf{K}=\operatorname{diag}(k_1,k_2,\cdots,k_N)$. Matching the coefficients of $x^{n}$ leads to the linear systems %
\begin{equation} \label{eqs: linear systems for alpha^(V) obtained by matching the coefficients of x^n}
\begin{aligned}
\mf{V}_N \mf{X}=\mathscr{T}_{e}u, \quad\text{and}\quad \mf{V}_N \mf{KY}=\mathscr{T}_{o}u,  \quad\text{where } \mf{X},\mf{Y}\in\mb{R}^N.
\end{aligned}
\end{equation}
Let $\mf{\boldsymbol{\alpha}}^{V}=(\mf{V}_N^{-1}\mathscr{T}_{e}u,\mf{K}^{-1}\mf{V}_N^{-1}\mathscr{T}_{o}u)^{\top}$ be the solution\footnote{Since $\{k_i\}_{i=1}^{N}$ are $i.i.d.$ continuous random variables, $\{k_i^{2}\}_{i=1}^{N}$ are almost surely distinct positive numbers. Hence, $\mf{V}_N$ and $\mf{K}$ are almost surely invertible.} of these equations. We construct the desired $u_N$ with coefficients $\tilde{\boldsymbol{\alpha}} = \boldsymbol{\alpha}^{V}\mf{1}_{A}$, where $A$ is an event to be determined later, on which $\boldsymbol{\alpha}^{V}$ is well-controlled. To estimate the error, we employ the natural decomposition
\begin{equation} \label{eq: split error into the parts on event A and A^c}
\begin{aligned}
\mb{E} \abs{u-u_{N}}_{W^{l,p}(\bar{\Omega})}
& = \mb{P}(A^{c})\abs{u}_{W^{l,p}(\bar{\Omega})} + \mb{E} \abs{u-u_{N}}_{W^{l,p}(\bar{\Omega})}\mf{1}_{A},  \\
\end{aligned}
\end{equation}
where on the event $A$, we use the triangle inequality to get
\begin{equation}
\label{preliminary error decomposition, before pf of Theorem: main theorem of RFM, Gevrey classes, exponential convergence}
\begin{aligned}
\abs{u-u_{N}}_{W^{l,p}(\bar{\Omega})}
&\le \abs{u - p_{u} - q_{u}}_{W^{l,p}(\bar{\Omega})} + \sum_{i=1}^{N}\abs{\alpha_{i}^{V}} \abs{p_{i}-\cos(k_{i}\cdot)}_{W^{l,p}(\bar{\Omega})} \\
&\quad + \sum_{i=1}^{N}\abs{\alpha_{N+i}^{V}} \abs{q_i -\sin(k_i\cdot)}_{W^{l,p}(\bar{\Omega})}  .
\end{aligned}
\end{equation}

First, Lemma \ref{Lemma:vandermonde} provides an explicit upper bound on $\boldsymbol{\alpha}^{V}$ by inverting the Vandermonde matrix. Second, Lemma \ref{Lemma: probability bound for events A_i^c under general distributions} characterizes the events $A_{i}$ and shows that $A := \bigcap_{i=1}^{N} A_{i}$ occurs with high probability. Third, Lemma \ref{Lemma: expectation bound for k_1^(2 N) (mcQ)^(N-1), extended distribution} controls the polynomial approximation error of features on $A$, which corresponds to the summation terms in (\ref{preliminary error decomposition, before pf of Theorem: main theorem of RFM, Gevrey classes, exponential convergence}). Consequently, a combination of these three lemmas leads to Theorem \ref{Theorem: main theorem of RFM, Gevrey classes, W^(l,p) error norm, extended distribution}.
\begin{lemma}\label{Lemma:vandermonde}
Let $u\in G^s(\bar{\Omega})$ and $\left\{k_i^2\right\}_{i=1}^{N}\subset\mb{R}_{+}$ be pairwise distinct. Then, the solutions of~\eqref{eqs: linear systems for alpha^(V) obtained by matching the coefficients of x^n} satisfy, for $1\le i \le N$
\begin{equation*} \label{VandermonddBD}
\begin{aligned}
\abs{X_i} &\le  M_{u}[(2 N-2)!]^{s} \prod_{j \neq i} \frac{C_{u}^{2}+k_j^2}{\abs{k_i^{2}-k_j^{2}}},\\
\abs{Y_i} &\le\frac{M_{u} C_{u}[(2 N-1)!]^s}{k_i} \prod_{j \neq i} \frac{C_{u}^{2}+k_j^2}
{\abs{k_i^2-k_j^2}} .
\end{aligned}
\end{equation*}
\end{lemma}

\begin{proof}[Proof.]
Set $\lambda_{i} = -k_{i}^{2}$ and Vandermonde matrix $\mf{V}_{N}=(\lambda_j^{i-1})_{1\le i,j\le N}$. For $1\le i\le N$, denote by $\sigma_m^{(i)}$ the elementary symmetric polynomial of degree $m$ in the variables $\{\lambda_\ell:1\le\ell\le N, \ \ell\ne i\}$, with $\sigma_0^{(i)}=1$. If $\mf{V}_{N}^{-1}=(v_{ij})$, the standard inverse formula~\citep[p.~306]{muir1960history} for the Vandermonde matrix gives
\begin{equation*}%
v_{ij}=(-1)^{j-1}\frac{\sigma_{N-j}^{(i)}}{\prod_{\ell\ne i}(\lambda_\ell-\lambda_i)} .
\end{equation*}
Then, it follows from $\mf{X}=\mf{V}_N^{-1}\mathscr{T}_{e}u$ and $u\in G^s(\bar{\Omega})$ that
\begin{equation*}
\begin{aligned}
\abs{X_{i}} & %
\le  \frac{\sum_{j=1}^{N}\left|\sigma_{N-j}^{(i)}\right| M_{u} C_{u}^{2 j-2}[(2 j-2)!]^{s}}{\prod_{\ell \neq i} \left|k_{i}^{2} - k_{\ell}^{2}\right|} . \\
\end{aligned}
\end{equation*}
Using $[(2 j-2)!]^{s} \le[(2 N-2)!]^{s}$ and the identity $\prod_{j \neq i}\left(x + k_{j}^{2}\right) = \sum_{j=0}^{N-1} \left|\sigma_{N-j}^{(i)}\right| x^{j}$ evaluated at $x=C_{u}^{2}$, we obtain the asserted estimate for $X_i$. The bound for $(\mf{K}\mf{Y})_i$ is derived by the same argument. Finally, applying $\mf{K}^{-1}$ yields the desired estimate for $Y_i$.
\end{proof}

\begin{lemma} \label{Lemma: probability bound for events A_i^c under general distributions}
Under Assumption \ref{Assumption: extended sampling distribution}, for any $C\ge0$ and $c>2$, define the events
\begin{equation} \label{eq: events A_i defined for control alpha, extended distribution}
\begin{aligned}
A_{i}{:}=\left\{\prod_{j \neq i} \frac{C^{2}+k_{j}^{2}}{\left|k_{i}^{2}-k_{j}^{2}\right|} < \left[ e^{c} \mc{Q}(|k_{i}|/\mu;C/\mu)\right]^{N-1}\right\},
\end{aligned}
\end{equation}
where $\mc{Q}(k;w)$ is defined as
\begin{equation*}
\left\{\begin{aligned}
& \varpi^{\frac{c}{c-1}}\max\left(\frac{1 + w^2}{2(k^2 - 1)}+\frac{w^2}{2k^2}, \frac{28(1 + w^2)}{23k^{2}}\right), & k>4/3,\\
& 2\varpi^{\frac{c}{c-1}} \max\left(\frac{w^2 + a^2}{k+a},\frac{w^2 + b^2}{k+b}, \frac{w^2+1}{k+1}\right), & k\le 4/3,
\end{aligned}\right.
\end{equation*}
and $a = \max(k-1/2,0)$, $b = a+1$. Then,
\[
\mb{P}\left(A_{i}^{c}\right) \leq e^{-[c-1-\ln (1+c)](N-1)}.
\]
\end{lemma}

\begin{lemma} \label{Lemma: expectation bound for k_1^(2 N) (mcQ)^(N-1), extended distribution}
Under Assumption \ref{Assumption: extended sampling distribution}, for any $C\ge0$, there holds
\begin{equation*}
\begin{aligned}
\mb{E}\left[k_{1}^{2 N} \mc{Q}^{N-1}(|k_{1}|/\mu;C/\mu) \right] %
& \le \mu^{2} \left[5\varpi^{\frac{c}{c-1}}\max(C, \mu)^2\right]^{N-1} .
\end{aligned}
\end{equation*}
\end{lemma}

The proofs of Lemma~\ref{Lemma: probability bound for events A_i^c under general distributions} and Lemma~\ref{Lemma: expectation bound for k_1^(2 N) (mcQ)^(N-1), extended distribution} are postponed to Appendix B. We are ready to prove Theorem \ref{Theorem: main theorem of RFM, Gevrey classes, W^(l,p) error norm, extended distribution}.
\begin{proof}[Proof of Theorem \ref{Theorem: main theorem of RFM, Gevrey classes, W^(l,p) error norm, extended distribution}.]
The case $N=1$ is trivial; we assume that $N\ge 2$. Let $A = \bigcap_{i=1}^{N} A_{i}$ with $A_i$ defined in \eqref{eq: events A_i defined for control alpha, extended distribution}.  $u_{N} = \mathcal{R}_N[\{k_{i}\}_{i = 1}^{N},\mu,C_u,c] u$ with
\begin{equation*}
\mathcal{R}_N[\{k_{i}\}_{i = 1}^{N},\mu,C,c] = S_N\circ (\mf{V}_N^{-1}\mathscr{T}_{e},\mf{K}^{-1}\mf{V}_N^{-1}\mathscr{T}_{o})^{\top}\mf{1}_{A}.
\end{equation*}
It follows from $u\in G^s(\bar{\Omega})$ and Lemma \ref{Lemma: probability bound for events A_i^c under general distributions} that
\begin{equation} \label{ineq: approximation error on A^c, extended distribution}
\begin{aligned}
\mb{P}(A^{c})\abs{u}_{W^{l,p}(\bar{\Omega})}  %
\le (2R)^{1/p} M_{u} C_{u}^{l}(l!)^s N e^{-(c-1-\ln(c+1))(N-1)} .
\end{aligned}
\end{equation}
Since $p_{u} + q_{u}$ is the Taylor polynomial of order $2N-1$ of $u$, for $l\le 2N$,
\begin{equation*}
\abs{u - p_{u} - q_{u}}_{W^{l,p}(\bar{\Omega})} \le  \frac{2^{1/p} R^{2N-l+1/p} M_u C_u^{2N} [(2N)!]^s}{[(2N-l)p+1]^{1/p} (2N-l)!} .
\end{equation*}
It follows from Lemma \ref{Lemma:vandermonde} and Taylor approximation that
\begin{equation*}
\begin{aligned}
\abs{\alpha_{i}^{V}} \abs{p_{i}-\cos(k_{i}\cdot)}_{W^{l,p}(\bar{\Omega})} & \le \frac{2^{1/p} R^{2N-l+1/p} M_u [(2N-2)!]^s}{[(2N-l)p+1]^{1/p} (2N-l)!} k_{i}^{2N} \prod_{j \neq i} \frac{k_j^2 + C_u^2}{|k_j^2 - k_i^2|},
\end{aligned}
\end{equation*}
where on the event $A$, according to Lemma \ref{Lemma: probability bound for events A_i^c under general distributions} and Lemma \ref{Lemma: expectation bound for k_1^(2 N) (mcQ)^(N-1), extended distribution},
\begin{equation*}
\begin{aligned}
\mb{E} \left[k_i^{2N} \prod_{j \neq i} \frac{k_j^2 + C_u^2}{|k_j^2 - k_i^2|} \mf{1}_{A} \right] & \le e^{c(N-1)}\mb{E}\left[k_{1}^{2 N} \mc{Q}(|k_{1}|/\mu;C_{u}/\mu)^{N-1} \right] \\
& \le \mu^{2} \left[5e^{c}\varpi^{\frac{c}{c-1}}\max(C_{u}, \mu)^2\right]^{N-1} .
\end{aligned}
\end{equation*}
Analogously, we obtain
\begin{equation*}
\begin{aligned}
\abs{\alpha_{N+i}^{V}} \abs{q_i -\sin(k_i\cdot)}_{W^{l,p}(\bar{\Omega})} %
& \lesssim_{l} \frac{R^{2N+1-l+1/p} M_{u} C_{u} [(2N-1)!]^s}{N^{1/p}(2N+1-l)!} k_{i}^{2N} \prod_{j \neq i} \frac{k_j^2 + C_u^2}{|k_j^2 - k_i^2|} .
\end{aligned}
\end{equation*}
Substituting the above four bounds into (\ref{preliminary error decomposition, before pf of Theorem: main theorem of RFM, Gevrey classes, exponential convergence}), we obtain
\begin{equation*}
\begin{aligned}
\mb{E} \abs{u-u_{N}}_{W^{l,p}(\bar{\Omega})}\mf{1}_{A} %
& \lesssim_{l}  \frac{M_u \varpi^{l}\max(C_{u}, \mu)^{l}R^{1/p}}{N^{(2-l)s+1/p}} e^{c(N-1)} \tau^{2N-l} (1+\tau) . \\
\end{aligned}
\end{equation*}
A combination of (\ref{eq: split error into the parts on event A and A^c}), (\ref{ineq: approximation error on A^c, extended distribution}) and the above estimate completes the proof.
\end{proof}
Matching the two terms in (\ref{error bound in Theorem:Theorem: main theorem of RFM, Gevrey classes, W^(l,p) error norm}) with a suitable constant $c$, we prove Corollary \ref{Corollary: Gevrey classes s<1, super exponential convergence for N, W^(l,p) error norm}.
\begin{proof}[Proof of Corollary \ref{Corollary: Gevrey classes s<1, super exponential convergence for N, W^(l,p) error norm}.]
We only need to prove the case when $\epsilon<1-s$. First, we prove the case when $\sigma=l$ is an integer. We choose %
$c = 3$ when $N \le N_{s,l} := 2^{-1}[e^{2(l+1)/({1-s}) +1} +l].$
Since $N$ is bounded, there exists $C$ independent of $N$ such that
\[
Ne^{2(\ln(2)-1)(N-1)} + N^{l+1}e^{3(N-1)} \tau^{2N-l} \le  C \Gamma(2N-l+1)^{-(1-s)/2 + \epsilon}.
\]
For $N > N_{s,l}$, we choose $c = (1-s)(2N-2)^{-1} \ln \Gamma(2 N-l+1).$ %
Then, $2 < c < 2\ln N$ by \eqref{ineq: two sided bounds for factorial, Stirling's formula}, and it follows from
$Ne^{N}(1+2\ln N)^{N} \Gamma(2 N-l+1)^{-\epsilon/4} \to 0$ as $N\to \infty$ that
\begin{equation*}
\begin{aligned}
Ne^{-[c-1-\ln(1+c)](N-1)} %
& \lesssim_{\epsilon,l} \Gamma(2 N-l+1)^{-(1-s)/2+\epsilon/{4}}.
\end{aligned}
\end{equation*}
Proceeding along a similar argument, we have
\begin{equation*}
\begin{aligned}
N^{l+1} e^{c(N-1)} \tau^{2N-l} & \le N^{l+1} \left[\sqrt{5}\varpi \max(C_{u}, \mu) R\right]^{2N-l} \Gamma(2N-l+1)^{-\frac{1}{2}(1-s)} \\
& \lesssim_{l,\epsilon,\varpi,C_{u},\mu,R} \Gamma(2 N-l+1)^{-(1-s)/2+\epsilon/{4}} .
\end{aligned}
\end{equation*}
A combination of the above two estimates yields
\begin{equation*}
\begin{aligned}
\mb{E} \abs{u-u_{N}}_{W^{l,p}(\bar{\Omega})}
& \lesssim_{l,\epsilon,\varpi,C_{u},\mu,R} \Gamma(2 N-l+1)^{-(1-s)/2+\epsilon/{4}} .
\end{aligned}
\end{equation*}
The same estimate for $\mb{E}\|u-u_{N}\|_{W^{l,p}(\bar{\Omega})}$ follows from summing over all semi-norms of order less than $l$. Note that $\Gamma(N)^2\le (2N)^l\Gamma(2N-l+1),$ %
which gives the desired estimate for integer $\sigma$. Finally, the estimate for the non-integer $\sigma$ follows from the Gagliardo-Nirenberg interpolation inequality in Theorem \ref{Theorem: Gagliardo–Nirenberg inequalities and non-inequalities: The full story}, which completes the proof.
\end{proof}

With the aid of Fourier truncation and Theorem \ref{Theorem: approximate W^(s,p)-func by G^0 func on Lipschitz domain or whole space with explicit error bound}, we prove Corollary~\ref{Corollary: spectral convergence when solution u in Analytic class, Gevrey class s>1 or Sobolev space, W^(l,p) error norm}. %
\begin{proof}[Proof of Corollary \ref{Corollary: spectral convergence when solution u in Analytic class, Gevrey class s>1 or Sobolev space, W^(l,p) error norm}.]
We first prove for the case $R=1$ and $\sigma\in\mb{N}$. Then we apply the interpolation theorem, Theorem~\ref{Theorem: Gagliardo–Nirenberg inequalities and non-inequalities: The full story}, and a scaling argument to obtain the desired estimate. In each case (a)-(c), we construct an approximation $u_{\mu}= \mathcal{T}_{\mu}u$ for $u$ with $u_{\mu} \in G^{0}(\bar{\Omega})$ and Gevrey constant $C_{u_{\mu}}\le\mu$, and construct $u_{N}$ as in Theorem \ref{Theorem: main theorem of RFM, Gevrey classes, W^(l,p) error norm, extended distribution} to approximate $u_{\mu}$. %
The error admits the following decomposition
\begin{equation*}
\begin{aligned}
\mb{E}\|u-u_{N}\|_{W^{\sigma,p}(\Omega)} & \le \|u-u_{\mu}\|_{W^{\sigma,p}(\Omega)} + \mb{E}\|u_{\mu}-u_{N}\|_{W^{\sigma,p}(\Omega)} .
\end{aligned}
\end{equation*}

The following proof is divided into four steps.

\textbf{Step 1. Construct $u_{\mu}$ and control its Gevrey constant.}
For (a), the assumption gives the required analytic extension. For (b), Proposition  \ref{Proposition: explicit Fourier decay for extended Gevrey functions} provides a compactly supported smooth extension, which we also denote by $u$. We construct $u_{\mu}$ as its low-frequency part $u_{\mu} := (2\pi)^{-1/2}\int_{-\mu}^{\mu} \hat{u}(\xi) e^{i \xi x} \mr{d} \xi$. For bandlimited function $u_{\mu}$, Lemma~\ref{Lemma: trigonometric functions, exponential functions, polynomials in G^0, dimension d} shows that
$u_{\mu}\in G^0(\bar{\Omega})$ with $M_{u_{\mu}}\lesssim\|\hat{u}\|_{L^{1}(\mb{R})}$ and $C_{u_{\mu}}\le \mu$.
For (c), we select $u_\mu$ as $u_{\varepsilon}$ provided by Theorem \ref{Theorem: approximate W^(s,p)-func by G^0 func on Lipschitz domain or whole space with explicit error bound} with $\varepsilon = (C_{s}\mu)^{-1}$. Then, $M_{u_{\mu}}\lesssim_{s,p} \mu^{1/p} \|u\|_{W^{s,p}(\Omega)}$ and $C_{u_{\mu}}\le \mu$.

\textbf{Step 2. Bound $\|u-u_{\mu}\|_{W^{\sigma,p}(\Omega)}$.}
For (a) and (b), it follows from Hausdorff-Young inequality that for any $l\in\mb{N}$,
\begin{equation*}
\begin{aligned}
|u-u_{\mu}|_{W^{l,p}(\Omega)}
& \lesssim_{l, p} |u-u_{\mu}|_{W^{l,\infty}(\mb{R})} \lesssim_{l, p} \||\cdot|^{l}\hat{u}\mf{1}_{\{|\cdot|>\mu\}}\|_{L^{1}(\mb{R})} . \\
\end{aligned}
\end{equation*}
Under the assumptions in (a), Paley-Wiener theorem~\citep[Theorem IX.13]{Reed1978Methods} %
ensures $v_{\rho} := e^{\rho |\cdot|} \hat{u} \in L^{2}(\mb{R})$. Then, H\"older's inequality gives
\begin{equation*}
\begin{aligned}
\||\cdot|^{l}\hat{u}\mf{1}_{\{|\cdot|>\mu\}}\|_{L^{1}(\mb{R})} %
& \lesssim_{l,p} \rho^{-l} e^{-\rho\mu/2} \|v_{\rho}\|_{L^{2}(\mb{R})}.
\end{aligned}
\end{equation*}
A combination of the above two estimates yields
\begin{equation} \label{ineq: W^(l,p)(-1,1) seminorm bound for u-u_mu, assumption (a)}
\begin{aligned}
|u-u_{\mu}|_{W^{l,p}(\Omega)}
& \lesssim_{l, p} \rho^{-l} e^{-\rho\mu/2} \|v_{\rho}\|_{L^{2}(\mb{R})} . \\
\end{aligned}
\end{equation}
Under the assumptions in (b), Proposition \ref{Proposition: explicit Fourier decay for extended Gevrey functions} shows the extension of $u$ satisfies $\abs{\hat{u}(\xi)}\le C e^{-\kappa |\xi|^{1/s}}$ with $C = 2 e^{7s/4} (1+\rho)M_{u}$ and $\kappa = (3s/4) (C_u + 2 \zeta_s/\rho)^{-1/s}$. Then, performing the change of variables $t = \kappa \xi^{1/s}$, we obtain
\begin{equation*}
\begin{aligned}
\||\cdot|^{l}\hat{u}\mf{1}_{\{|\cdot|>\mu\}}\|_{L^{1}(\mb{R})} & \le C \||\cdot|^{l} e^{-\kappa |\cdot|^{1/s}}\mf{1}_{\{|\cdot|>\mu\}}\|_{L^{1}(\mb{R})} \\
& \lesssim_{s,l,p} C\kappa^{-s-sl} e^{-\kappa \mu^{1/s}/2} , \\
\end{aligned}
\end{equation*}
which implies %
\begin{equation} \label{ineq: W^(l,p)(-1,1) seminorm bound for u-u_mu, assumption (b)}
\begin{aligned}
|u-u_{\mu}|_{W^{l,p}(\Omega)}
& \lesssim_{s, l, p} M_{u} \kappa^{-s-sl} e^{-\kappa \mu^{1/s}/2}  .
\end{aligned}
\end{equation}
Under the assumption of (c), $\|u-u_{\mu}\|_{W^{\sigma,p}(\Omega)}$ has been estimated in Theorem \ref{Theorem: approximate W^(s,p)-func by G^0 func on Lipschitz domain or whole space with explicit error bound}. %

\textbf{Step 3. Bound $\mb{E}\|u_{\mu}-u_{N}\|_{W^{\sigma,p}(\Omega)}$.}
We apply Theorem \ref{Theorem: main theorem of RFM, Gevrey classes, W^(l,p) error norm, extended distribution} to approximate $u_{\mu}$ with the choice $c = 29/10$. By the bound \eqref{ineq: two sided bounds for factorial, Stirling's formula}, for $l \le N$, %
$$\tau %
\le \frac{\sqrt{5} e}{34},
\quad  e^{-c+1+\ln(c+1)} < e^{-0.53} \quad \text{and} \quad e^{c} \tau^{2} < \frac{5 e^{2+c}}{34^{2}}< e^{-0.54}.$$
which guarantees %
\begin{equation} \label{ineq: the error in approximating G^0 part under W^(l,p)(-1,1) seminorm}
\begin{aligned}
\mb{E} \abs{u_{\mu}-u_{N}}_{W^{l,p}(\Omega)}
& \lesssim_{l}  M_{u_{\mu}} N^{l+1} \left(e^{-0.53N} + N^{-1/p}e^{-0.54N}\right) \\
& \lesssim_{l}  M_{u_{\mu}} e^{-N/2} .\\
\end{aligned}
\end{equation}

\textbf{Step 4. Combine the error bounds from the previous two steps.}
We are now ready to derive the convergence rate when $R=1$ and $\sigma$ is an integer. For (a), Lemma \ref{Lemma: trigonometric functions, exponential functions, polynomials in G^0, dimension d} (4) ensures $M_{u_{\mu}}\lesssim\|\hat{u}\|_{L^{1}(\mb{R})} \lesssim \rho^{-1/2} \|v_{\rho}\|_{L^{2}(\mb{R})}$. Then, a combination of \eqref{ineq: W^(l,p)(-1,1) seminorm bound for u-u_mu, assumption (a)} and \eqref{ineq: the error in approximating G^0 part under W^(l,p)(-1,1) seminorm} gives %
\begin{equation*}
\begin{aligned}
\mb{E}|u-u_{N}|_{W^{l,p}(\Omega)} %
& \lesssim_{l, p} \left(\rho^{-l} + 1\right)  \rho^{-1/2}\|v_{\rho}\|_{L^{2}(\mb{R})} e^{-\min\left(\frac{\rho}{34\varpi},1\right)\frac{N}{2}}   ,
\end{aligned}
\end{equation*}
for $l\le N$. With $l=0$ and $l=\sigma$, we obtain the desired estimate for an integer $\sigma$. For the non-integer $\sigma$, it follows from the Gagliardo-Nirenberg inequality in Theorem \ref{Theorem: Gagliardo–Nirenberg inequalities and non-inequalities: The full story} that
\begin{equation} \label{ineq: Gagliardo-Nirenberg inequality and Holder}
\begin{aligned}
\mb{E}\|u-u_{N}\|_{W^{\sigma,p}(\Omega)} & \lesssim_{\sigma, p} \mb{E}\left(\|u-u_{N}\|_{W^{\lceil\sigma\rceil,p}(\Omega)}^{\theta} \|u-u_{N}\|_{L^{p}(\Omega)}^{1-\theta}\right) \\
& \lesssim_{\sigma, p} \left(\mb{E}\|u-u_{N}\|_{W^{\lceil\sigma\rceil,p}(\Omega)}\right)^{\theta} \left(\mb{E}\|u-u_{N}\|_{L^{p}(\Omega)}\right)^{1-\theta},
\end{aligned}
\end{equation}
where $\theta = \sigma/\lceil\sigma\rceil$. With the fact
$$\sup_{t> 0}\frac{2^{1-\theta}(1 + t^{\lceil\sigma\rceil})^{\theta}}{1 + t^{\sigma}} =2^{1-\theta}\le 2,$$
we extend the desired estimate from integer $\sigma$ to non-integer $\sigma$.

For (b), Lemma \ref{Lemma: trigonometric functions, exponential functions, polynomials in G^0, dimension d} (4) ensures  $M_{u_{\mu}}\lesssim\|\hat{u}\|_{L^{1}(\mb{R})} \lesssim_{s} M_{u}\kappa^{-s} $. Then, it follows from \eqref{ineq: W^(l,p)(-1,1) seminorm bound for u-u_mu, assumption (b)} and \eqref{ineq: the error in approximating G^0 part under W^(l,p)(-1,1) seminorm} that for $l\le N$,
\begin{equation*}
\begin{aligned}
\mb{E}|u-u_{N}|_{W^{l,p}(\Omega)} %
&  \lesssim_{s, l, p} \left(\kappa^{-sl} e^{-\kappa \mu^{1/s}/2}  + e^{-N/2} \right) \kappa^{-s} M_{u} .
\end{aligned}
\end{equation*}
Since $ -N + \kappa(34\varpi)^{-1/s} N^{1/s}\le \exp \left(s \kappa^{\frac{s}{s-1}} \right)$ for all $N\ge 0$ and $\kappa\lesssim_s 1$, we have $e^{-N/2} \lesssim_{s} e^{-\kappa \mu^{1/s}/2}$ and
\begin{equation*}
\begin{aligned}
\mb{E}|u-u_{N}|_{W^{l,p}(\Omega)} &  \lesssim_{s, l, p} \left(\kappa^{-sl} + 1 \right) \kappa^{-s} M_{u} e^{-\kappa \mu^{1/s}/2} ,
\end{aligned}
\end{equation*}
which yields the desired estimate in (b) for integer $\sigma.$ Proceeding along the same argument to prove (a), we extend the desired estimate from integer $\sigma$ to non-integer $\sigma$ utilizing Theorem \ref{Theorem: Gagliardo–Nirenberg inequalities and non-inequalities: The full story}. For (c), Theorem \ref{Theorem: approximate W^(s,p)-func by G^0 func on Lipschitz domain or whole space with explicit error bound} demonstrated $M_{u_{\mu}}\lesssim_{s,p} \mu^{1/p} \|u\|_{W^{s,p}(\Omega)}$. Since $\mu\le N$, combining Estimate 1 of Theorem \ref{Theorem: approximate W^(s,p)-func by G^0 func on Lipschitz domain or whole space with explicit error bound} with \eqref{ineq: the error in approximating G^0 part under W^(l,p)(-1,1) seminorm} gives for $l\le N$,
\begin{equation*}
\begin{aligned}
\mb{E}\|u-u_{N}\|_{W^{l,p}(\Omega)}
& \lesssim_{s,l,p} \varpi^{s-l}N^{l-s} \|u\|_{W^{s,p}(\Omega)}  .
\end{aligned}
\end{equation*}
Similarly, we interpolate this desired estimate from integer $l$ to all $\sigma\le N$.

Finally, for general $R>0$, %
a standard scaling argument completes the proof.
\end{proof}

\subsection{Proofs of convergence rates with respect to patch size}
\label{Subsection: Proofs of convergence rates with respect to patch size}

In this part, we prove Theorem \ref{Theorem: error bound for RFM combined with PUM, Gevrey class s leq 1}, Corollary \ref{Corollary: error bound for RFM combined with PUM, Gevrey class, bounded hidden weights} and Corollary \ref{Corollary: convergence for PUM when solution u in Analytic class, Gevrey class s>1 or Sobolev space}. The proof consists of two steps: First, we show how the global approximation error is controlled by the local errors on individual patches, as stated in Proposition \ref{Proposition: Error decomposition of PUM from global to local, 1 dimension}.\footnote{This proposition is adapted from \cite{Melenk2005meshless}; the general $1\le p\le \infty$ case follows from Theorem 4.1 and Remark 4.1 therein, or can be deduced via direct calculation following their proof strategy.} Second, we apply Theorem \ref{Theorem: main theorem of RFM, Gevrey classes, W^(l,p) error norm, extended distribution} to bound the local approximation error.
\begin{proposition}[\citealp{Melenk2005meshless}] \label{Proposition: Error decomposition of PUM from global to local, 1 dimension}
Suppose $u \in W^{m,p}(\Omega)$ with $m \in \mathbb{N}$ and Assumption \ref{Assumption: PoU with Finite Overlapping Condition and Derivative Bounds} holds. $u_{RF}$ is of the form (\ref{random feature model with partition of unity}). Denote the local error on the patch $\Omega_j$ as $e_j = v_j - u$. Then, for $1 \le p < \infty$, %
\begin{equation*}
\begin{aligned}
& \|u-u_{RF}\|_{W^{m,p}(\Omega)} %
\lesssim_{m,p} C_{\phi} M^{1-1/{p}} \left(\sum_{j=1}^J \sum_{l = 0}^{m} r_j^{(l-m)p} |e_j|_{W^{l, p}(\Omega_j)}^{p}\right)^{1/p}, \\
& \|u-u_{RF}\|_{W^{m,\infty}(\Omega)} %
\lesssim_{m} C_{\phi} M \max_{1 \le j \le J} \max_{0\le l \le m} r_j^{l-m} |e_j|_{W^{l, \infty}(\Omega_j)}.
\end{aligned}
\end{equation*}
\end{proposition}

\begin{proof}[Proof of Theorem \ref{Theorem: error bound for RFM combined with PUM, Gevrey class s leq 1}.]
We begin by taking the expectation of the global error decomposition bounds provided in Proposition \ref{Proposition: Error decomposition of PUM from global to local, 1 dimension}. By the linearity of expectation, we have
\begin{equation*}
\begin{aligned}
& \mb{E}\|u-u_{RF}\|_{W^{m,p}(\Omega)} \lesssim_{m,p} C_{\phi} M^{1-\frac{1}{p}} \sum_{j=1}^J \sum_{l = 0}^{m} r_j^{l-m} \mb{E}|e_j|_{W^{l,p}(\Omega_j)} . \\
\end{aligned}
\end{equation*}
Next, we apply the local approximation bound from Theorem \ref{Theorem: main theorem of RFM, Gevrey classes, W^(l,p) error norm, extended distribution} to each patch $\Omega_j$. Note that for a fixed number of features, the proof and conclusion of Theorem \ref{Theorem: main theorem of RFM, Gevrey classes, W^(l,p) error norm, extended distribution} hold for all Gevrey indices $s\ge0$, with constants allowed to depend on the number of features. Let the exponent constant $K_j = c_j-1-\ln(1+c_j)$. It follows from $\varpi^l \max(C_{u, \bar{\Omega}_j}, \mu_j)^l = A_j^l r_j^{-l}$ that
\begin{equation*}
\begin{aligned}
\mb{E}|e_j|_{W^{l,p}(\Omega_j)} %
&\lesssim_{l,\mathbf{N},s,p} M_u r_j^{1/p-l} \left( A_j^l e^{-K_j(N_j-1)}\mf{1}_{\{N_j\ge 2\}}+  e^{c_j(N_j-1)} A_j^{2N_j} (1+A_j) \right).
\end{aligned}
\end{equation*}
Substituting this local bound into the global $W^{m,p}(\Omega)$ estimate yields
\begin{equation*}
\begin{aligned}
\mb{E}\|u-u_{RF}\|_{W^{m,p}(\Omega)} %
&\lesssim_{m,p,s,\mathbf{N}} C_{\phi} M^{1-\frac{1}{p}} M_u \sum_{j=1}^J r_j^{\frac{1}{p}-m} \left( e^{-K_j(N_j-1)} \mf{1}_{\{N_j\ge 2\}} + e^{c_j(N_j-1)} A_j^{2N_j} \right) . \\
\end{aligned}
\end{equation*}
To optimize this bound, we select $c_j = -\frac{N_j}{N_j-1} \ln A_j$ for $N_j\ge2$, and $c_j = 3$ for $N_j=1$. This choice gives
\begin{equation*}
e^{-K_j(N_j-1)}\mf{1}_{\{N_j\ge 2\}} + e^{c_j(N_j-1)} A_j^{2N_j}   \lesssim_{\mathbf{N}} A_j^{N_j} |\ln A_j|^{N_j-1}.
\end{equation*}
Finally, invoking the overlapping condition $\sum_{j=1}^J r_j \le MR$, we conclude
\begin{equation*}
\begin{aligned}
\mb{E}\|u-u_{RF}\|_{W^{m,p}(\Omega)}  %
&\lesssim_{m,p,s,\mathbf{N}} C_{\phi} M^{2-\frac{1}{p}} R M_u \max_{1\le j\le J}  r_j^{\frac{1}{p}-m-1} A_j^{N_j} |\ln A_j|^{N_j-1} ,
\end{aligned}
\end{equation*}
which completes the proof.
\end{proof}

\begin{proof}[Proof of Corollary \ref{Corollary: error bound for RFM combined with PUM, Gevrey class, bounded hidden weights}.]
We require a modified version of Theorem \ref{Theorem: main theorem of RFM, Gevrey classes, W^(l,p) error norm, extended distribution}. We define the modified RF approximation $u_{N} = \tilde{\mathcal{R}}_N[\{k_{i}\}_{i = 1}^{N},\mu,C_u,c] u$ by assigning the coefficients
\begin{equation*}
\boldsymbol{\alpha} = \begin{pmatrix} \mf{V}_N^{-1}\mathscr{T}_{e}u \\ \mf{K}^{-1}\mf{V}_N^{-1}\mathscr{T}_{o}u \end{pmatrix} \mf{1}_{A} + \begin{pmatrix} J_1\mathscr{T}_{e}u \\ \mf{K}^{-1}J_1\mathscr{T}_{o}u \end{pmatrix} \mf{1}_{A^c},
\end{equation*}
where $J_1 \in \mb{R}^{N \times N}$ is the matrix with $1$ at the $(1,1)$-position and zeros elsewhere. This means on the rare event $A^c$, $u_N$ uses the features corresponding to $k_1$ to match the zeroth- and first-order Taylor expansion of $u$, i.e.,
\begin{equation*}
u_N(x)\mf{1}_{A^c} = \left[ u(0)\cos(k_1 x) + k_1^{-1}u'(0)\sin(k_1 x) \right] \mf{1}_{A^c}.
\end{equation*}
It suffices to re-estimate $\mb{E} \left[ \abs{u-u_{N}}_{W^{l,p}(\bar{\Omega})}\mf{1}_{A^c} \right]$. For $l \le 1$, analogously to (\ref{preliminary error decomposition, before pf of Theorem: main theorem of RFM, Gevrey classes, exponential convergence}), by employing $T_1(x) = u(0) + u'(0)x$ as an intermediate approximation, we obtain $$\abs{u-u_N}_{W^{l,p}} \lesssim_l M_u \max(C_u, S)^2 R^{2-l+1/p} (1+R).$$ For $l \ge 2$, the triangle inequality yields $$\abs{u-u_N}_{W^{l,p}} \le \abs{u}_{W^{l,p}} + \abs{u_N}_{W^{l,p}} \lesssim_l M_u \max(C_u, S)^l R^{1/p}.$$ Combining both cases, we arrive at
\begin{equation*}
\abs{u-u_N}_{W^{l,p}(\bar{\Omega})} \mf{1}_{A^c} \lesssim_l M_u \max(C_u, S)^{\max(l, 2)} R^{\max(2-l, 0)+1/p} (1+R) \mf{1}_{A^c},
\end{equation*}
which yields a modified version of Theorem \ref{Theorem: main theorem of RFM, Gevrey classes, W^(l,p) error norm, extended distribution}. In the proof of Theorem \ref{Theorem: error bound for RFM combined with PUM, Gevrey class s leq 1}, applying this modified version of Theorem \ref{Theorem: main theorem of RFM, Gevrey classes, W^(l,p) error norm, extended distribution} to bound the local error on each patch gives
\begin{equation*}
\begin{aligned}
\mb{E}|e_j|_{W^{l,p}(\Omega_j)}
&\lesssim r_j^{1/p-l} \left( r_j^{\max(2, l)} e^{-K_j(N_j-1)}\mf{1}_{\{N_j\ge 2\}}+  e^{c_j(N_j-1)} r_j^{2N_j} \right),
\end{aligned}
\end{equation*}
where the implied constant is independent of $r_j$. Proceeding exactly as in the proof of Theorem \ref{Theorem: error bound for RFM combined with PUM, Gevrey class s leq 1}, we obtain the desired estimate.
\end{proof}

\begin{proof}[Proof of Corollary \ref{Corollary: convergence for PUM when solution u in Analytic class, Gevrey class s>1 or Sobolev space}.]
As in Corollary \ref{Corollary: spectral convergence when solution u in Analytic class, Gevrey class s>1 or Sobolev space, W^(l,p) error norm}, the proof relies on a smoothing argument combined with an interpolation technique. We first bound the error using the triangle inequality
\begin{equation*}
\mb{E}\|u - u_{RF}\|_{W^{\sigma,p}(\Omega)} \le \|u - u_\varepsilon\|_{W^{\sigma,p}(\Omega)} + \mb{E}\|u_\varepsilon - u_{RF}\|_{W^{\sigma,p}(\Omega)},
\end{equation*}
where $u_\varepsilon \in G^0(\bar{\Omega})$ is a smooth approximation of $u$ constructed via Theorem \ref{Theorem: approximate W^(s,p)-func by G^0 func on Lipschitz domain or whole space with explicit error bound}, satisfying $M_{u_\varepsilon}\lesssim_{s,\sigma,p,R}  \varepsilon^{-1/p}\|u\|_{W^{s,p}(\Omega)}$ and $C_{u_\varepsilon,\bar\Omega_j}\lesssim_{s} \varepsilon^{-1}$.

Next, we approximate $u_\varepsilon$ using Theorem \ref{Theorem: error bound for RFM combined with PUM, Gevrey class s leq 1}. Setting $\varepsilon := \nu^{-1} r^{\chi}$ with $\chi:=1-(s+1)/(N_{\mathrm{loc}}+s-\sigma+1/p)$, the condition $t\ge -\chi$ ensures that $\varepsilon^{-1}\ge \nu r^t$. By quasi-uniformity, the local scaling factors satisfy $A_j \lesssim_{s, \varpi, \nu} r^{1-\chi}$. Consequently, for sufficiently small $r$, we have $A_j\le e^{-2}$ and $|\ln A_j| \lesssim_{s} |\ln (\varpi\nu r^t)|$ for all $j$. For any integer $m \le N_{\mathrm{loc}}-1$, we use Theorem \ref{Theorem: error bound for RFM combined with PUM, Gevrey class s leq 1} to obtain
\begin{equation} \label{ineq: RFM error integer m}
\begin{aligned}
\mb{E}\|u_\varepsilon-u_{RF}\|_{W^{m,p}(\Omega)}
& \lesssim r^{\sigma-m+\frac{(N_{\mathrm{loc}}+1/p-\sigma-1)(s-\sigma)}{N_{\mathrm{loc}}+s-\sigma+1/p}}
|\ln r|^{N_{\mathrm{loc}}-1}\|u\|_{W^{s,p}(\Omega)} ,
\end{aligned}
\end{equation}
where the implied constant is independent of $r$. If $\sigma$ is an integer, we simply take $m=\sigma$. If $\sigma \in (0, N_{\mathrm{loc}}-1)$ is not an integer, recalling \eqref{ineq: Gagliardo-Nirenberg inequality and Holder}, Theorem \ref{Theorem: Gagliardo–Nirenberg inequalities and non-inequalities: The full story} and H\"older's inequality ensure that the above estimate holds with $m$ replaced by $\sigma$, which together with the approximation error bound~\eqref{ineq: Approximation Error Bound for Gevrey approximation} yields the desired estimate \eqref{eq in Corollary: convergence for PUM when solution u in Sobolev space}.
\end{proof}

\section{Proofs for RFM solvers: error estimates and conditioning} \label{section: Exponential ill conditionality of random feature matrices}
\subsection{Proofs of abstract error estimates}
Theorem \ref{Theorem: abstract error estimate for Strong-form RFM} follows from Proposition \ref{Proposition: control H^2 error by loss} and Lemma \ref{Lemma: Lv and Bv bounded by H2 and H1 norm of v}.

\begin{proposition} \label{Proposition: control H^2 error by loss}
Under the assumptions of Theorem~\ref{Theorem: abstract error estimate for Strong-form RFM}, there exists a constant $C>0$, independent of $u^*$ and $u_N$, such that
\[
\inf_{v \in V} \|u^\ast+v-u_N\|_{H^2(\Omega)} \le C\sqrt{\mc{L}(u_N)} .
\]
\end{proposition}

\begin{proof}[Proof.]
Set $L_0v:=a^{-1}Lv=v''+({b}/{a})v'+({c}/{a})v$. Since $1/a\in L^\infty(\Omega)$ and $b,c\in L^2(\Omega)$, one has $b/a,c/a\in L^2(\Omega)$. Hence, by~\citet[Theorem~1]{Gnyp2015Fredholm} with $p=2$, $n=0$, $r=2$ and $m=1$, the operator $T_0:=(L_0,B):H^2(\Omega)\to L^2(\Omega)\times\mathbb R^2$ is Fredholm of index zero. Since $a\neq0$ a.e. in $\Omega$, $\ker T_0=V$. As $T_0$ is Fredholm, $\operatorname{Ran}T_0$ is closed, and the induced operator $\widetilde T_0:H^2(\Omega)/V\to\operatorname{Ran}T_0$, $\widetilde T_0([w])=T_0w$, is an isomorphism of Banach spaces. The bounded inverse theorem gives
\[
\inf_{v\in V}\|w-v\|_{H^2(\Omega)}
\le C\left(\|L_0w\|_{L^2(\Omega)}+\|Bw\|_{L^2(\partial\Omega)}\right),
\]
where $\|L_0w\|_{L^2(\Omega)}\le \|1/a\|_{L^\infty(\Omega)}\|Lw\|_{L^2(\Omega)}$. Taking $w=u_N-u^*$ and using $Lu^*=f$, $Bu^*=g$, we complete the proof.
\end{proof}

As a direct consequence of the trace theorem, we get
\begin{lemma} \label{Lemma: Lv and Bv bounded by H2 and H1 norm of v}
For any $v \in H^{2}(\Omega)$, it holds that $\|L v\|_{L^{2}(\Omega)} \lesssim_{R} \sqrt{\Lambda_{1}^{2}+\Lambda_{2}^{2}+\Lambda_{3}^{2}} \|v\|_{H^{2}(\Omega)}$ and $\|B v\|_{L^{2}(\partial\Omega)} \leq \sqrt{\Lambda_{4}^{2}+\Lambda_{5}^{2}} \|v\|_{H^{1}(\partial\Omega)} ,$ where  $\|v\|_{H^{1}(\partial\Omega)}^2 := \|v\|_{L^2(\partial\Omega)}^2 + \|v^{\prime}\|_{L^2(\partial\Omega)}^2$.
\end{lemma}

\subsection{Proofs of ill-conditioning and PUM mitigation}
\label{Subsection: Proofs of ill-conditioning and PUM mitigation}

In this part, we prove the assertions in \S \ref{subsection of main result: Exponential ill conditionality of random feature matrices}.
\begin{proof}[Proof of Theorem \ref{Theorem: decay rate of random feature matrices and their exponential ill conditionality}.]

Note that $\sigma_{m}(\mf{\Phi}) \leq \|\mf{\Phi}\|_{F} /\sqrt{m}$ holds universally for any matrix and any $m$. We first establish the upper bound for $m \geq 3$. By the Courant-Fisher theorem,
\begin{equation*}
\sigma_{m} = \min _{\operatorname{dim} H=2N-m+1} \max _{\boldsymbol{\alpha} \in H} \frac{\|\mf{\Phi} \boldsymbol{\alpha}\|_{2}}{\|\boldsymbol{\alpha}\|_{2}} ,
\end{equation*}
where $H$ is a linear subspace of $\mb{R}^{2N}$. Our strategy is to construct $2N-m+1$ linearly independent vectors $\{\boldsymbol{\alpha}^{r}\}_{r = 1}^{2N-m+1} \subset \mb{R}^{2N}$ such that $\|\boldsymbol{\alpha}^{r}\|_{2} \geq 1$ and $\|\mf{\Phi}\boldsymbol{\alpha}^{r}\|_{2}$ is sufficiently small. By choosing $H = \operatorname{span}\{\boldsymbol{\alpha}^{r}\}$, we can bound the quotient uniformly over $H$, which immediately yields the desired bound on $\sigma_{m}$.

Let $p=\lfloor(2 N+1-m)/2\rfloor$ and $q=\lceil(2 N+1-m)/2\rceil$. To construct $\boldsymbol{\alpha}^{r}$ for $1\leq r \leq p$, we follow the proof of Theorem \ref{Theorem: main theorem of RFM, Gevrey classes, W^(l,p) error norm, extended distribution}, and approximate $\cos (k_{N} x)$ using $\{\cos (k_{j} x)\}_{j=p+1}^{N-1}$ and $\cos (k_{r} x)$. Let $p_j(x)$ denote the degree-$2(N-p-1)$ Maclaurin polynomial of $\cos(k_j x)$.

Set $\alpha_{N}^{r} = -1$, define the index set $I_{r}^{p} = \{r, p+1, p+2, \ldots, N-1\}$, and set $\alpha_{j}^{r} = 0$ for $j \notin I_{r}^{p} \cup \{N\}$. We determine the remaining coefficients by enforcing $\sum_{j\in I_{r}^{p}} \alpha_{j}^{r} p_{j}(x) = p_{N}(x)$. Matching the coefficients of $x^{2i}$ yields the linear system
\begin{equation*}
\sum_{j\in I_{r}^{p}} k_{j}^{2i} \alpha_{j}^{r}  = k_{N}^{2i}, \quad 0\le i \le N-p-1.
\end{equation*}%
By the properties of Vandermonde determinants, the unique solution is
\begin{equation*}
\alpha_{j}^{r} = \prod_{i \in I_{r}^{p}, i\neq j} \frac{k_{i}^{2} - k_{N}^{2}}{k_{i}^{2} - k_{j}^{2}},  \quad j \in I_{r}^{p}.
\end{equation*}
Clearly, $\alpha_{j}^{r} \neq 0$ for all $j \in I_{r}^{p} \cup \{N\}$. Since $\alpha_{r}^{r} \neq 0$ while $\alpha_{j}^{r} = 0$ for all $j \leq p$ with $j \neq r$, the vectors $\{\boldsymbol{\alpha}^{r}\}_{r = 1}^{p}$ are linearly independent.

Similarly, let $q_j(x)$ be the degree-$(2N-2q-1)$ Maclaurin polynomial of $\sin(k_j x)$. We construct $\boldsymbol{\alpha}^{r}$ for $p+1\leq r \leq p+q$ by setting $\alpha^{r}_{2N} = -1$ and $\alpha^{r}_{j} = 0$ for $j-N \notin I_{r-p}^{q} \cup \{N\}$. The remaining coefficients $\alpha^{r}_{j}$ (for $j-N \in I_{r-p}^{q}$) are determined by enforcing $\sum_{j-N \in I_{r-p}^{q}} \alpha_{j}^{r} q_{j-N}(x) = q_{N}(x)$. Solving this system yields
\[
\alpha^r_{j+N}=\frac{k_N}{k_j}\prod_{\substack{i\in I_{r-p}^q, i\ne j}}\frac{k_i^2-k_N^2}{k_i^2-k_j^2},\qquad j\in I_{r-p}^q .
\]
which ensures \(\alpha^r_{j+N}\ne0\) for all \(j\in I_{r-p}^q\cup\{N\}\). Consequently, the vectors $\{\boldsymbol{\alpha}^{r}\}_{r = p+1}^{p+q}$ are linearly independent. Furthermore, $\{\boldsymbol{\alpha}^{r}\}_{r = 1}^{p+q}$ are linearly independent because the last $N$ components of $\{\boldsymbol{\alpha}^{r}\}_{r = 1}^{p}$ and the first $N$ components of $\{\boldsymbol{\alpha}^{r}\}_{r = p+1}^{p+q}$ are identically zero.

Now, consider any $\boldsymbol{\beta} \in H := \operatorname{span}\{\boldsymbol{\alpha}^{r}\}_{r=1}^{p+q}$. We  decompose $\boldsymbol{\beta} = \boldsymbol{\beta}^{p} + \boldsymbol{\beta}^{q}$, where $\boldsymbol{\beta}^{p} = \sum_{r=1}^{p} \mu_{r} \boldsymbol{\alpha}^{r}$ and $\boldsymbol{\beta}^{q} = \sum_{r=p+1}^{p+q} \mu_{r} \boldsymbol{\alpha}^{r}$. By the construction of $\boldsymbol{\beta}^{p}$, we have
\begin{equation*}
\sum_{j=1}^{N} \beta^{p}_{j} \widetilde{L} p_{j}(x) = \sum_{r=1}^{p} \mu_{r} \widetilde{L} \left(\sum_{j=1}^{N}  \alpha_{j}^{r} p_{j}(x)\right) = 0,
\end{equation*}
which implies $\left(\mf{\Phi}\boldsymbol{\beta}^{p}\right)_{i} = \sum_{j=1}^{N} \beta^{p}_{j} \widetilde{L} \left(\cos(k_{j} x_{i})- p_{j}(x_{i})\right)$.

Following a similar argument for $\boldsymbol{\beta}^{q}$, we obtain $\sum_{j=1}^{N} \beta^{q}_{j+N} \widetilde{L} q_{j}(x) = 0$, which yields $\left(\mf{\Phi}\boldsymbol{\beta}^{q}\right)_{i} = \sum_{j=1}^{N} \beta^{q}_{j+N} \widetilde{L} \left(\sin(k_{j} x_{i})- q_{j}(x_{i})\right)$. By Cauchy-Schwarz inequality, %
\begin{equation*}
\left\|\mf{\Phi}\boldsymbol{\beta}\right\|_{2}^{2} \leq \left\|\boldsymbol{\beta}\right\|_{2}^{2} \sum_{i=1}^{n} \sum_{j=1}^{N}\left[ \left|\widetilde{L} \left(\cos(k_{j} x_{i})-p_{j}(x_{i})\right)\right|^{2} + \left|\widetilde{L} \left(\sin(k_{j} x_{i})-q_{j}(x_{i})\right)\right|^{2} \right].
\end{equation*}
For \(m=3\), set \(\tau=SR\); for \(m>3\), set $\tau = SR\,\Gamma(2N-2p-1)^{-1/(2N-2p-2)}.$ For all \(m\ge3\), set $\widetilde{\tau}
= SR\,\Gamma(2N-2q)^{-1/(2N-2q-1)}.$
With these definitions, we obtain
\begin{equation*}
\begin{aligned}
\left|\widetilde{L}\left(\cos(k_{j} x)-p_{j}(x)\right)\right|^{2}
&\leq \begin{cases}
3\left[a^{2}(x) k_{j}^{4} + b^{2}(x) k_{j}^{2} \tau^{2} + c^{2}(x) \tau^{4}\right] \tau^{4(N-p-1)}, & x \in \Omega, \\
2\gamma\left[g_{1}^{2}(x) k_{j}^{2} \tau^{2} + g_{2}^{2}(x) \tau^{4}\right] \tau^{4(N-p-1)}, & x \in \partial\Omega,
\end{cases} \\%
\left|\widetilde{L}\left(\sin(k_{j} x)-q_{j}(x)\right)\right|^{2}
&\leq \begin{cases}
3\left[a^{2}(x) k_{j}^{4} + b^{2}(x) k_{j}^{2} \tilde{\tau}^{2} + c^{2}(x) \tilde{\tau}^{4}\right] \tilde{\tau}^{4(N-q)-2}, & x \in \Omega, \\
2\gamma\left[g_{1}^{2}(x) k_{j}^{2} \tilde{\tau}^{2} + g_{2}^{2}(x) \tilde{\tau}^{4}\right] \tilde{\tau}^{4(N-q)-2}, & x \in \partial\Omega.
\end{cases}
\end{aligned}
\end{equation*}
Since $\min (2 N-2 p-2,2 N-2 q-1)=m-3$ and $r\mapsto \Gamma(r+1)^{-1/r}$ is decreasing, it follows that $\tau\le \hat{\tau}$ and $\tilde{\tau}\le \hat{\tau}$. Combining the above three bounds, we obtain
\begin{equation*} \label{bound for frac|Zbeta|_2^2 |beta|_2^2}
\begin{aligned}
\left\|\mf{\Phi}\boldsymbol{\beta}\right\|_{2}^{2} %
& \leq \left\|\boldsymbol{\beta}\right\|_{2}^{2}\sum_{j=1}^{N}  \left[3A^{2} k_{j}^{4} + (3D^{2} + 2\tilde{\Lambda}_{4}^{2}) k_{j}^{2} \hat{\tau}^{2} + (3C^{2} + 2\tilde{\Lambda}_{5}^{2}) \hat{\tau}^{4}\right] (1+ \hat{\tau}^{2})\hat{\tau}^{2(m-3)} .
\end{aligned}
\end{equation*}
Hence, the minimax principle yields (\ref{ineq: fast decay of sigma_m(Z)}). Moreover, %
\[
\sigma_{m}^{2}
\leq  6 \sum_{j=1}^{N}\left( A^{2} k_{j}^{4}  + C^{2} + D^{2} k_{j}^{2} +  \tilde{\Lambda}_{4}^{2} k_{j}^{2} + \tilde{\Lambda}_{5}^{2} \right) \max(1,\hat{\tau}^{6}) \hat{\tau}^{2(m-3)}.
\]

Next, we seek a lower bound for the condition number. A direct calculation yields
\begin{equation*}
\begin{aligned}
\|\mf{\Phi}\|_{F}^{2} %
& = \sum_{j=1}^{N}\left(A^{2} k_{j}^{4} - 2 B k_{j}^{2} + C^{2}  + D^{2} k_{j}^{2} + \tilde{\Lambda}_{4}^{2} k_{j}^{2} + \tilde{\Lambda}_{5}^{2}  \right) .
\end{aligned}
\end{equation*}
Note that $\sigma_{1}(\mf{\Phi}) \geq \|\mf{\Phi}\|_{F} / \sqrt{M}$. Denote the random variable
\begin{equation*}
\begin{aligned}
\rho %
& = 1 - \frac{2 B\sum_{j=1}^{N} k_{j}^{2}}{\sum_{j=1}^{N}  \left(A^{2} k_{j}^{4}  + C^{2} + D^{2} k_{j}^{2} +  \tilde{\Lambda}_{4}^{2} k_{j}^{2} + \tilde{\Lambda}_{5}^{2}\right) } .
\end{aligned}
\end{equation*}
Then, the lower bound for condition number reads as
\begin{equation*}
\begin{aligned}
\kappa(\mf{\Phi}) %
& \geq \frac{\|\mf{\Phi}\|_{F}}{\sqrt{M}\sigma_{M}} \geq \frac{\sqrt{\rho}}{\sqrt{6M} \max(1,\hat{\tau}^{3}) \hat{\tau}^{M-3} } .
\end{aligned}
\end{equation*}

It remains to find the lower bound of $\rho$. If $B \leq 0$, then $\rho\geq 1$. Otherwise, the lower bound follows by noting that
$A^{2} k_{j}^{2} + (C^{2} + \tilde{\Lambda}_{5}^{2})k_{j}^{-2} + D^{2} + \tilde{\Lambda}_{4}^{2}$
decreases monotonically before attaining its global minimum at $k_{j}^{2} = A^{-1} \sqrt{C^{2} + \tilde{\Lambda}_{5}^{2}}$. Finally, we derive the high-probability bound. Since $k_{j} \overset{i.i.d.}{\sim} \operatorname{Unif}(0,S)$, let $\hat{k} \sim \operatorname{Unif}(0,1)$, for $\delta > 1, \lambda > 0$,
\begin{equation*}
\begin{aligned}
\mb{P}\left\{\sum_{j=1}^{N}\!\left(\!- A^{2} k_{j}^{4} - C^{2} + 2 \delta B k_{j}^{2} \right) \!\geq 0 \right\} & = \mb{P}\left\{\sum_{j=1}^{N}\!\left(\!- \left(\frac{k_{j}}{S}\right)^{4} - \frac{C^{2}}{A^{2}S^{4}}\! +  \frac{2 \delta B}{A^{2}S^{2}} \left(\frac{k_{j}}{S}\right)^{2} \right)\! \geq 0 \right\} \\
& \leq \left[\mb{E} e^{- \lambda \hat{k}^{4} - \lambda C^{2}A^{-2}S^{-4} + 2 \lambda \delta CA^{-1}S^{-2}\hat{k}^{2}}\right]^{N}.
\end{aligned}
\end{equation*}
Using Lemma \ref{Lemma: integral to bound the expectation} and taking $\lambda = 1/[8(\delta^{2}-1)]$, we control the expectation term by
\begin{equation*}
\begin{aligned}
\max_{\alpha \in \mb{R}} \int_{0}^{1} e^{- \lambda \hat{k}^{4} - \lambda \alpha^{2} + 2 \lambda \delta \alpha \hat{k}^{2}} \mr{d} \hat{k}
& \leq  \tilde{c} \left(\delta^{2}-1\right)^{1/4} ,
\end{aligned}
\end{equation*}
where $\tilde{c}\le 2.994$. %
Setting $\delta = 200/199$ yields $\tilde{c} \sqrt[4]{\delta^{2}-1} \leq 0.95$. Consequently, $\mb{P}(\rho \leq 1-1/\delta)\leq 0.95^{N}, $ which completes the proof.
\end{proof}

\begin{proof}[Proof of Theorem \ref{Theorem: The effect of PUM on the singular value of RFMtx}.]
The first inequality in \eqref{ineq: basic lower and upper bounds for RFMtx with PUM} follows from the fact that the tall matrix $\mf{B}$ may be obtained by deleting the rows of $\mf{\Phi}$ where $\mf{A}_{j}$ and $\mf{C}_{j}$ lie. By column permutation, we partition $\mf{\Phi}$ into two submatrices based on the parity of the patch indices: $\mf{\Phi}_{o} = [\mf{\Phi}_{j}]_{j \text{ is odd}}$ and $\mf{\Phi}_{e} = [\mf{\Phi}_{j}]_{j \text{ is even}}.$ It suffices to upper bound the singular values of $[\mf{\Phi}_{o}, \mf{\Phi}_{e}].$ By Weyl's inequality, $\sigma_{m}^{2}(\mf{\Phi}) \leq  \sigma_{k}^{2}(\mf{\Phi}_{o}) + \sigma_{m+1-k}^{2}(\mf{\Phi}_{e})$ for any $k$ such that $1\leq k \leq N_{o} + 1$ and $1\leq m+1-k \leq N_{e} + 1$. Since there is no overlap between the patches of odd indices, $\mf{\Phi}_{o}$ and $\mf{D}_{o}$ have exactly the same singular values, and similarly, $\mf{\Phi}_{e}$ and $\mf{D}_{e}$ have exactly the same singular values, which imply the second inequality in \eqref{ineq: basic lower and upper bounds for RFMtx with PUM}.  The proof is completed by noting that the singular values of $\mf{B}$ are the union of the singular values of its diagonal blocks $\mf{B}_j$, and the same logic applies to $\mf{D}_{o}$ and $\mf{D}_{e}$.
\end{proof}

\section{Conclusion and outlook} \label{Section: Conclusion and outlook}
We estimate the one-dimensional approximation error of RFM for target functions in the Gevrey and Sobolev classes and apply these estimates to general second-order elliptic PDEs and eigenvalue problems. We prove spectral convergence with respect to the number of features and establish the patch-size convergence rate, highlighting the significant impact of incorporating PUM into the RFM. Additionally, we demonstrate the severe ill-conditioning of the RFMtx and clarify that PUM mitigates the rapid singular-value decay of the RFMtx.

A standard RF expansion employs basis functions of the form $h_i(x) = \eta(k_i x)$. %
While this paper focuses on trigonometric functions for $\eta$, several of our results extend to analytic activation functions. If $\eta$ has nonzero derivatives of all orders at $x=0$ and a convergent Taylor series in $\{x\in\mb{R}:|x|< R \max(C_{u}, \mu)\}$, a similar decomposition as in (\ref{preliminary error decomposition, before pf of Theorem: main theorem of RFM, Gevrey classes, exponential convergence}) holds and the coefficient vector $\boldsymbol{\alpha}$ solves the linear system $\mf{V}_N\boldsymbol{\alpha} = \mf{F}$, where $\mf{V}_N = (k_{j}^{i-1})_{1\le i,j\le N}$ and $F_{i} = u^{(i-1)}(0)/\eta^{(i-1)}(0)$. The proof of Theorem \ref{Theorem: main theorem of RFM, Gevrey classes, W^(l,p) error norm, extended distribution} proceeds similarly. Corollary \ref{Corollary: spectral convergence when solution u in Analytic class, Gevrey class s>1 or Sobolev space, W^(l,p) error norm} requires the Taylor series of $\eta$ to have a sufficiently large radius of convergence. The singular-value analysis in Section~\ref{subsection of main result: Exponential ill conditionality of random feature matrices} treats one feature as the target to be approximated by others. Hence, the argument extends to analytic activation functions.

Several open questions deserve further investigation.

\textbf{(1)} Spectral convergence of the RFM has been numerically observed in higher dimensions~\citep{Chen2022BridgingTA,wang2024extreme}, particularly for $d=2, 3$. We conjecture that the $d$-dimensional convergence rates scale with $N^{1/d}$ instead of $N$, yielding rates of the form $e^{-\kappa N^{1/(sd)}}$ for Gevrey classes and $N^{-s/d}$ for Sobolev spaces with smoothness index $s$.

\textbf{(2)} Leveraging a priori knowledge of the solution enables the design of adaptive features and optimal sampling distributions, which improves RFM performance. For example, optimization-based approaches yield high-quality, transferable feature spaces~\citep{zhang2024transferable}, while variable-depth networks have proven effective in resolving sharp gradients~\citep{dang2024adaptive}.

\textbf{(3)} Effective preconditioners for the RFMtx are crucial for improving the efficiency and accuracy of the associated least-squares solvers.
Recently, randomized preconditioning techniques have emerged as a promising solution to this challenge~\citep{chen2024highprecisionrim,shang2024overlappreconditioner}.

\appendix
\renewcommand{\theHsection}{appendix.\Alph{section}}
\renewcommand{\theHsubsection}{appendix.\Alph{section}.\arabic{subsection}}
\renewcommand{\theHsubsubsection}{appendix.\Alph{section}.\arabic{subsection}.\arabic{subsubsection}}
\renewcommand{\theHequation}{appendix.\Alph{section}.\arabic{equation}}
\renewcommand{\theHtheorem}{appendix.\Alph{section}.\arabic{theorem}}

\section{Bounds for the Gamma function and sampling constants}  \label{Appendix section: Bounds for the Gamma function and sampling constants}
In this part, we prove a two-sided bound for the $\Gamma$-function (Lemma \ref{Lemma: sqrt(2 pi x)(x/e)^x Gamma(x+1)^(-1) is increasing}) and establish Lemma \ref{Lemma: bound for varpi when k_1 follows a Student's t distribution}, which estimates the distribution parameter $\varpi$ of the Student's t-distribution. We recall a useful two-sided bound for the $\Gamma$-function from \citet[Chapter~3]{Artin1964gamma}:
\begin{equation} \label{ineq: two sided Stirling's formula for the Gamma function}
\sqrt{2\pi}x^{x-1/2}e^{-x} \le \Gamma(x)\le \sqrt{2\pi}x^{x-1/2}e^{-x+1/(12x)}\quad \text{for all $x>0$.}
\end{equation}
\begin{lemma} \label{Lemma: sqrt(2 pi x)(x/e)^x Gamma(x+1)^(-1) is increasing}
The function $\sqrt{2 \pi x} \left(x/e\right)^{x}  \Gamma(x+1)^{-1}$ is increasing for $x>0$ and
\begin{equation} \label{ineq: two sided bounds for factorial, Stirling's formula}
\begin{aligned}
\sqrt{2 \pi}e^{-1} \le \sqrt{2 \pi x} \left(x/e\right)^{x}  \Gamma(x+1)^{-1} \le 1 \quad \text{ for } x\geq 1.
\end{aligned}
\end{equation}
\end{lemma}

\begin{proof}[Proof.]
To prove the monotonicity of $f(x)=\sqrt{2 \pi x} \left(x/e\right)^{x} \Gamma(x+1)^{-1}$, we consider its logarithm $h(x) =\ln f(x)$. By $\Gamma(x+1) = x\Gamma(x)$ and differentiating $h(x)$, we obtain $h'(x) = \ln x - 1/{2x} - \psi(x),$ where $\psi(x) = \Gamma'(x)/\Gamma(x)$ is the digamma function. The monotonicity follows directly from \citet[eq.~(2.2)]{Alzer1997gamma}. Then, $f(x) \ge f(1)$ for $x\ge1$. It follows from  \eqref{ineq: two sided Stirling's formula for the Gamma function} that $\lim_{x \to \infty} f(x) = 1$, which completes the proof.
\end{proof}

\begin{lemma} \label{Lemma: bound for varpi when k_1 follows a Student's t distribution}
For $\nu > 2$, the function $f(\nu)$ is strictly decreasing, while the function $g(\nu)$ is strictly increasing, where
\[
    f(\nu) = \frac{2\Gamma((\nu+1)/2)}{\sqrt{(\nu-2)\pi}\, \Gamma(\nu/2)}, \quad\text{and}\quad g(\nu) = 3\sqrt{3} f(\nu) [1 + 3/(\nu-2)]^{-\frac{\nu+1}{2}}.
\]
\end{lemma}

\begin{proof}[Proof.]
To prove the monotonicity of $f(\nu)$, we employ the representation~\citep[eq. (2.10)]{Artin1964gamma} of the digamma function $\psi(x) = \Gamma'(x)/\Gamma(x)$ as
\begin{equation} \label{eq: series representation of digamma function}
\psi(x) = - C - \frac{1}{x} + \sum_{i=1}^\infty \left( \frac{1}{i} - \frac{1}{x+i} \right),
\end{equation}
where $C$ is Euler's constant. Let $L_f(\nu) = \ln f(\nu)$. Its derivative is
\[
    L_f'(\nu) = \frac{1}{2} \left[ \psi\left(\frac{\nu+1}{2}\right) - \psi\left(\frac{\nu}{2}\right) - \frac{1}{\nu-2} \right].
\]
Using \eqref{eq: series representation of digamma function} and a telescoping series,
\begin{equation*}
\begin{aligned}
\psi\left(\frac{\nu+1}{2}\right) - \psi\left(\frac{\nu}{2}\right) & = \sum_{k=0}^\infty \frac{2}{(\nu+2k)(\nu+2k+1)} \\
& < \sum_{k=0}^\infty \frac{2}{(\nu+2k-1)(\nu+2k+1)} = \frac{1}{\nu-1},
\end{aligned}
\end{equation*}
which implies $L_f'(\nu)<0$.

For the monotonicity of $g(\nu)$, let $h(\nu) = 2(\ln g(\nu))^{\prime}$, where
\[
    h(\nu) = \psi\left(\frac{\nu+1}{2}\right) - \psi\left(\frac{\nu}{2}\right) - \ln\left(\frac{\nu+1}{\nu-2}\right) + \frac{2}{\nu-2}.
\]
Since $\lim_{\nu \to \infty} h(\nu) = 0$, it remains to prove $h'(\nu) < 0$ with
\begin{equation*}
\begin{aligned}
h'(\nu) & = \frac{1}{2} \sum_{k=0}^\infty \left[ \frac{1}{[(\nu+1)/{2}+k]^2} - \frac{1}{(\nu/{2}+k)^2} \right] + \frac{\nu-8}{(\nu+1)(\nu-2)^2} \\
& \le -\sum_{k=0}^\infty \frac{1}{[(\nu+1)/{2}+k]^3} + \frac{\nu-8}{(\nu+1)(\nu-2)^2} ,
\end{aligned}
\end{equation*}
where we have used the mean value theorem. Finally, we bound the summation by the integral $\int_0^\infty [(\nu+1)/{2}+t]^{-3} \mr{d}t$, which implies  $h'(\nu) < 0$. This completes the proof.
\end{proof}
\section{Proof of Lemma~\ref{Lemma: probability bound for events A_i^c under general distributions} and Lemma~\ref{Lemma: expectation bound for k_1^(2 N) (mcQ)^(N-1), extended distribution}}
\label{Appendix section: Probabilistic estimates for random feature approximation}
\begin{proof}[Proof of Lemma \ref{Lemma: probability bound for events A_i^c under general distributions}.]
Let $X$ be a random variable with the same distribution as $|k_{1}|/\mu$ and
\[
Y := \frac{X^{2} + w^{2}}{X^{2} - k^{2}}, \quad\text{for some fixed $w, k> 0$.}
\]
Let $\{Y_j\}_{j=1}^{N-1}$ be independent copies of $Y$. We claim that for any $c>2$,
\begin{equation} \label{ineq: rescaled version of Lemma: probability bound for events A_i^c under general distributions}
\begin{aligned}
\mb{P}\left(\prod_{j=1}^{N-1} \left|Y_{j}\right| \ge \left[ e^{c} \mc{Q}(k;w)\right]^{N-1}\right)
\le  e^{-[c-1-\ln (1+c)](N-1)}.
\end{aligned}
\end{equation}
We prove the claim by the moment method. For any $t\in\mb{R}$ and $\lambda \in (1/2, 1)$,
\begin{equation} \label{ineq: basic preprocessing by the moment method}
\begin{aligned}
\mb{P}\left(\prod_{j=1}^{N-1} \left|Y_{j}\right| \geq e^{(N-1) t}\right) \leq e^{-\lambda (N-1) t} \mb{E} \left(\prod_{j=1}^{N-1}  \left|Y_{j}\right|^{\lambda}\right)
=  \left(e^{-\lambda t}\mb{E}|Y|^{\lambda}\right)^{N-1} .
\end{aligned}
\end{equation}
where we may take $t = \ln\mc{Q}(k;w) + c$. We observe that $Y$ becomes very large as $X$ approaches $k$. To estimate $\mb{E}|Y|^{\lambda}$, we decompose it into three parts
\begin{equation*}
\begin{aligned}
\mb{E}|Y|^{\lambda} & = \mb{E}|Y|^{\lambda}\mf{1}_{\{X\le\eta k\}} + \mb{E}|Y|^{\lambda}\mf{1}_{\{\eta k<X\le(2-\eta)k\}} + \mb{E}|Y|^{\lambda}\mf{1}_{\{X>(2-\eta)k\}}, \\
\end{aligned}
\end{equation*}
where we set $\eta=3/4$. We estimate expectations in two cases, $k\ge1/\eta$ and $k<1/\eta$.

When $k\ge1/\eta$, for the first part, since $x\le\eta k$ implies $k^2-x^2\ge (1-\eta^2)k^2$, we have
\begin{equation*}
\begin{aligned}
\mb{E}|Y|^{\lambda}\mf{1}_{\{1<X\le\eta k\}} & \le  \frac{\varpi}{(1-\eta^2)^\lambda k^{2\lambda}} \int_{1}^{\eta k} (x^2 + w^2)^\lambda x^{-3}\mr{d}x \\
& \le  \frac{\varpi}{(1-\eta^2)^\lambda k^{2\lambda}} \left[\frac{1-(\eta k)^{2\lambda-2}}{2-2\lambda} + \frac{w^{2\lambda}}{2} -\frac{w^{2\lambda}}{2\eta^{2} k^{2}}\right] . \\
\end{aligned}
\end{equation*}
It follows from Jensen's inequality and the convexity of $(x^2+w^2)/(k^2-x^2)$ on $x<k$ that
\begin{equation*}
\begin{aligned}
\mb{E}|Y|^{\lambda}\mf{1}_{\{X\le1\}} & \le \left(\varpi\int_0^1 \frac{x^2 + w^2}{k^2 - x^2}\mr{d}x\right)^{\lambda}
\le \frac{\varpi^\lambda}{2^\lambda} \left(\frac{1 + w^2}{k^2 - 1}+\frac{w^2}{k^2}\right)^{\lambda}   . \\
\end{aligned}
\end{equation*}
For the second part, since $(x^2 + w^2)^{\lambda}(x + k)^{-\lambda}x^{-3}$ is convex and decreasing,
\begin{equation*}
\begin{aligned}
\mb{E}|Y|^{\lambda}\mf{1}_{\{\eta k<X\le(2-\eta)k\}} & \le \varpi \int_{\eta k}^{k}  \left(\frac{(x^2 + w^2)^\lambda}{(x+k)^\lambda x^3} + \frac{|(2k-x)^2 + w^2|^\lambda}{(3k-x)^\lambda (2k-x)^3}\right)\frac{\mr{d}x}{\left|x - k\right|^\lambda}  \\
& \le \frac{\varpi(1-\eta)^{1-\lambda}}{(1-\lambda)k^{2+2\lambda}}\left(\frac{(\eta^2 k^2 + w^2)^\lambda}{\eta^3(1+\eta)^\lambda} + \frac{|(2-\eta)^2k^2 + w^2|^\lambda}{(2-\eta)^3(3-\eta)^\lambda}\right) .
\end{aligned}
\end{equation*}
For the third part, $x\ge(2-\eta)k$ ensures $x^2 - k^2 \ge [(2-\eta)^2-1]x^2/(2-\eta)^2$ and then
\begin{equation*}
\begin{aligned}
\mb{E}|Y|^{\lambda}\mf{1}_{\{X>(2-\eta)k\}} & \le  \frac{(2-\eta)^{2\lambda} \varpi}{[(2-\eta)^2-1]^\lambda}\int_{(2-\eta)k}^{\infty} \frac{(x^2 + w^2)^\lambda}{x^{3+2\lambda}} \mr{d}x \\
& \le \frac{(2-\eta)^{2\lambda-2} \varpi}{[(2-\eta)^2-1]^\lambda}\left[\frac{1}{2 k^2} + \frac{w^{2\lambda}}{2+2\lambda}\frac{(2-\eta)^{-2\lambda}}{k^{2+2\lambda}}\right] .
\end{aligned}
\end{equation*}
A combination of the above four bounds and applying Jensen's inequality yield
\begin{equation*}
\begin{aligned}
\mb{E}|Y|^{\lambda} %
& \le \frac{\varpi^\lambda}{2^\lambda} \left(\frac{1 + w^2}{k^2 - 1}+\frac{w^2}{k^2}\right)^{\lambda}
+ \frac{\varpi}{(1-\lambda) k^{2\lambda}} \left[A(\lambda;\eta) + B(\lambda;\eta)  w^{2\lambda}\right] \\
& \le \frac{\varpi^\lambda}{2^\lambda} \left(\frac{1 + w^2}{k^2 - 1}+\frac{w^2}{k^2}\right)^{\lambda}
+ \frac{\varpi(1+w^2)^\lambda}{(1-\lambda) k^{2\lambda}} \left[A(\lambda;\eta)^{\frac{1}{1-\lambda}}+B(\lambda;\eta)^{\frac{1}{1-\lambda}}\right]^{1-\lambda} , \\
\end{aligned}
\end{equation*}
where
\begin{equation*}
\begin{aligned}
& A(\lambda;\eta) = \frac{1/\eta-1}{(1-\eta^2)^\lambda} + \frac{(2/\eta-1)^{2\lambda-2}}{(1-\eta)^\lambda(3-\eta)^\lambda}\left(\frac{1-\eta}{2-\eta} + \frac{1-\lambda}{2}\right), \\
& B(\lambda;\eta) = \frac{1/\eta-1}{(1-\eta^2)^\lambda} + \frac{(2/\eta-1)^{-2}}{(1-\eta)^\lambda(3-\eta)^\lambda}\left(\frac{1-\eta}{2-\eta} + \frac{1-\lambda}{2(1+\lambda)}\right).
\end{aligned}
\end{equation*}
With $\eta=3/4,$ the derivative of $[23A(\lambda)/28]^{\frac{1}{1-\lambda}}$ has the same sign as
$$\ln(23A/28)+(1-\lambda)A^{\prime}/A,$$
whose derivative has the same sign as $A^{\prime\prime}A-(A^{\prime})^2.$ Direct calculations verify that $A^{\prime\prime}A-(A^{\prime})^2\le0$ for $1/2\le\lambda\le1$ and $\ln(23A/28)+(1-\lambda)A^{\prime}/A<0$ at $\lambda=1/2$. Thus, $[23A(\lambda)/28]^{\frac{1}{1-\lambda}}$ is monotonically decreasing. A similar argument applies to $B(\lambda)$ with $B^{\prime\prime}B-(B^{\prime})^2\ge0$ for $1/2\le\lambda\le1$ and $\ln(23B/28)<0$ at $\lambda=1$, which proves that $[23B(\lambda)/28]^{\frac{1}{1-\lambda}}$ is monotonically decreasing. Therefore, $$A(\lambda)^{\frac{1}{1-\lambda}}+B(\lambda)^{\frac{1}{1-\lambda}} \le (28/23)^{\frac{2\lambda-1}{1-\lambda}}\left[A(1/2)^2+B(1/2)^2\right] \le (28/23)^{\frac{\lambda}{1-\lambda}},$$ and %
$$\mb{E}|Y|^{\lambda} \le \frac{2-\lambda}{1-\lambda} \varpi\max\left(\frac{1 + w^2}{2(k^2 - 1)}+\frac{w^2}{2k^2}, \frac{28(1 + w^2)}{23k^{2}}\right)^\lambda .$$
Substituting the above estimate into (\ref{ineq: basic preprocessing by the moment method}) and taking $\lambda = 1-c^{-1}$ prove (\ref{ineq: rescaled version of Lemma: probability bound for events A_i^c under general distributions}) for $k\ge1/\eta$.

When $k\le 1/\eta$, we set $a = \max(k-1/2,0)$ and $b = a+1$. Since $(x^{2} + w^{2})/|x^{2} - k^{2}|$ increases when $x<k$ and decreases when $x>k$, $$\mb{E}|Y|^{\lambda}\mf{1}_{\{X\le a\}} + \mb{E}|Y|^{\lambda}\mf{1}_{\{X>b\}} \le \max\left(\frac{w^2 + a^2}{k^2 - a^2}\mf{1}_{\{k>\frac{1}{2}\}},\frac{w^2 + b^2}{b^2 - k^2}\right)^\lambda.$$ Since $(x^2 + w^2)/(x+k)$ first decreases and then increases, and $(x^2 + w^2)^{\lambda}(x + k)^{-\lambda}x^{-3}$ is decreasing, we split $|Y|^{\lambda}\mf{1}_{\{a<X\le b\}}$ into $|Y|^{\lambda}\mf{1}_{\{a<X\le 1\}} + |Y|^{\lambda}\mf{1}_{\{1<X\le b\}}$ and obtain
\begin{equation*}
\begin{aligned}
\mb{E}|Y|^{\lambda}\mf{1}_{\{a<X\le b\}} %
& \le \frac{2^{\lambda}\varpi}{1-\lambda} \max\left(\frac{w^2+a^2}{k+a}, \frac{w^2+1}{k+1}\right)^\lambda.
\end{aligned}
\end{equation*}
A combination of the above two bounds gives
\begin{equation*}
\begin{aligned}
\mb{E}|Y|^{\lambda} %
& \le\frac{2-\lambda}{1-\lambda} \varpi\cdot 2^{\lambda}\max\left(\frac{w^2 + a^2}{k+a},\frac{w^2 + b^2}{k+b}, \frac{w^2+1}{k+1}\right)^\lambda.
\end{aligned}
\end{equation*}
The claim is proved by substituting the above estimate into (\ref{ineq: basic preprocessing by the moment method}) and taking $\lambda = 1-c^{-1}$.

Let $w = C/\mu$ and recall the independence of $\{k_i\}_{i=1}^{N}$. Conditional on $|k_i|/\mu = k$, the random variable $\prod_{j \neq i} \frac{C^{2}+k_{j}^{2}}{|k_{i}^{2}-k_{j}^{2}|}$ has the same distribution as $\prod_{j=1}^{N-1} |Y_{j}|$. Therefore, by (\ref{ineq: rescaled version of Lemma: probability bound for events A_i^c under general distributions}),
\begin{equation*}
\begin{aligned}
\mb{P}\left(A_{i}^{c}\right) = \mb{E}\left(\mb{E}\left[\mf{1}_{A_{i}^{c}} \mid |k_i|/\mu = k\right]\right) \le e^{-[c-1-\ln (1+c)](N-1)},
\end{aligned}
\end{equation*}
which completes the proof.
\end{proof}

\begin{proof}[Proof of Lemma \ref{Lemma: expectation bound for k_1^(2 N) (mcQ)^(N-1), extended distribution}.]
Denote $k=|k_{1}|/\mu$ and $w=C/\mu$. Since $\mb{E}k^{2} = 1$, we have
\begin{equation*}
\begin{aligned}
\mb{E}\left[k_{1}^{2 N} \mc{Q}(|k_{1}|/\mu;C/\mu)^{N-1} \right]
\le \mu^{2 N} \left[\sup_{k> 0}k^{2} \mc{Q}(k;w)\right]^{N-1}    .
\end{aligned}
\end{equation*}
For $k>4/3$, since $k^2/(k^2-1)$ decreases with respect to $k$,
\begin{equation*}
\begin{aligned}
k^2 \max\left(\frac{1 + w^2}{2(k^2 - 1)}+\frac{w^2}{2k^2}, \frac{28(1 + w^2)}{23k^{2}}\right) %
& \le \max\left(\frac{23w^2}{14}+\frac{8}{7}, \frac{28(1 + w^2)}{23}\right). \\
\end{aligned}
\end{equation*}
For $k\le4/3$, recall $a = \max(k-1/2,0)$ and $b = a+1$. Since $k^2/(k + a)$, $k^2a^2/(k + a)$, $k^2/(k + b)$ and $k^2b^2/(k + b)$ are all non-decreasing in $k$, we have
\begin{equation*}
\begin{aligned}
k^{2}\max\left(\frac{w^2 + a^2}{k+a},\frac{w^2 + 1}{k+1},\frac{w^2 + b^2}{k+b}\right) & = 8\max\left(\frac{4w^2}{39}+\frac{25}{351},\frac{2w^2}{21}+\frac{1}{21},\frac{4w^2}{57}+\frac{121}{513}\right) .
\end{aligned}
\end{equation*}
A combination of the above two bounds yields
\begin{equation*}
\sup_{k> 0}k^{2} \mc{Q}(k;w) %
\le \varpi^{\frac{c}{c-1}}\max\left(\frac{23w^2}{14}+\frac{8}{7}, \frac{64w^2}{57}+\frac{1936}{513}\right),
\end{equation*}
which implies the desired estimate.
\end{proof}
\section{Properties, Cutoffs, and Fourier Decay of Gevrey Functions}
\label{Appendix section: Properties, Cutoffs, and Fourier Decay of Gevrey Functions}
This section provides properties of $G^s$ functions. We work in $d$ dimensions because these estimates are of independent interest.
\subsection{Basic calculus and examples} \label{Appendix section: Basic calculus and examples}

We prove quantitative estimates for the Gevrey parameters associated with Lemma \ref{Lemma: preliminary facts about G^s(Omega)}. For $x \in \mathbb{R}^d$ and $\alpha \in \mathbb{N}^d$, we adopt the notation $x^{\alpha} = \prod_{j=1}^d x_{j}^{\alpha_j}$ and let $|x|_p$ be the $\ell^p$-norm of $x$.
\begin{lemma} \label{Lemma: G^s is a linear space and is closed to the product, dimension d}
Let $s \ge 0$ and $\Omega \subset \mathbb{R}^d$. For $f, g \in G^{s}(\Omega)$ and any compact subset $K \subset \Omega$, \\
(1) \textbf{Linearity:} For any $a, b \in \mathbb{R}$, $af + bg \in G^s(\Omega)$ with $M_{af+bg, K} \le |a|M_{f,K} + |b|M_{g,K}$ and $C_{af+bg, K} \le \max(C_{f,K}, C_{g,K});$ \\
(2) \textbf{Multiplication:} $fg \in G^s(\Omega)$ with $M_{fg, K} \le M_{f,K} M_{g,K}$ and $C_{fg, K} \le C_{f,K} + C_{g,K}$; \\
(3) \textbf{Differentiation:} For any $1\le j \le d$, $\partial_{x_j} f \in G^s(\Omega)$ with $M_{\partial_j f, K} \le M_{f,K} C_{f,K}$ and $C_{\partial_j f, K} \le 2^s C_{f,K}$.
\end{lemma}

Lemma~\ref{Lemma: G^s is a linear space and is closed to the product, dimension d} follows by tracking the constants in the proof of \citet[Proposition~1.4.5]{rodino1993linear}; we omit the details.

The next lemma gives several examples in $G^{0}(\bar{\Omega})$, which can be proved by direct calculations.
\begin{lemma} \label{Lemma: trigonometric functions, exponential functions, polynomials in G^0, dimension d}
Let $w \in \mathbb{R}^d$, $\kappa \in \mathbb{N}^d$ and $\Omega = (-R, R)^d$ with $R > 0$. The following functions belong to the Gevrey class $G^0(\bar{\Omega})$. \\
(1) \textbf{Trigonometric functions} $f(x) = \cos(w^{\top} x)$ or $\sin(w^{\top} x)$: $M_{f, \bar{\Omega}} \le 1$ and $C_{f, \bar{\Omega}} \le |w|_\infty$; \\
(2) \textbf{Exponential function} $f(x) = e^{w^{\top} x}$: $M_{f, \bar{\Omega}} \le e^{R |w|_1}$ and $C_{f, \bar{\Omega}} \le |w|_\infty$; \\
(3) \textbf{Monomials} $f(x) = x^\kappa$: $M_{f, \bar{\Omega}} \le R^{|\kappa|}$ and $C_{f, \bar{\Omega}} \le |\kappa|_\infty/R$; \\
(4) \textbf{Bandlimited functions:} If $\operatorname{supp}\hat{f} \subset [-S, S]^d$ and $\hat{f} \in L^1(\mathbb{R}^d)$, then
$$M_{f, \bar{\Omega}} \le (2\pi)^{-d/2} \|\hat{f}\|_{L^1(\mathbb{R}^d)} \quad \text{and} \quad C_{f, \bar{\Omega}} \le S.$$
\end{lemma}

\subsection{Gevrey cutoff construction via Hörmander mollifiers}
In this subsection, we construct a cutoff function $\varphi \in G^s_0(\mb{R}^d)$ for $s>1$ with explicit Gevrey constants. Our construction relies on the optimal sequence of mollifiers introduced by~\cite{hormander2003analysis}. This cutoff function is used to extend $G^s$ functions and to derive a sharp subexponential Fourier decay estimate.
\begin{lemma}[\!{\citealp[Theorem~1.3.5]{hormander2003analysis}}] \label{Lemma: Hormander's mollifier}
Let $\{\rho_j\}_{j=0}^{\infty}$ be a positive non-increasing sequence with $\rho=\sum_{j=0}^{\infty}\rho_j<\infty$. There exists a non-negative function $\psi \in C_0^\infty([-\rho/2, \rho/2])$ satisfying  %
\begin{equation*}
    \text{$\int_{\mathbb{R}} \psi(x) \, \mr{d}x = 1$ \quad and} \quad |\psi^{(n)}(x)| \le \frac{2^n}{\prod_{j=0}^{n} \rho_j}, \quad \text{for all $x\in\mb{R}$ and $n \ge 0.$}
\end{equation*}
\end{lemma}

\begin{proposition} \label{Proposition: cutoff functions used to extend G^s([-R, R]^d) functions to G^s(mbR^d)}
Let $s > 1$, $R\ge0$, $\rho > 0$ and $\zeta_s = \sum_{k=1}^\infty k^{-s}$. %
There exists a non-negative cutoff function $\varphi\in G^s_0(\mb{R}^d)$ such that $\varphi(x) = 1$ on $[-R, R]^d$ and $\operatorname{supp} \varphi \subset [-R-\rho, R+\rho]^d$ with Gevrey constants $M_{\varphi} \le 1$ and $C_{\varphi} \le 2 \zeta_s/\rho.$
\end{proposition}

\begin{proof}[Proof.]
We construct a mollifier $\psi$ given by Lemma \ref{Lemma: Hormander's mollifier} with the Gevrey sequence $\rho_j := \rho\zeta_s^{-1} (j+1)^{-s}$. Then, we define the one-dimensional cutoff function $\varphi_1$ by the convolution $\varphi_1 := \mf{1}_{[-R-\rho/2, R+\rho/2]} * \psi,$ and define $\varphi(x) := \prod_{i=1}^d \varphi_1(x_i)$. It is easy to verify that $0\le\varphi(x)\le1$, $\varphi(x) = 1$ on $[-R, R]^d$ and $\operatorname{supp} \varphi \subset [-R-\rho, R+\rho]^d$. It remains to demonstrate
\begin{equation*}
    \sup_{x \in \mathbb{R}^d} |D^\alpha \varphi(x)| \le \left(2 \zeta_s/\rho\right)^{|\alpha|} (\alpha!)^s, \quad \text{for all $\alpha \in \mathbb{N}^d$}.
\end{equation*}
It suffices to prove the one-dimensional case, as the multidimensional bound follows from the product rule. For $\alpha=0$, $|\varphi_1| \le 1$. For $\alpha \ge 1$, the derivative is given by
\[
    \varphi_1^{(\alpha)}(x) %
    = \psi^{(\alpha-1)}(x + R + \rho/2) - \psi^{(\alpha-1)}(x - R - \rho/2).
\]
It follows from the derivative bound in Lemma \ref{Lemma: Hormander's mollifier} that
\[
    |\varphi_1^{(\alpha)}(x)| \le 2 \sup_{x\in\mb{R}} |\psi^{(\alpha-1)}(x)| \le \frac{2 \cdot 2^{\alpha-1}}{\prod_{j=0}^{\alpha-1} \rho_j} ,
\]
which yields the desired estimates.
\end{proof}

\subsection{Fourier decay of compactly supported Gevrey functions}
To establish sharp error bounds with explicit dependence on all key parameters, we prove an explicit subexponential Fourier decay estimate for $G_0^s(\mb{R}^d)$-functions in Lemma \ref{Lemma: explicit Fourier decay for G_0^s functions}.
\begin{lemma}
\label{Lemma: explicit Fourier decay for G_0^s functions}
Let $s > 1$ and $\varphi \in G_0^s(\mb{R}^d)$ be supported in a compact set $K$ of volume $V$. Then, $|\hat{\varphi}(\xi)| \le  e^{7s/4} V M_{\varphi,K} e^{-\kappa|\xi|^{1/s} }$ for all $\xi \in \mb{R}^d$, where $\kappa = (3s/{4})d^{-1/(2s)} C_{\varphi,K}^{-1/s}.$
\end{lemma}

\begin{proof}[Proof.]\;\;
Using the identity $\xi^\alpha \hat{\varphi}(\xi) = (-i)^{|\alpha|} \widehat{D^\alpha \varphi}(\xi)$, we have
\[
    |\xi|^N |\hat{\varphi}(\xi)| \le d^{N/2} \max_{|\alpha|=N} |\xi^\alpha \hat{\varphi}(\xi)| = d^{N/2} \max_{|\alpha|=N} |\widehat{D^\alpha \varphi}(\xi)|.
\]
Since $\operatorname{supp} \varphi \subset K$ with volume $V$, Hausdorff-Young inequality ensures
\[
    |\widehat{D^\alpha \varphi}(\xi)| \le \|D^\alpha \varphi\|_{L^1(\mb{R}^d)} \le V \|D^\alpha \varphi\|_{L^{\infty}(K)} .
\]
A combination of the above two estimates and $\varphi \in G_0^s(\mb{R}^d)$ gives
\begin{equation} \label{ineq: bound for Fourier(varphi)(xi), varphi in G_0^s, N to be optimized}
    |\hat{\varphi}(\xi)| \le d^{N/2} V M_{\varphi} \left(C_{\varphi}/|\xi|\right)^N (N!)^s.
\end{equation}
Denote $N_* = C_{\varphi}^{-1/s}d^{-1/(2s)}|\xi|^{1/s}$. We take $N=\lfloor N_* \rfloor$ so that $N_*-1\le N\le N_*$. If $N_*<1$, setting $N=0$ in \eqref{ineq: bound for Fourier(varphi)(xi), varphi in G_0^s, N to be optimized} gives the desired estimate. If $N_*\ge1$, we employ \eqref{ineq: two sided bounds for factorial, Stirling's formula} to bound $N!$ and obtain $|\hat{\varphi}(\xi)|\le e^s V M_{\varphi} N^{s/2}e^{-sN}.$ Bounding $N^{s/2}e^{-sN}$ by $e^{-3sN/4}$ completes the proof.
\end{proof}

We extend $G^s$ functions on compact sets to $\mb{R}^d$ and establish the subexponential Fourier decay.
\begin{proposition} \label{Proposition: explicit Fourier decay for extended Gevrey functions}
Let $s, R, \rho, \zeta_s$ and $\varphi$ be defined in Proposition \ref{Proposition: cutoff functions used to extend G^s([-R, R]^d) functions to G^s(mbR^d)}, and $\Omega=[-R-\rho, R+\rho]^d$. Given $u \in G^s(\Omega)$, define the cutoff extension $v$ by $v=\varphi u$ on $\Omega$ and $v=0$ elsewhere. Then, $|\hat{v}(\xi)| \le Ce^{-\kappa|\xi|^{1/s} }$ with $C = 2^d e^{7s/4} (R+\rho)^d M_{u}$ and $\kappa = (3s/{4})d^{-1/(2s)} \left(C_u + 2 \zeta_s/\rho\right)^{-1/s}.$
\end{proposition}

\begin{proof}[Proof.]\;\;
It follows from Proposition \ref{Proposition: cutoff functions used to extend G^s([-R, R]^d) functions to G^s(mbR^d)} and Lemma \ref{Lemma: G^s is a linear space and is closed to the product, dimension d} that the Gevrey constants $M_v \le M_u$ and $C_v \le C_u + 2 \zeta_s/\rho.$ Applying Lemma \ref{Lemma: explicit Fourier decay for G_0^s functions} with these bounds, and noting that the volume $|\operatorname{supp} v|$ is bounded by $2^d (R+\rho)^d$ completes the proof.
\end{proof}

\section{Gevrey approximation of Sobolev functions}
To prove sharp RFM approximation rates for general Sobolev functions in the $W^{\sigma,p}$ norm ($1 \le p \le \infty$) in Corollaries \ref{Corollary: spectral convergence when solution u in Analytic class, Gevrey class s>1 or Sobolev space, W^(l,p) error norm} and \ref{Corollary: convergence for PUM when solution u in Analytic class, Gevrey class s>1 or Sobolev space}, we construct a bandlimited $G^0(\Omega)$ approximator for such functions and prove the explicit approximation error that includes the Gevrey constants, bandwidth, and norms of the approximants. The third estimate in Theorem \ref{Theorem: approximate W^(s,p)-func by G^0 func on Lipschitz domain or whole space with explicit error bound} is not used in this paper; we include it because it is useful in RKHS-based approximation estimates.

\begin{theorem} \label{Theorem: approximate W^(s,p)-func by G^0 func on Lipschitz domain or whole space with explicit error bound}
Let $\Omega = \mb{R}^d$ or $\Omega \subset \mb{R}^d$ be a bounded Lipschitz domain. Let $s \ge 0$, $1 \le p \le \infty$, and $0<\varepsilon\le1$. For any $u \in W^{s,p}(\Omega)$, there exists a function $u_\varepsilon \in G^0(\mb{R}^d)$ such that:

1. \textbf{Approximation Error:} For all $t \in [0, s]$,
    \begin{equation} \label{ineq: Approximation Error Bound for Gevrey approximation}
        \|u - u_\varepsilon\|_{W^{t,p}(\Omega)} \lesssim_{s,p,d,\Omega} \varepsilon^{s-t} \|u\|_{W^{s,p}(\Omega)} .
    \end{equation}

2. \textbf{Gevrey Regularity Estimates:} For all $x\in\mathbb{R}^d$ and $\beta\in\mathbb{N}^d$,
    \begin{equation} \label{ineq: Gevrey Regularity Estimates for Gevrey approximation}
        |D^\beta u_\varepsilon(x)| \lesssim_{s,p,d,\Omega} \left( C_{s}\varepsilon\right)^{\min(\lfloor s\rfloor-|\beta|,0)-d/p} \|u\|_{W^{s,p}(\Omega)}.
    \end{equation}

3. \textbf{Bandlimitedness:} $\operatorname{supp} \hat{u}_\varepsilon \subset B_{1/\varepsilon}(0)$. If, in addition, $\Omega$ is bounded, then  \begin{equation*} %
    \|u_\varepsilon\|_{L^{2}(\mb{R}^d)} \lesssim_{s,p,d,\Omega} \varepsilon^{d\min(1/2-1/p,0)} \|u\|_{W^{s,p}(\Omega)}.
\end{equation*}
\end{theorem}

\begin{proof}[Proof.]
If $s=0$, take $u_\varepsilon=0$ and the two estimates are immediate. Hence we assume $s>0$ below. We start with the case $\Omega = \mb{R}^d$. First, we construct the approximant $u_\varepsilon$ by convolving with a mollifier $\phi_\varepsilon(x) = \varepsilon^{-d} \phi(x/\varepsilon)$. We employ a bandlimited $\phi$ with sub-exponential decay. Specifically, let $\hat{\phi} \in G_0^2(B_1(0))$, normalized such that $\int_{\mb{R}^d} \phi(x) \, \mr{d}x = 1$. By \citet[Theorem~1.6.1]{rodino1993linear}, there exist constants $C_1, C_2$ depending only on $\phi$ such that
\begin{equation} \label{ineq: sub-exponential decay of mollifier phi}
|\phi(x)| \le C_1 e^{-C_2 \sqrt{|x|}}, \quad \text{for all } x \in \mb{R}^d.
\end{equation}
Denote the kernel $K_\alpha(z) := z^\alpha \phi_\varepsilon(z)/(\alpha!)$ for $\alpha\in\mb{N}^d$. The approximant is defined as
\begin{equation*} %
u_\varepsilon(x) = \sum_{|\alpha| \le m_{s}}  \int_{\mb{R}^d} D^\alpha u(y) K_\alpha(x-y) \, \mr{d}y, \quad\text{with}\quad m_{s} := \lceil s\rceil-1.
\end{equation*}
For $p < \infty$, by a standard density argument, it suffices to establish the two estimates for $u \in C_0^\infty(\mb{R}^d)$. For $p = \infty$, we recall the embedding $W^{s,\infty}(\mb{R}^d) \subset C^{\lceil s \rceil - 1}(\mb{R}^d)$.

Second, we estimate the approximation error $R := u - u_\varepsilon$ in $L^p$-norm. When $s\in \mb{N}_+$, we apply Taylor's formula with integral remainder %
to  $f(t) = u(tx + (1-t)y)$ and obtain
\begin{equation*} %
u(x) = \sum_{|\alpha| \le s-1} \frac{1}{\alpha!} D^\alpha u(y) (x-y)^\alpha
+ \sum_{|\alpha|=s} \frac{s(x-y)^\alpha}{\alpha!} \int_0^1 \frac{D^\alpha u(tx + (1-t)y)}{(1-t)^{1-s}} \, \mr{d}t.
\end{equation*}
It follows from the above formula and the normalization property of $\phi_\varepsilon$ that
\begin{equation*}
\begin{aligned}
R(x) & = \sum_{|\alpha|=s} \frac{s}{\alpha!} \int_{\mb{R}^d} \int_0^1 (1-t)^{s-1} D^\alpha u(tx+(1-t)y) (x-y)^\alpha \phi_\varepsilon(x-y) \, \mr{d}t \, \mr{d}y \\
& = \sum_{|\alpha|=s} \frac{s}{\alpha!} \int_0^1 \int_{\mb{R}^d} (1-t)^{-1-d} D^\alpha u(z) (x-z)^\alpha \phi_\varepsilon\left(\frac{x-z}{1-t}\right) \, \mr{d}z \, \mr{d}t ,
\end{aligned}
\end{equation*}
where we exchange the order of integration and perform the change of variables $z = tx + (1-t)y$ in the second line. %
Then, we estimate
\begin{equation*}
\begin{aligned}
|R(x)| & \le \sum_{|\alpha|=s} \frac{s}{\alpha!} \int_{\mb{R}^d} \left|D^\alpha u(z)\right| |x-z|^s \left(\int_0^1 (1-t)^{-1-d}\left|\phi_\varepsilon\left(\frac{x-z}{1-t}\right)\right| \, \mr{d}t\right) \, \mr{d}z.
\end{aligned}
\end{equation*}
Using (\ref{ineq: sub-exponential decay of mollifier phi}), the inner integral may be bounded by %
\begin{equation} \label{ineq: inner integral bound from sub-exponential decay}
\begin{aligned}
\frac{C_1}{|x-z|^d} \int_{|x-z|/\varepsilon}^\infty v^{d-1} e^{-C_2 \sqrt{v}} \, \mr{d}v.
\end{aligned}
\end{equation}
Substituting this bound and applying Young's inequality, we obtain
\begin{equation} \label{ineq: L^p bound for R when s is an integer with integral to be controlled}
\begin{aligned}
\|R\|_{L^p(\mb{R}^d)} & \le \sum_{|\alpha|=s} \frac{sC_1}{\alpha!} \|D^\alpha u\|_{L^p(\mb{R}^d)} \int_{\mb{R}^d}  |x|^{s-d} \int_{|x|/\varepsilon}^\infty v^{d-1} e^{-C_2 \sqrt{v}} \, \mr{d}v \, \mr{d}x .
\end{aligned}
\end{equation}
A straightforward computation shows that for any $s, q>0$,
\begin{equation} \label{ineq: integral estimate for |x|^(sq - d) under subexponentially decaying weights}
\int_{\mathbb{R}^d} |x|^{sq - d} \left( \int_{|x|/\varepsilon}^\infty v^{d-1} e^{-C_2 \sqrt{v}} \, \mathrm{d}v \right)^q \mathrm{d}x
\lesssim_{s,q,d} \varepsilon^{sq}.
\end{equation}
Applying \eqref{ineq: integral estimate for |x|^(sq - d) under subexponentially decaying weights} to \eqref{ineq: L^p bound for R when s is an integer with integral to be controlled} with $q=1$, we get
\begin{equation} \label{ineq: bound approximation error R in L^p-norm}
\begin{aligned}
\|R\|_{L^p(\mb{R}^d)} & \lesssim_{s,d} \varepsilon^{s} \sum_{|\alpha|=s} \|D^\alpha u\|_{L^p(\mb{R}^d)}\quad \text{ for all $1\le p\le\infty$}.
\end{aligned}
\end{equation}

When $s\notin \mb{N}$, our goal is to prove that  for all $1\le p\le\infty$,
\begin{equation}  \label{ineq: bound approximation error R in L^p-norm, s notin mbN}
\begin{aligned}
\|R\|_{L^p(\mb{R}^d)} & \lesssim_{s,p,d} \varepsilon^{s} \sum_{|\alpha|=k} [D^\alpha u]_{W^{\sigma,p}(\mb{R}^d)} .
\end{aligned}
\end{equation}
For $s<1$, it follows from Minkowski's inequality and the standard translation estimate \(\|u(\cdot-h)-u\|_{L^p(\mathbb R^d)}
\le |h|^s [u]_{W^{s,p}(\mathbb R^d)}\)  for all \(h\in\mathbb R^d\) that $\|R\|_{L^p(\mb{R}^d)}$ is bounded by
\begin{equation*}
\begin{aligned}
\int_{\mathbb R^d}|\phi(z)|
\|u(\cdot)-u(\cdot-\varepsilon z)\|_{L^p(\mb{R}^d)}\,dz
\lesssim_{s,d}\varepsilon^s [u]_{W^{s,p}(\mb{R}^d)}
\int_{\mathbb R^d}|\phi(z)|\,|z|^s\,dz.
\end{aligned}
\end{equation*}
For $s>1$, let $s=k+\sigma$ with $k\in \mb{N}$ and $\sigma\in(0,1)$. We similarly apply Taylor's formula with integral remainder and extract one more term to obtain
\begin{equation*} %
\begin{aligned}
u(x) & = \sum_{|\alpha| \le k} \frac{D^\alpha u(y)}{\alpha!}  (x-y)^\alpha + \sum_{|\alpha|=k} \frac{k(x-y)^\alpha}{\alpha!} \int_0^1 \frac{D^\alpha u(tx + (1-t)y) - D^\alpha u(y)}{(1-t)^{1-k}}  \, \mr{d}t.
\end{aligned}
\end{equation*}
It follows from the above formula and the normalization property of $\phi_\varepsilon$ that $R(x)$ may be decomposed into $\operatorname{I}(x) + \operatorname{II}(x)$ with
\begin{equation*}
\begin{aligned}
& \operatorname{I}(x) = \sum_{|\alpha|=k} \frac{k}{\alpha!} \int_{\mb{R}^d} \int_0^1 \frac{D^\alpha u(tx+(1-t)y) - D^\alpha u(x) }{(1-t)^{1-k}} (x-y)^\alpha \phi_\varepsilon(x-y) \, \mr{d}t \, \mr{d}y, \\
& \operatorname{II}(x) = \sum_{|\alpha|=k} \frac{k}{\alpha!} \int_{\mb{R}^d} \int_0^1 \frac{D^\alpha u(x) - D^\alpha u(y)}{(1-t)^{1-k}} (x-y)^\alpha \phi_\varepsilon(x-y) \, \mr{d}t \, \mr{d}y .
\end{aligned}
\end{equation*}
To control $\operatorname{I}$, similarly, we perform the change of variables $z = tx + (1-t)y$ and obtain
\begin{equation*}
\begin{aligned}
|\operatorname{I}(x)| & = \sum_{|\alpha|=k} \frac{k}{\alpha!} \int_{\mb{R}^d} \frac{|D^\alpha u(z) - D^\alpha u(x)|}{|x-z|^{d/p+\sigma}} \left(\int_0^1 \frac{|x-z|^{d/p+s}}{(1-t)^{d+1}} \left|\phi_\varepsilon\left(\frac{x-z}{1-t}\right)\right| \, \mr{d}t\right) \, \mr{d}z \\
& \lesssim_{s} \sum_{|\alpha|=k} \int_{\mb{R}^d} \frac{|D^\alpha u(z) - D^\alpha u(x)|}{|x-z|^{d/p+\sigma}} \left(|x-z|^{s-d/q} \int_{|x-z|/\varepsilon}^\infty v^{d-1} e^{-C_2 \sqrt{v}} \, \mr{d}v \right) \, \mr{d}z, \\
\end{aligned}
\end{equation*}
where we used \eqref{ineq: inner integral bound from sub-exponential decay} in the second line and $1/q = 1-1/p.$ If $q=\infty$, note that
\begin{equation*}
\begin{aligned}
\sup_{x\in\mb{R}^d}|x|^{s} \int_{|x|/\varepsilon}^\infty v^{d-1} e^{-C_2 \sqrt{v}} \, \mr{d}v = \varepsilon^{s}\sup_{r\in\mb{R}}r^{s} \int_{r}^\infty v^{d-1} e^{-C_2 \sqrt{v}}\, \mr{d}v.
\end{aligned}
\end{equation*}
An application of H\"older's inequality, \eqref{ineq: integral estimate for |x|^(sq - d) under subexponentially decaying weights} and the above relation gives
\begin{equation} \label{ineq: bound for |operatornameI(x)| as a part of R(x)}
|\operatorname{I}(x)| \lesssim_{s,p,d}
    \begin{cases}
         \varepsilon^{s} \sum_{|\alpha|=k} \left(\int_{\mb{R}^d} \frac{|D^\alpha u(z) - D^\alpha u(x)|^p}{|x-z|^{d+\sigma p}}\,\mr{d}z\right)^{1/p}, & \text{if } 1 \le p < \infty, \\[10pt]
        \varepsilon^{s} \sum_{|\alpha|=k} \sup_{z\in\mathbb R^d,\ z\ne x} \frac{|D^\alpha u(z)-D^\alpha u(x)|}{|z-x|^\sigma}, & \text{if } p = \infty.
    \end{cases}
\end{equation}
The control of $\operatorname{II}$ proceeds in an analogous way where we integrate with respect to $t$, apply H\"older's inequality and bound the inner integral. Then, we bound $|\operatorname{II}(x)|$ by the right-hand side of \eqref{ineq: bound for |operatornameI(x)| as a part of R(x)} even with $\varepsilon^{s}$ replaced by $\varepsilon^{s+d}$. Therefore, for non-integer $s$, we achieve \eqref{ineq: bound approximation error R in L^p-norm, s notin mbN}.

Next, we estimate the derivatives of $u_\varepsilon$. %
For any multi-index $\beta$,
\[
D^\beta u_\varepsilon(x) = \sum_{|\alpha| \le m_{s}} \int_{\mb{R}^d} D^\beta K_\alpha(x-y) D^\alpha u(y) \, \mr{d}y.
\]
We decompose $\beta = \gamma + \nu$ with $\gamma, \nu \in \mb{N}^d$ such that $|\gamma| = \max(|\alpha|+|\beta|-\lfloor s\rfloor,0)$. Since $\frac{\partial}{\partial x_i}K_\alpha(x-y) = -\frac{\partial}{\partial y_i}K_\alpha(x-y) ,$ integration by parts gives
\begin{equation} \label{eq: representation of D^beta u_varepsilon after integration by parts}
D^\beta u_\varepsilon(x) = \sum_{|\alpha| \le m_{s}} \int_{\mb{R}^d} D^\gamma K_\alpha(x-y) D^{\alpha+\nu} u(y) \, \mr{d}y.
\end{equation}
Recall $\hat{\phi} \in G_0^2(B_1(0))$, which ensures $\widehat{\cdot^\alpha\phi} \in G_0^2(B_1(0))$. Then, by Bernstein Lemma~\citep[Lemma 2.1]{Bahouri2011Fourier}, there exists a constant $C_{s} > 0$ such that for all $|\alpha|\le s$, $\gamma \in \mb{N}^d$ and $1\le r\le \infty$,
\begin{equation*}
\begin{aligned}
\|D^\gamma K_\alpha\|_{L^r(\mb{R}^d)} & = (\alpha!)^{-1}\varepsilon^{|\alpha|-|\gamma|-d+d/r}\|D^\gamma (\cdot^\alpha\phi)\|_{L^r(\mb{R}^d)} \\
& \lesssim_{s} \varepsilon^{|\alpha|-|\gamma|-d+d/r} C_{s}^{|\gamma|}\|\cdot^\alpha\phi\|_{L^1(\mb{R}^d)}, \\
\end{aligned}
\end{equation*}
where \eqref{ineq: sub-exponential decay of mollifier phi} implies $\|\cdot^\alpha\phi\|_{L^1(\mb{R}^d)} \lesssim_{s,d} 1$. Let $\bar{p}$ be such that $p\le\bar{p}\le\infty$. With $1/r = 1+1/\bar{p}-1/p,$ it follows from  Young’s inequality that
\begin{equation} \label{ineq: L^(bar(p)) estimate of D^beta u_varepsilon}
\begin{aligned}
\|D^\beta u_\varepsilon\|_{L^{\bar{p}}(\mb{R}^d)} & \le \sum_{|\alpha| \le m_{s}} \|D^\gamma K_\alpha\|_{L^{r}(\mb{R}^d)} \|D^{\alpha+\nu} u\|_{L^p(\mb{R}^d)} \\
& \lesssim_{s,d} \sum_{|\alpha| \le m_{s}} \varepsilon^{|\alpha|-|\gamma|-d+d/r} C_{s}^{|\gamma|} \|D^{\alpha+\nu} u\|_{L^p(\mb{R}^d)} .
\end{aligned}
\end{equation}
The estimate \eqref{ineq: Gevrey Regularity Estimates for Gevrey approximation} follows directly from \eqref{ineq: L^(bar(p)) estimate of D^beta u_varepsilon} with $\bar{p}=\infty$. For $|\beta|\le \lfloor s\rfloor$, we have $|\alpha|-|\gamma|\ge0$ and $|\gamma|\le s$. With $\bar{p}=p$, \eqref{ineq: L^(bar(p)) estimate of D^beta u_varepsilon} implies the stability estimate
\begin{equation} \label{ineq: stability estimate for integer s}
\begin{aligned}
\|D^\beta u_\varepsilon\|_{L^{p}(\mb{R}^d)} & \lesssim_{s,d} \sum_{|\alpha| \le m_{s}} \|D^{\alpha+\nu} u\|_{L^p(\mb{R}^d)} ,
\end{aligned}
\end{equation}
For non-integer $s$, it follows from \eqref{eq: representation of D^beta u_varepsilon after integration by parts} that
\begin{equation*}
\begin{aligned}
D^\beta u_\varepsilon(x) - D^\beta u_\varepsilon(\bar{x}) = \sum_{|\alpha| \le m_{s}} \int_{\mb{R}^d} D^\gamma K_\alpha(y) [D^{\alpha+\nu} u(x-y) - D^{\alpha+\nu} u(\bar{x}-y)] \, \mr{d}y .
\end{aligned}
\end{equation*}
For $p<\infty$, an application of H\"older's inequality gives
\begin{equation*}
\begin{aligned}
\left|D^\beta u_\varepsilon(x) - D^\beta u_\varepsilon(\bar{x})\right|^p & \lesssim_{s,p} \sum_{|\alpha| \le m_{s}} \|D^\gamma K_\alpha\|_{L^1(\mb{R}^d)}^{p-1}\cdot\\
& \qquad \int_{\mb{R}^d} |D^\gamma K_\alpha(y)| |D^{\alpha+\nu} u(x-y) - D^{\alpha+\nu} u(\bar{x}-y)|^p \, \mr{d}y .
\end{aligned}
\end{equation*}
For $|\beta| = \lfloor s\rfloor$, we have $|\gamma|=|\alpha|\le \lfloor s\rfloor$, which implies $\|D^\gamma K_\alpha\|_{L^1(\mb{R}^d)}\lesssim_{s}1$ and
\begin{equation*}
\begin{aligned}
[D^\beta u_\varepsilon]_{W^{\sigma,p}(\mb{R}^d)} \lesssim_{s,p,d} \sum_{|\alpha| \le m_{s}} [D^{\alpha+\nu} u]_{W^{\sigma,p}(\mb{R}^d)}.
\end{aligned}
\end{equation*}
For $p=\infty$, the above estimate follows analogously. %
A combination of the above estimate with \eqref{ineq: stability estimate for integer s} yields $\|u_\varepsilon\|_{W^{s,p}(\mb{R}^d)} \lesssim_{s,p,d} \|u\|_{W^{s,p}(\mb{R}^d)},$ and then
\begin{equation*}
\begin{aligned}
\|R\|_{W^{s,p}(\mb{R}^d)} \lesssim_{s,p,d} \|u\|_{W^{s,p}(\mb{R}^d)}.
\end{aligned}
\end{equation*}
From \eqref{ineq: bound approximation error R in L^p-norm} and the above bound, by the Gagliardo-Nirenberg interpolation inequality (Theorem \ref{Theorem: Gagliardo–Nirenberg inequalities and non-inequalities: The full story}) between $L^p(\mb{R}^d)$ and $W^{s,p}(\mb{R}^d)$, we proved \eqref{ineq: Approximation Error Bound for Gevrey approximation} when $\Omega=\mb{R}^d.$

Finally, when $\Omega$ is a bounded Lipschitz domain~\citep[Chapters 5 and 7]{Adams2003Sobolev}, there exists a bounded linear extension operator $E: W^{s,p}(\Omega)\to W^{s,p}(\mb{R}^d)$ such that $\operatorname{supp} (Eu) \subset K_{\Omega}$ for some compact set $K_{\Omega}$ depending only on $\Omega$. We construct $u_\varepsilon := (Eu)_\varepsilon$ exactly as in the $\mb{R}^d$ case. The estimates \eqref{ineq: Approximation Error Bound for Gevrey approximation} and \eqref{ineq: Gevrey Regularity Estimates for Gevrey approximation} are then easily verified by $\|u - u_\varepsilon\|_{W^{t,p}(\Omega)} \le \|Eu - u_\varepsilon\|_{W^{t,p}(\mb{R}^d)}$ and $\|Eu\|_{W^{s,p}(\mb{R}^d)} \lesssim_{s,p,\Omega} \|u\|_{W^{s,p}(\Omega)}$. The bandlimitedness of $u_\varepsilon$ is ensured by $\hat{\phi} \in G_0^2(B_1(0))$.  Moreover, by \eqref{ineq: L^(bar(p)) estimate of D^beta u_varepsilon} and the fact that $\operatorname{supp} (Eu) \subset K_{\Omega}$, we have
\begin{equation*}
\|u_\varepsilon\|_{L^{2}(\mb{R}^d)}
\lesssim_{s,d} \sum_{|\alpha| \le m_{s}} \varepsilon^{|\alpha|+d\min(1/2-1/p,0)}\|D^{\alpha} (Eu)\|_{L^{\min(p,2)}(\mb{R}^d)} ,
\end{equation*}
and $\|D^{\alpha} (Eu)\|_{L^{\min(p,2)}(\mb{R}^d)}\lesssim_{s,p,d,\Omega}\|u\|_{W^{s,p}(\Omega)}$ for $|\alpha| \le m_{s}$, which completes the proof.
\end{proof}

To bound the approximation error in fractional Sobolev spaces, we recall the following sharp version of the Gagliardo--Nirenberg interpolation inequality.
\begin{theorem}[\citealp{Brezis2018GagliardoNirenberg}] \label{Theorem: Gagliardo–Nirenberg inequalities and non-inequalities: The full story}
Let $\Omega \subset \mathbb{R}^d$ be either $\mathbb{R}^d$ or a bounded Lipschitz domain. Let $0 \le s_1 \le s \le s_2$, $1 \le p, p_1, p_2 \le \infty$, and $\theta \in (0, 1)$ such that $s = \theta s_1 + (1 - \theta)s_2$ and $1/p = \theta/p_1 + (1 - \theta)/p_2$. Then, for any $u \in W^{s_1, p_1}(\Omega) \cap W^{s_2, p_2}(\Omega)$,
\[ \|u\|_{W^{s,p}(\Omega)} \le C \|u\|_{W^{s_1, p_1}(\Omega)}^{\theta} \|u\|_{W^{s_2, p_2}(\Omega)}^{1-\theta} \]
holds if and only if the following exceptional case is excluded: $s_2 \in \mathbb{N}_+$, $p_2 = 1$, and $0 < s_2 - s_1 \le 1 - 1/p_1$. The constant $C$ depends on $p$, $p_1$, $p_2$, $s$, $s_1$, $s_2$, $\theta$, $\Omega$, but not on $u$.
\end{theorem}

\section{Auxiliary estimates for RFMtx conditioning} \label{Appendix section: Auxiliary result in the proof of Theorem: decay rate of random feature matrices and their exponential ill conditionality}
\begin{lemma} \label{Lemma: integral to bound the expectation}
For $\lambda > 0$ and $\delta > 1$, there holds
\begin{equation*}
\begin{aligned}
\max_{\alpha \in \mb{R}} J(\alpha) := \max_{\alpha \in \mb{R}} \int_{0}^{1} e^{- \lambda \hat{k}^{4} - \lambda \alpha^{2} + 2 \lambda \delta \alpha \hat{k}^{2}} \mr{d} \hat{k} & \leq c  \lambda^{-1/4} e^{2 \lambda\left(\delta^{2}-1\right)},
\end{aligned}
\end{equation*}
where $c = \Gamma\left(5/4\right) \left[2\operatorname{erf}\left(\pi \Gamma\left(5/4\right)^{-2}/16\right) + 1 \right]$ and  $\operatorname{erf}(x) = \frac{2}{\sqrt{\pi}}\int_{0}^{x} e^{-t^{2}} \mr{d} t$.
\end{lemma}
\begin{proof}[Proof.]
Denote $\mr{I}(\alpha) := \int_{0}^{1} e^{- \lambda\left(\hat{k}^{2}-\alpha\right)^{2}}  \mr{d} \hat{k}$ for $\alpha\in\mb{R}$. A direct calculation yields
\begin{equation*}
\begin{aligned}
\frac{\mr{d}}{\mr{d} \alpha} J(\alpha) & = 2\lambda \int_{0}^{1} (\delta\hat{k}^{2} - \alpha) e^{- \lambda \hat{k}^{4} - \lambda \alpha^{2} + 2 \lambda \delta \alpha \hat{k}^{2}} \mr{d} \hat{k},
\end{aligned}
\end{equation*}
which is positive when $\alpha \leq 0$ and negative when $\alpha \geq \delta$, and for $\alpha \geq 1$,
\begin{equation*}
\begin{aligned}
\frac{\mr{d}}{\mr{d} \alpha} \mr{I}(\alpha) & = 2\lambda \int_{0}^{1} (\hat{k}^{2} - \alpha) e^{- \lambda\left(\hat{k}^{2}-\alpha\right)^{2}}  \mr{d} \hat{k} \leq 0.
\end{aligned}
\end{equation*}
The above two facts indicate that $\max_{\alpha \in \mb{R}} J(\alpha)  \leq \max_{0\leq \alpha \leq \delta} J(\alpha)$ and
\begin{equation}
\begin{aligned} \label{ineq: the first ineq to control the integral}
\max_{0\leq \alpha \leq \delta} J(\alpha) & \leq e^{2 \lambda\left(\delta^{2}-1\right)} \max_{0\leq \alpha \leq \delta} \mr{I}(\alpha)
\leq e^{2 \lambda\left(\delta^{2}-1\right)} \max_{0\leq \alpha \leq 1} \mr{I}(\alpha) .
\end{aligned}
\end{equation}
For $\alpha \in [0, 1]$, we decompose $\mr{I}(\alpha)$ into $\mr{I}_{1}(\alpha):= \int_{0}^{\sqrt{\alpha}} e^{- \lambda\left(\hat{k}^{2}-\alpha\right)^{2}}  \mr{d} \hat{k}$ and $\mr{I}_{2}(\alpha) := \mr{I}(\alpha) - \mr{I}_{1}(\alpha).$
$\mr{I}_{1}(\alpha)$ may be bounded by
\begin{equation*}
\begin{aligned}
\mr{I}_{1} & \leq \int_{0}^{\sqrt{\alpha}} e^{- \lambda \alpha\left(\hat{k}-\sqrt{\alpha}\right)^{2}}  \mr{d} \hat{k}
= \frac{1}{2}\sqrt{\frac{\pi}{\lambda\alpha}}
\operatorname{erf}(\sqrt{\lambda} \alpha) .\\
\end{aligned}
\end{equation*}
We control $\mr{I}_{2}(\alpha)$ in two ways as
\begin{equation*}
\begin{aligned}
\mr{I}_{2} & \leq \int_{\sqrt{\alpha}}^{1} e^{- \lambda\left(\hat{k}-\sqrt{\alpha}\right)^{4}}  \mr{d} \hat{k}
\leq \lambda^{-1/4}\int_{0}^{\infty} e^{-x^{4}} \mr{d}x = \lambda^{-1/4}\Gamma\left(5/4\right),\\
\mr{I}_{2} & \leq \int_{\sqrt{\alpha}}^{\infty} e^{- 4\lambda\alpha\left(\hat{k}-\sqrt{\alpha}\right)^{2}}  \mr{d} \hat{k}
 =\frac{1}{4}\sqrt{\frac{\pi}{\lambda\alpha}}.
\end{aligned}
\end{equation*}
Hence, a combination of the estimates for $\mr{I}_{1}$ and $\mr{I}_{2}$ yields
\begin{equation*}
\begin{aligned}
\mr{I}
& \leq\frac{1}{2}\sqrt{\frac{\pi}{\lambda\alpha}}\operatorname{erf}(\sqrt{\lambda} \alpha)  +  \min \left(\lambda^{-1/4}\Gamma\left(5/4\right),\frac{1}{4}\sqrt{\dfrac{\pi}{\lambda\alpha}}\right)
=:  g(\alpha). \\
\end{aligned}
\end{equation*}
We claim that $g(\alpha)$ attains its maximum at $\alpha_{0} = \pi \lambda^{-1/2} \Gamma(5/4)^{-2}/16$ because for $\alpha < \alpha_{0}$,
\begin{equation*}
\begin{aligned}
g^{\prime}(\alpha) & =\frac{1}{\sqrt{\alpha}}\left(e^{-\lambda \alpha^{2}}-\frac{1}{2 \sqrt{\lambda} \alpha} \int_{0}^{\sqrt{\lambda} \alpha} e^{-x^{2}} \mr{d} x\right)
\geq \frac{1}{\sqrt{\alpha}}\left(e^{-\lambda \alpha_{0}^{2}}-\frac{1}{2}\right)>0 , \\
\end{aligned}
\end{equation*}
and for $\alpha > \alpha_{0}$,
\begin{equation*}
\begin{aligned}
g^{\prime}(\alpha) & =\frac{1}{\sqrt{\alpha}}\left(e^{-\lambda \alpha^{2}}-\frac{1}{2 \sqrt{\lambda} \alpha} \int_{0}^{\sqrt{\lambda} \alpha} e^{-x^{2}} \mr{d} x-\frac{\sqrt{\pi}}{8 \sqrt{\lambda} \alpha}\right)  \\
& \leq \frac{1}{2 \alpha\sqrt{\lambda\alpha}}\left(\sqrt{\lambda} \alpha e^{-\lambda \alpha^{2}}-\frac{\sqrt{\pi}}{4}\right) < 0, \\
\end{aligned}
\end{equation*}
where the last inequality follows from $x e^{-x^{2}}\le (2e)^{-1/2}$ for $x\ge0$. Hence, %
\begin{equation*}
\begin{aligned}
\max_{0\leq \alpha \leq 1} \mr{I}(\alpha) \leq g(\alpha_{0}) %
& = \lambda^{-1/4} \Gamma\left(5/4\right) \left[2\operatorname{erf}\left(\pi \Gamma\left(5/4\right)^{-2}/16\right) + 1 \right].
\end{aligned}
\end{equation*}
Substituting the above inequality into (\ref{ineq: the first ineq to control the integral}), we complete the proof.
\end{proof}

\begin{figure}[ht]
    \centering
    \includegraphics[width=0.36\textwidth]{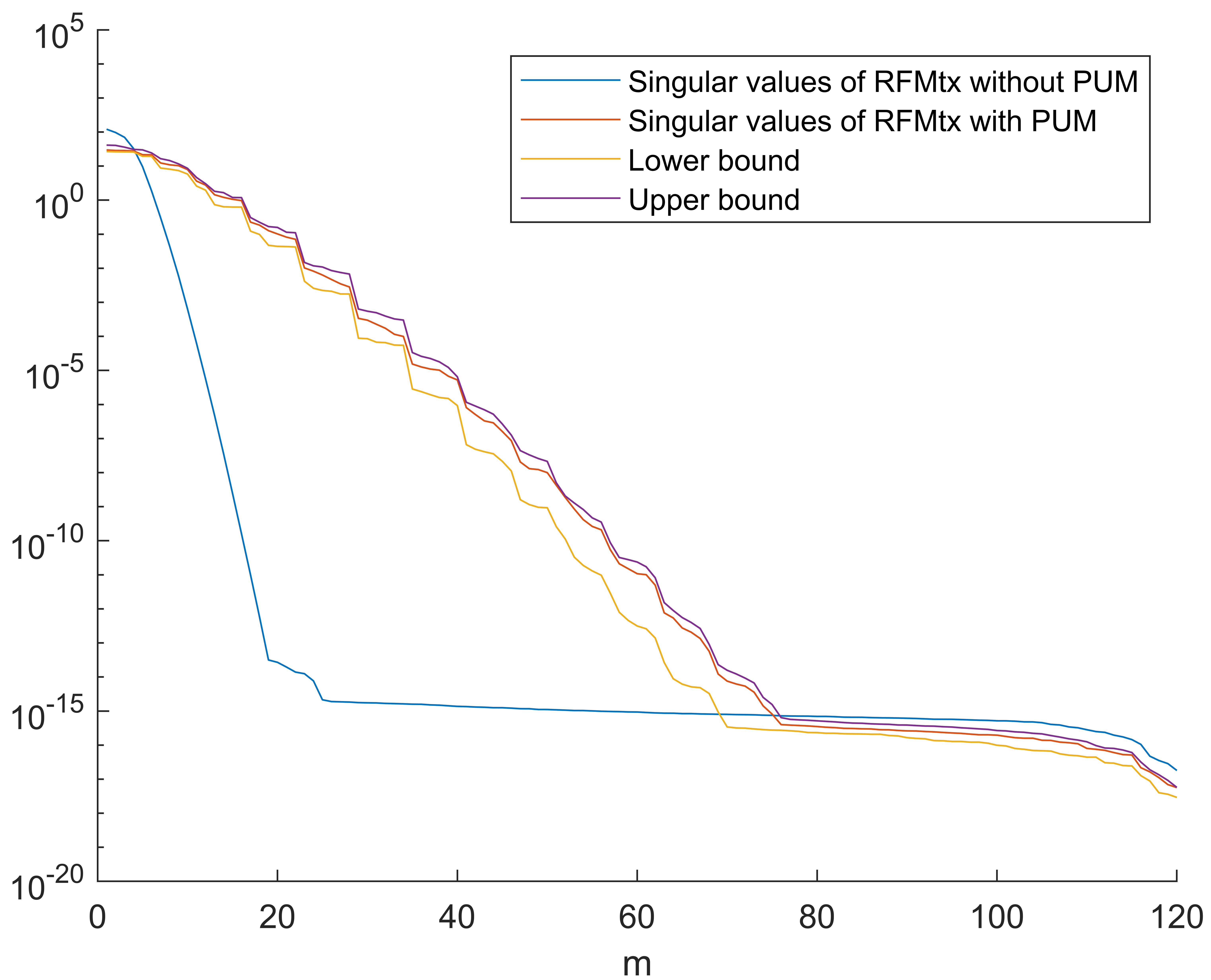}
    \includegraphics[width=0.4\textwidth]{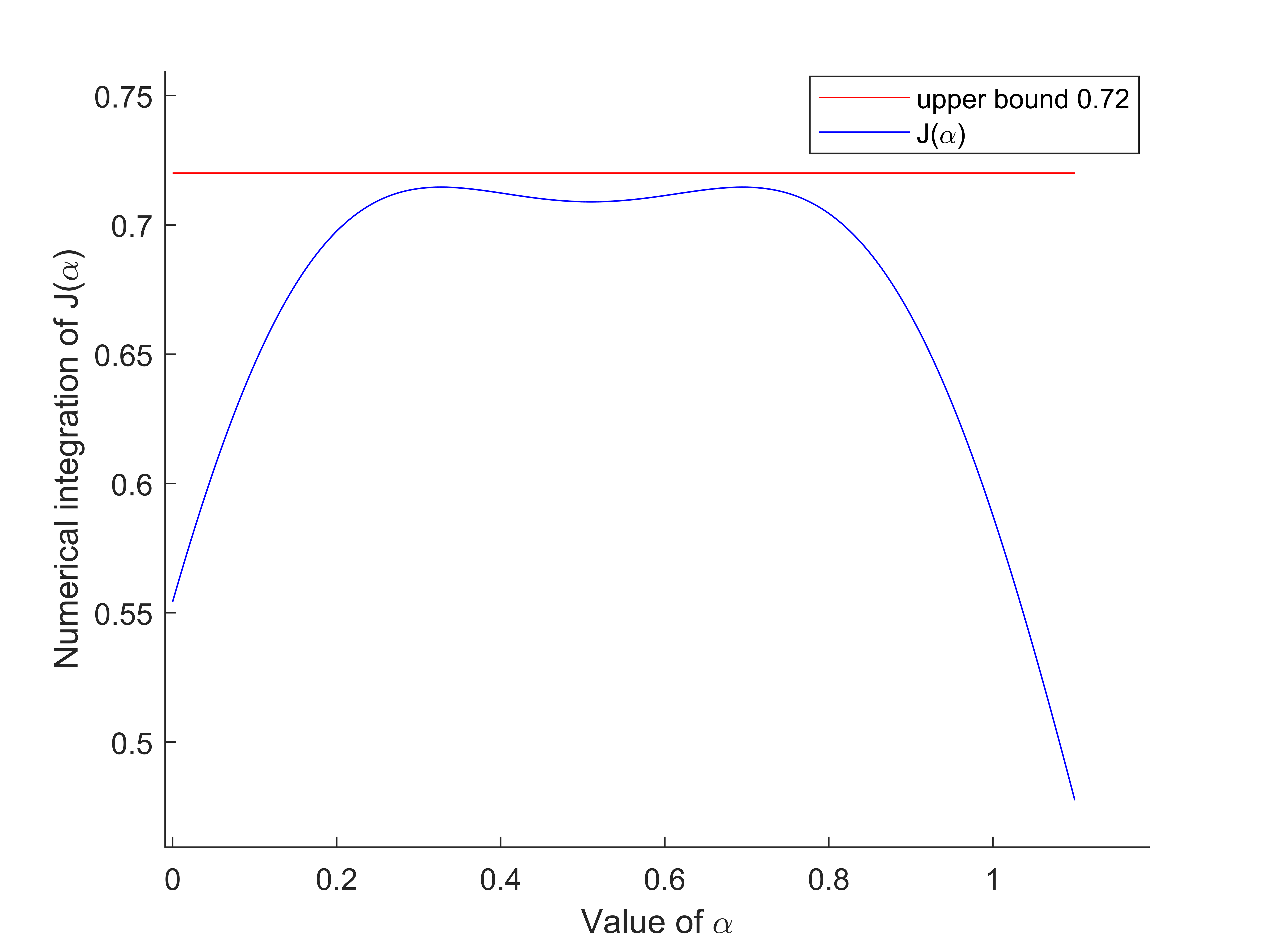}
    \caption{\textbf{Left:} The singular values of RFMtx with and without PUM, and the upper and lower bounds given by (\ref{ineq: basic lower and upper bounds for RFMtx with PUM}). $\widetilde{L}$ is taken as the identity operator. \textbf{Right:} Numerical integration of $\mr{J}(\alpha)$ and its upper bound $0.72$. }  \label{Figure: integralbound}
\end{figure}

The estimate in the above lemma is far from optimal. With the choice of $\delta = 1.1$ and $\lambda = 7.15$, by $\max_{\alpha \in \mb{R}} J(\alpha)  \leq \max_{0\leq \alpha \leq \delta} J(\alpha)$, we calculate the integral in $\mr{J}(\alpha)$ numerically for $0\leq \alpha\leq \delta$ (see Fig.~\ref{Figure: integralbound}). It shows that $\mr{J}(\alpha)$ has a clear upper bound 0.72, which implies
$\mb{P}\left\{\rho \leq  1/11 \right\} \leq 0.72^{N} $ as in the proof of Theorem \ref{Theorem: decay rate of random feature matrices and their exponential ill conditionality}.

\section*{Acknowledgments}
The authors thank the anonymous reviewers for their valuable suggestions. This work is supported by the National Natural Science Foundation of China under Grant No. 12371438. The AI-driven experiments, simulations, and model training were conducted on the Robotics AI Scientist Platform of the Chinese Academy of Sciences.

\end{document}